\documentclass[preprint,12pt,sort&compress]{elsarticle}

\usepackage{geometry}
\usepackage{lineno}
\modulolinenumbers[5]

\usepackage[utf8x]{inputenc}
\usepackage[english]{babel}
\usepackage[T1]{fontenc}

\usepackage{amsmath}
\usepackage{amssymb}
\usepackage{graphicx}
\graphicspath{{Figures/}}
\usepackage{subfig}

\usepackage[pdfborder={0 0 0},colorlinks,allcolors=blue]{hyperref}
\usepackage{xcolor}

\usepackage{pgfplots}
\usetikzlibrary{arrows}
\usetikzlibrary{calc}
\pgfplotsset{compat=1.16}

\usepackage{amsthm}
\theoremstyle{definition}

\newtheorem{remark}{Remark}[section]

\usepackage{tcolorbox}

\usepackage{siunitx} 

\makeatletter
\newcommand*\bigcdot{\mathpalette\bigcdot@{.8}}
\newcommand*\bigcdot@[2]{\mathbin{\vcenter{\hbox{\scalebox{#2}{$\m@th#1\bullet$}}}}}
\makeatother

\usepackage{bm} 

\usepackage{cleveref}
\hypersetup{colorlinks=true,allcolors=blue,pdfauthor=author}
\crefname{equation}{Equation}{Equations}
\crefname{figure}{Figure}{Figures}
\crefname{section}{Section}{Sections}
\crefname{remark}{Remark}{Remarks}
\crefname{appendix}{}{}

\usepackage{xparse}

\NewDocumentCommand \an{ m }{%
    \langle {#1} \rangle%
}

\newcommand*{\medcap}{\mathbin{\scalebox{1.5}{\ensuremath{\cap}}}}%

\newcommand\notype[1]{\unskip}

\title{IGA-PD Penalty-Based Coupling for Immersed Air-Blast Fluid--Structure Interaction: A Simple and Effective Solution for Fracture and Fragmentation}

\begin{document}

\begin{frontmatter}



\author[brown]{Masoud Behzadinasab\corref{mycorrespondingauthor}}
\ead{masoud\_behzadinasab@brown.edu}
\author[psu]{Michael Hillman}
\author[brown]{Yuri Bazilevs\corref{mycorrespondingauthor}}
\ead{yuri\_bazilevs@brown.edu}

\address[brown]{School of Engineering, \\ Brown University, Providence, RI 02912, USA}
\address[psu]{Department of Civil and Environmental Engineering, \\ The Pennsylvania State University, University Park, PA 16802, USA}

\cortext[mycorrespondingauthor]{Corresponding authors.}

\begin{abstract}
We present a novel formulation for the immersed coupling of Isogeometric Analysis (IGA) and Peridynamics (PD) for the simulation of fluid–-structure interaction (FSI). We focus on air-blast FSI and address the computational challenges of immersed FSI methods in the simulation of fracture and fragmentation by developing a weakly volume-coupled FSI formulation by means of a simple penalty approach. We show the mathematical formulation and present several numerical examples of inelastic ductile and brittle solids that clearly demonstrate the power and robustness of the proposed methodology.
\end{abstract}

\begin{keyword}
Air blast \sep Fluid--structure interaction \sep Immersed methods \sep Penalty coupling \sep Isogeometric analysis \sep Meshfree methods \sep Peridynamics \sep Fracture mechanics
\end{keyword}

\end{frontmatter}


\section{Introduction}
\label{sec:intro}

The present work is a continuation of the efforts first reported in~\cite{bazilevs2017new1,bazilevs2017new2}, where the authors formulated an immersed fluid--structure interaction (FSI) framework for air blast. The framework made use of the fluid and structural mechanics governing equations posed in the weak form and two meshes, background and foreground, to discretize the resulting coupled formulation. The coupled FSI solution was approximated using the basis functions of a fixed background mesh. The foreground mesh was mainly used to carry out the solid-domain quadrature, track the solid current configuration, and store the solid constitutive law history variables. It was demonstrated that the use of high-order accurate and smooth background dicsretizations, such as those in NURBS-based Isogeometric Analysis (IGA)~\cite{hughes2005isogeometric,cottrell2009isogeometric}, resulted in much higher solution quality, especially for the solid, due to the higher-order accuracy of the strain-rate approximation and the circumvention of cell-crossing instabilities occurring in the traditional Material-Point Methods (MPMs)~\cite{steffen2008analysis,moutsanidis2020iga}. The proposed approach resulted in an a priori strongly coupled formulation and eliminated mesh distortion issues associated with Lagrangian or other moving-mesh techniques for this class of problems.

While the resulting FSI framework and its IGA-based implementation present a very promising approach and a clear advance over the existing methods, it was pointed out in~\cite{behzadinasab2021coupling} that the proposed methodology could benefit from further improvements. 

Discretizing the solid on the background domain precludes the direct use of the existing, well established solid and structural mechanics solvers in the proposed FSI framework. This lack of modularity was recently addressed in~\cite{behzadinasab2021coupling}, where the authors developed a strongly coupled FSI formulation that introduced linear constraints between the background and foreground discrete function spaces. The kinematic compatibility was enforced through the constraint on the solution spaces, while that same constraint on the test functions gave a consistent definition of the discrete structural force vector on the nodes (or control points in the case of IGA) of the background grid. The resulting formulation was successfully demonstrated on the coupling of IGA with a state-based Peridynamic (PD) solid~\cite{silling2000reformulation,silling2007peridynamic,behzadinasab2020semi}. 

Constraining the foreground solid to the background fluid kinematics as in ~\cite{behzadinasab2021coupling} gives the desired modularity together with strong coupling, however, the modeling of fracture and fragmentation in the immersed FSI simulations remains a challenge. While the foreground discretization such as PD can easily support discontinuous kinematic fields by locally breaking bonds between material points~\cite{bobaru2015cracks,jafarzadeh2018peridynamic,behzadinasab2018peridynamics,kramer2019third,chen2019peridynamic,behzadinasab2019third,gao2020peridynamics,butt2021peridynamic,behzadinasab2021shell}, the smooth background discretization of IGA~\cite{bazilevs2006isogeometric,bazilevs2008isogeometric,auricchio2012simple,hsu2014fluid,kamensky2015immersogeometric,casquero2015a,kamensky2018hyperbolic,moutsanidis2018hyperbolic,zhu2020immersogeometric} is not designed to excel in approximating discontinuous kinematics. Thus, constraining the foreground solution to its background counterpart results in an overly smooth foreground solution and, when coupled with continuum-damage (or phase-field~\cite{kamensky2018hyperbolic,moutsanidis2018hyperbolic}) approaches to model fracture and fragmentation, results in the size of damage zones that scales with that of the background mesh. As a result, unless the background mesh is sufficiently fine, the damage bands appear to be artificially thick and often predict non-physical behavior, especially in the brittle fracture regime. The size of the damage zones may be reduced by refining the background mesh, however, this leads to a significant increase in the computational costs. 

In order to directly address the issue of modeling fracture and fragmentation, an alternative approach is needed where the coupling of background and foreground solutions is weakened, especially in the presence of damage growth and fracture propagation. To develop such a coupling strategy, we took inspiration from the recent work~\cite{wang2021consistent} where an elaborate volume coupling methodology using a Nitsche technique~\cite{nitsche1971variationsprinzip} was formulated in the context of composite materials with a complex microstructure. In the present effort we develop a relatively simple weakly coupled formulation where we only retain the penalty terms of the volumetric Nitsche approach. The main goal here is to assess whether the concept of weak coupling, in our case of IGA and PD, is effective in addressing the challenges involved in the modeling of fracture and fragmentation in the immersed FSI. Although in the present work PD is used for the modeling of the solid on the foreground domain, other meshfree methods like RKPM~\cite{liu1995reproducing,chen1996reproducing,hillman2020generalized} may be naturally employed for that purpose.

The paper is outlined as follows. In \cref{sec:weakFlSol}, we summarize the governing equations of fluid and structural mechanics at the continuum level. In \cref{sec:discFSI}, we revisit the strongly-coupled FSI formulation from~\cite{behzadinasab2021coupling} and develop a new weakly coupled approach by means of a carefully designed volumetric penalty operator. In \cref{sec:implementation}, several 2D numerical examples are presented to compare the performance of the strongly and weakly coupled immersed FSI formulations, to study the effect of the penalty parameter choices on the solution quality, and to demonstrate the effectiveness of weakly coupled approaches in dealing with fracture and fragmentation in brittle and ductile solids subjected to blast loading. In \cref{sec:conclusions}, we make concluding remarks and outlines future research directions.

\section{Fluid and Structural Mechanics Governing Equations}
\label{sec:weakFlSol}

Let $\Omega$ denote the combined fluid and solid domain, and let $\Omega^f$ and $\Omega^s$ denote the individual fluid and solid subdomains in the spatial configuration, such that $\Omega^f\bigcup\Omega^s = \Omega$ and $\Omega^f\bigcap\Omega^s = \emptyset$. Both the fluid and solid problems are stated in terms of weak or variational forms, which are provided in what follows. 

The fluid mechanics problem is governed by the Navier--Stokes equations of compressible flows. The following semilinear forms and linear functionals comprise the weak form of the fluid problem
\begin{align}
\label{MfCont}
M^f_{\omega} (\mathbf{W},\mathbf{Y}) = \int_{\omega} \mathbf{W} \cdot \mathbf{A}^f_0 \mathbf{Y}_{,t}~d\omega, 
\end{align}
\begin{align}
\label{BfCont}
B^f_{\omega} (\mathbf{W},\mathbf{Y}) = \int_{\omega} \mathbf{W} \cdot \mathbf{A}^f_i  \mathbf{Y}_{,i}~d\omega - \int_{\omega} \mathbf{W}_{,i} \cdot (\mathbf{F}^{p}_{i} - \mathbf{F}^{d}_{i})~d\omega,
\end{align}
\begin{align}
\label{FfCont}
F^f_{\omega} (\mathbf{W}) = \int_{\omega} \mathbf{W} \cdot \mathbf{S}^f~d\omega + \int_{\Gamma^f_H} \mathbf{W} \cdot \mathbf{H}^f~d\Gamma.
\end{align}
Here, $\mathbf{Y}$ denotes a set of pressure-primitive variables~\cite{hauke1998comparative,hauke2001simple},
\begin{equation}
\mathbf{Y} = 
\begin{bmatrix}
p \\
\mathbf{v} \\
T
\end{bmatrix},
\label{eqn:Y}
\end{equation}
where $p$ is the pressure, $\mathbf{v}$ is the material-particle velocity vector, and $T$ is the temperature. $\mathbf{Y}$ and $\mathbf{W}$, the vector-valued trial and test functions, respectively, are the members of $\mathcal{S}$ and $\mathcal{V}$, the corresponding trial and test function spaces, respectively, \textit{defined on all of} $\Omega$. $\Gamma^f_H$ is the subset of the fluid-domain boundary where natural boundary conditions are imposed, and $\mathbf{H}^f$ contains the prescribed values of the natural boundary conditions.  $\mathbf{A}_0$ and $\mathbf{A}_i$ are the so-called Euler jacobian matrices, $\mathbf{F}^{p}_{i}$ and $\mathbf{F}^{d}_{i}$ are the pressure and viscous/thermal fluxes, respectively, and $\mathbf{S}$ is the volume source. (See~\cite{bazilevs2017new1} and references therein for the details.) The subscript $\omega$ on the semilinear forms and linear functionals denotes the domain of integration, comma denotes partial differentiation with respect to the spatial coordinates, and $i=1,\dots,d$, where $d=2,3$ is the space dimension. The compressible-flow equations are complemented with the ideal gas law as the equation of state.

With the above definitions, following~\cite{bazilevs2017new1}, the continuous weak form the compressible flow problem may be written as: Find $\mathbf{Y} \in \mathcal{S}$, such that $\forall \mathbf{W} \in \mathcal{V}$,
\begin{equation}
\label{EqFluidCont}
 M^f_{\Omega^f}(\mathbf{W},\mathbf{Y}) + B^f_{\Omega^f}(\mathbf{W},\mathbf{Y}) - F^f_{\Omega^f}(\mathbf{W}) = 0.
\end{equation}

The structural mechanics is modeled as an isothermal large-deformation inelastic solid in the framework of correspondence-based PD. Let $\tilde{\mathcal{S}}$ and $\tilde{\mathcal{V}}$ denote the trial and test function spaces for the PD formulation defined on $\Omega^s$. Because the mass balance is satisfied in the Lagrangian form and no thermal coupling is assumed in the solid, the pressure and temperature slots of the vector-valued trial ($\tilde{\mathbf{Y}}$) and test ($\tilde{\mathbf{W}}$) functions are set to zero, while the slots $2:2+(d-1)$ are occupied by the solid velocity trial functions ($\tilde{\mathbf{v}}$) and the momentum-equation test functions ($\tilde{\mathbf{w}}$), respectively. That is,
\begin{equation}
\tilde{\mathbf{Y}} = 
\begin{bmatrix}
0 \\
\tilde{\mathbf{v}} \\
0
\end{bmatrix},
\label{eqn:Ytilde}
\end{equation}
and
\begin{equation}
\tilde{\mathbf{W}} = 
\begin{bmatrix}
0 \\
\tilde{\mathbf{w}} \\
0
\end{bmatrix}.
\label{eqn:Ytilde}
\end{equation}
Following~\cite{behzadinasab2021coupling} we define the semilinear forms and linear functionals that comprise the weak form of the solid problem:
\begin{equation}
\label{EqWFPDM}
\mathcal{M}^s_{\Omega^s}(\tilde{\mathbf{W}},\tilde{\mathbf{Y}}) = \int_{\Omega^s}  \tilde{\mathbf{w}} \cdot \rho^s \, \dot{\tilde{\mathbf{v}}} \, {\rm d}\Omega^s,
\end{equation}
\begin{equation}
\label{EqWFPDB}
\mathcal{B}^s_{\Omega^s}(\tilde{\mathbf{W}},\tilde{\mathbf{Y}}) = \int_{\Omega^s}  \tilde{\mathbf{w}} \cdot \int_{\mathcal{H}} \left( \underline{\mathbf{T}} - \underline{\mathbf{T}'} \right) \, {\rm d}\mathcal{H} \, {\rm d}\Omega^s,
\end{equation}
and
\begin{equation}
\label{EqWFPDF}
\mathcal{F}^s_{\Omega^s}(\tilde{\mathbf{W}}) = \int_{\Omega^s}  \tilde{\mathbf{w}} \cdot \mathbf{s} \, {\rm d}\Omega^s.
\end{equation}
Here, $\rho^s$ is the solid mass density in the current configuration, $\mathbf{s}$ is the volume source term, and the overdot symbol denotes the material time derivative. $\mathcal{H}(\mathbf{x})$ is the PD family set (also known as the PD horizon) of the spatial location $\mathbf{x}$ defined as
\begin{equation} 
\mathcal{H}(\mathbf{x}) = \left\{ \mathbf{x}' \ | \ \mathbf{x}' \in \mathcal{H}(\mathbf{x}) \medcap \Omega^s, \, 0 < | \mathbf{x}' - \mathbf{x} | \leq \delta \right\} ,
\label{eqn:family}
\end{equation}
where $\delta$ is the horizon size, ${\an{\mathbf{x}-\mathbf{x}'}}$ denotes a PD bond between $\mathbf{x}$ and $\mathbf{x}'$, $\underline{\mathbf{T}} = \mathbf{T} \an{\mathbf{x}-\mathbf{x}'}$ is the so-called PD force state (with units of force per unit volume squared), and $\underline{\mathbf{T}'} = \mathbf{T} \an{\mathbf{x}'-\mathbf{x}}$. Here, the bond-associated fields are marked using underscores.

In the correspondence PD framework, the spatial gradients are computed using the integral form. Using the computed velocity gradient, and its symmetrized version, a classical constitutive law may be employed to evaluate the Cauchy stress tensor at the bond level. The relationship between the force state $\underline{\mathbf{T}}$ and the Cauchy stress at the bond level depends on the details of the velocity-gradient definition and evaluation. In the present work, we adopt the formulation detailed in~\cite[Appendix A]{behzadinasab2021coupling}, and the reader is encouraged to consult this reference as well as~\cite{behzadinasab2020semi,behzadinasab2020peridynamic} for further information.

With the above definitions, the continuous, weak correspondence-based PD formulation of the solid becomes: Find $\tilde{\mathbf{Y}} \in \tilde{\mathcal{S}}$, such that $\forall \tilde{\mathbf{W}} \in \tilde{\mathcal{V}}$,
\begin{equation}
\label{EqSolidCont}
 \mathcal{M}^s_{\Omega^s}(\tilde{\mathbf{W}},\tilde{\mathbf{Y}}) + \mathcal{B}^s_{\Omega^s}(\tilde{\mathbf{W}},\tilde{\mathbf{Y}}) - \mathcal{F}^s_{\Omega^s}(\tilde{\mathbf{W}}) = 0.
\end{equation}

\section{Coupled FSI Formulations in a Discrete Form}
\label{sec:discFSI}

In this section we present two discrete forms of the volume-coupled FSI problem. The first approach is presented recently in~\cite{behzadinasab2021coupling}, where the FSI coupling is carried out by explicitly constraining the fluid and solid velocity degrees-of-freedom (DOFs). We refer to this approach as \textit{strong coupling}. As an alternative, we also explore a \textit{weak coupling} approach based on a volumetric penalty formulation. 

In both cases, we first define $\mathcal{S}^h \subset \mathcal{S}$ and $\mathcal{V}^h \subset \mathcal{V}$, the background-domain finite-dimensional trial- and test-function spaces, respectively. The discrete trial and test functions may be expressed as
\begin{equation}
\mathbf{Y}^h(\mathbf{x}) = \sum_{B=1}^{\mathcal{N}_{cp}} \mathbf{Y}_B \, N_B(\mathbf{x})
\label{eqn:backgroundtri}
\end{equation}
and
\begin{equation}
\mathbf{W}^h(\mathbf{x}) = \sum_{A=1}^{\mathcal{N}_{cp}} \mathbf{W}_A \, N_A(\mathbf{x}),
\label{eqn:backgroundtest}
\end{equation}
where $\mathbf{Y}_B$ and $\mathbf{W}_A$ are the vector-valued control-point DOFs and weights, respectively, $N(\mathbf{x})$'s in this work are assumed to be the B-Spline basis functions defined everywhere in $\Omega$, and $\mathcal{N}_{cp}$ is the dimension of the B-Spline space. Note that all the components of the trial and test-function vectors are approximated using the same basis functions. We also define the finite-dimensional trial and test function spaces $\tilde{\mathcal{S}}^h$ and $\tilde{\mathcal{V}}^h$, respectively, for the PD solid. The discrete trial and test functions may be expressed as
\begin{equation}
\tilde{\mathbf{Y}}^h(\mathbf{x}) = \sum_{Q=1}^{\mathcal{N}_{mp}} \tilde{\mathbf{Y}}_Q \, \chi_Q(\mathbf{x})
\label{eqn:foregroundtri}
\end{equation}
and
\begin{equation}
\tilde{\mathbf{W}}^h(\mathbf{x}) = \sum_{P=1}^{\mathcal{N}_{mp}} \tilde{\mathbf{W}}_P \, \chi_P(\mathbf{x}),
\label{eqn:foregroundtest}
\end{equation}
where the solid domain is represented using a finite number $\mathcal{N}_{mp}$ of material points or PD nodes, $\tilde{\mathbf{Y}}_Q$ and $\tilde{\mathbf{W}}_P$ are the nodal DOFs and weights, respectively, and $\chi_P(\mathbf{x})$ is a characteristic function of a PD node $P$ that attains unity at $\mathbf{x}_P \in \Omega^s$, the spatial location of the PD node, is zero at all other PD nodes, and satisfies
\begin{equation}
\int_{\Omega^s} \chi_P(\mathbf{x}) \, {\rm d}\Omega^s = V_P,
\end{equation}
where $V_P$ is the local volume of the PD node.

For the strong coupling approach we define a mapping $\Pi$ between the discrete background and foreground test and trial functions using nodal interpolation. That is, we constrain the PD nodal DOFs and weights to the background IGA test and trial functions, respectively, as 
\begin{equation}
\tilde{\mathbf{Y}}^h (\mathbf{x}) = \Pi \mathbf{Y}^h(\mathbf{x}) = \sum_{Q=1}^{\mathcal{N}_{mp}} \left(~\sum_{B=1}^{\mathcal{N}_{cp}} \mathbf{Y}_B \, N_B(\mathbf{x}_Q)~\right) \, \chi_Q(\mathbf{x})
\label{eq:trialprojfun}
\end{equation}
and
\begin{equation}
\tilde{\mathbf{W}}^h(\mathbf{x}) = \Pi \mathbf{W}^h(\mathbf{x}) = \sum_{P=1}^{\mathcal{N}_{mp}} \left(~\sum_{A=1}^{\mathcal{N}_{cp}} \mathbf{W}_A \, N_A(\mathbf{x}_P)~\right) \, \chi_P(\mathbf{x}) .
\label{eq:testprojfun}
\end{equation}

With these definitions, the spatially discretized, immersed, strongly coupled FSI formulation may now be stated solely in terms of the background domain unknowns as: Find $\mathbf{Y}^h \in \mathcal{S}^h$, such that $\forall \mathbf{W}^h \in \mathcal{V}^h$, 
\begin{align}
\label{eqn:FSIDisc}
M^f_{\Omega}(\mathbf{W}^h,\mathbf{Y}^h) & + B^f_{\Omega}(\mathbf{W}^h,\mathbf{Y}^h) - 
F^f_{\Omega}(\mathbf{W}^h) + B^{st}_{\Omega}(\mathbf{W}^h,\mathbf{Y}^h) + B^{dc}_{\Omega}(\mathbf{W}^h,\mathbf{Y}^h) \nonumber \\
& \textcolor{red}{+} \nonumber \\
\mathcal{M}^s_{\Omega^s}(\tilde{\mathbf{W}}^h,\tilde{\mathbf{Y}}^h) & + \mathcal{B}^s_{\Omega^s}(\tilde{\mathbf{W}}^h,\tilde{\mathbf{Y}}^h) - \mathcal{F}^s_{\Omega^s}(\tilde{\mathbf{W}}^h) \nonumber \\
& \textcolor{red}{-} \nonumber \\
(~M^f_{\Omega^s}(\mathbf{W}^h,\mathbf{Y}^h) & + B^f_{\Omega^s}(\mathbf{W}^h,\mathbf{Y}^h) - 
F^f_{\Omega^s}(\mathbf{W}^h) + B^{st}_{\Omega^s}(\mathbf{W}^h,\mathbf{Y}^h) + B^{dc}_{\Omega^s}(\mathbf{W}^h,\mathbf{Y}^h)~) \nonumber \\
& \textcolor{red}{=} \nonumber \\ 
&0.
\end{align}
Note that, in the space-discrete case, we introduce the SUPG stabilization ($B^{st}_{\Omega}(\mathbf{W}^h,\mathbf{Y}^h)$)~\cite{brooks1982streamline,le1993supg,tezduyar2006stabilization,hughes2010stabilized}  and  discontinuity-capturing ($B^{dc}_{\Omega}(\mathbf{W}^h,\mathbf{Y}^h)$)~\cite{hughes1986new,tezduyar2006computation,rispoli2009computation,rispoli2015particle} operators to the compressible flow formulation to address the convective instability of the Galerkin technique in the regime of convection dominance and to provide additional dissipation in the shock regions. More details on the definition of these operators that are employed in the present formulation of compressible flows may be found in~\cite{hauke1998comparative,xu2017compressible,bazilevs2021gas}. Also note that the integration over the fluid mechanics domain is replaced by the integration over the combined domain \textit{minus} that over the solid domain. This form of the coupled problem statement is convenient for the numerical integration of the semi-discrete forms of the coupled FSI problem (see, e.g.,~\cite{casquero2015a}).

\begin{remark}
Constraining the trial and test functions of the PD domain to that of the background domain results in the following algorithmic approach to the strong FSI coupling. At the beginning of the step, the background fluid solution $\mathbf{Y}^h$ is interpolated to the PD mesh using Equation~(\ref{eq:trialprojfun}), resulting in the field $\tilde{\mathbf{Y}}^h$. At this stage, $\mathbf{Y}^h$ is used to evaluate the discrete residual vector (also often called the nodal force vector) for the background domain, while $\tilde{\mathbf{Y}}^h$ is used to evaluate the discrete residual vector for the foreground domain. The foreground-domain residual vector is then distributed to the background-domain DoFs using a linear transformation induced by Equation~(\ref{eq:testprojfun}). The reader is referred to~\cite{behzadinasab2021coupling} for the details of this transformation. At this stage the background-DoF residuals are added, and the solution increment is computed on the background mesh. The procedure repeats if multiple passes per step are employed.
\end{remark}

\begin{remark}
The resulting strong coupling methodology resembles the classical Immersed Boundary Method~\cite{peskin1972flow} and a more recent Immersed Finite Element Method~\cite{zhang2004immersed}, but without the use of {\em ad hoc} smoothed delta functions to distribute the foreground-domain residual vector to the background DoFs. The foreground-domain residual vector distribution on the background DoFs is defined to be consistent with the test-function constraints given by Equation~(\ref{eq:testprojfun}). This presents a clear benefit of using a variational formulation in the background domain, which is associated with the fluid mechanics part of the problem.
\end{remark}

For the weak coupling approach, the fluid and structural equations are discretized independently and a volumetric penalty term is added to the formulation to penalize the deviation between the fluid and structural velocities. The resulting formulation is now stated in terms of both the background- and foreground-mesh unknowns as: Find $\mathbf{Y}^h \in \mathcal{S}^h$ and $\tilde{\mathbf{Y}}^h \in \tilde{\mathcal{S}}^h$, such that for all $\mathbf{W}^h \in \mathcal{V}^h$ and $\tilde{\mathbf{W}}^h \in \tilde{\mathcal{V}}^h$,
\begin{align}
\label{eqn:FSIDiscPen}
M^f_{\Omega}(\mathbf{W}^h,\mathbf{Y}^h) & + B^f_{\Omega}(\mathbf{W}^h,\mathbf{Y}^h) - 
F^f_{\Omega}(\mathbf{W}^h) + B^{st}_{\Omega}(\mathbf{W}^h,\mathbf{Y}^h) + B^{dc}_{\Omega}(\mathbf{W}^h,\mathbf{Y}^h) \nonumber \\
& \textcolor{red}{+} \nonumber \\
\mathcal{M}^s_{\Omega^s}(\tilde{\mathbf{W}}^h,\tilde{\mathbf{Y}}^h) & + \mathcal{B}^s_{\Omega^s}(\tilde{\mathbf{W}}^h,\tilde{\mathbf{Y}}^h) - \mathcal{F}^s_{\Omega^s}(\tilde{\mathbf{W}}^h) \nonumber \\
& \textcolor{red}{+} \nonumber \\
\int_{\Omega^s}  (\mathbf{w}^h - \tilde{\mathbf{w}}^h) & \cdot C_{\rm pen} \, (\mathbf{v}^h - \tilde{\mathbf{v}}^h) \, {\rm d}\Omega^s \nonumber \\
& \textcolor{red}{=} \nonumber \\ 
&0.
\end{align}
The penalty parameter $C_{\rm pen}$ in Equation~(\ref{eqn:FSIDiscPen}) needs a careful design to ensure a proper coupling between the fluid and structural systems and to not produce an overly stiff method with significant limitations on the stable time step size. It is also important to note that in the regions where the solid and fluid overlap it is not necessary to generate an accurate fluid mechanics solution because, from the standpoint of the fluid problem, this region is completely fictitious. On the other hand, it is imperative that the solid solution is accurate and stable in this region since this is the actual solid domain. For these reasons we choose $C_{\rm pen}$ to scale with the internal work terms of the solid formulation, which, with the aid of scaling arguments, yields the following definition:
\begin{equation}
 C_{\rm pen} = \beta~\frac{E\Delta t}{h^2}.   
\end{equation}
Here $E$ is the local elastic modulus, $h$ is local mesh size, $\Delta t$ is the time step size, and $\beta$ is a dimensionless positive constant. The latter may be chosen just large enough to ensure that the penalty terms do not dominate the stable time step size, a common practice in contact and impact simulations using penalty methods (see~\cite{benson2010isogeometric,benson2011large,alaydin2021updated}). In the case the material damage is modeled, we propose to further modify the penalty parameter as
\begin{equation}
 C_{\rm pen} = \beta~\frac{(1-d)E\Delta t}{h^2},
\end{equation}
where $d$ is a damage variable (or a phase field variable~\cite{kamensky2018hyperbolic,moutsanidis2018hyperbolic}) with $d=1$ corresponding to a complete loss of material stiffness. We will examine the effects of the penalty parameter choice on the resulting coupled FSI solutions presented in the Numerical Examples section.

\begin{remark}
The penalty parameter $C_{\rm pen}$ is present only in the last integral of Equation~(\ref{eqn:FSIDiscPen}). This integral is evaluated using numerical quadrature associated with the nodes of PD mesh, where the damage field is readily available. 
\end{remark}

\begin{remark}
The resulting weak coupling methodology is similar in structure to Immersogeometric FSI (IMGA-FSI)~\cite{hsu2014fluid,kamensky2015immersogeometric,kamensky2017immersogeometric,kamensky2017projection}, which is a new class of immersed FSI formulations that was developed for the coupling of incompressible flow with a Kirchhoff--Love shell~\cite{kiendl2009isogeometric} using a combination of penalty and Lagrange multiplier techniques. The present approach couples an IGA-based compressible flow formulation to a PD solid by means of a volumetric penalty only. A more elaborate Nitsche-like technique may be formulated in the future to make the approach more robust with respect to the selection of the penalty parameter.
\end{remark}

The resulting semi-discrete FSI equations, for both the strongly and weakly coupled formulations, are integrated in time using the lumped-mass explicit Generalized-$\alpha$ predictor-multi-corrector procedure~\cite{chung1993time} adopted for immersed FSI and detailed in~\cite{bazilevs2017new1}.

\section{Numerical Examples}
\label{sec:implementation}

Three 2D numerical simulations are provided here to test and demonstrate the performance of the IGA-PD framework for FSI with a penalty-based coupling approach for blast loading and fragmentation applications. The first example serves as a verification case for the developed formulation in solving problems involving large deformation and plasticity. The next two cases are demonstrative examples that include fracture propagation and fragmentation in brittle and ductile solids. In all the presented computations, $C^1$-continuous quadratic NURBS functions are used for the background solution. RK functions with quadratic consistency, rectangular support, and bond-associative stabilization~\cite{breitzman2018bond,behzadinasab2020semi} are employed in the PD formulation. The PD support size $\delta$ is chosen with respect to the mesh size $h$ and the quadratic order of the method, i.e., $\delta = 2.5 h$~\cite{behzadinasab2021unifiedI}. The fluid is assumed to have properties of air at room temperature, namely, initial density $\rho = 1.0 \, {\rm kg}/{\rm m}^3$ viscosity  $\mu = 1.81 \times 10^{-5} \, {\rm kg}/({\rm m \, s})$, Prandtl number 0.72, and adiabatic index $\gamma=1.4$. In the examples involving damage, we report a \textit{normalized solid mass loss}, which we define as
\begin{equation}
 L_{\Omega^s} = \frac{\int_{\Omega^s} \chi(d) \, \rho^s \, {\rm d}\Omega^s}{\int_{\Omega^s} \rho^s \, {\rm d}\Omega^s} \, , \qquad \chi(d) = \begin{cases}
0 \quad & \text{if} \ d<0.99 \\
1 & \text{otherwise}
\end{cases}	
,
\end{equation}
where $d$ is the scalar damage field.

\subsection{Chamber Detonation}
\label{sec:detonation}

\begin{figure*}[!hbpt]
  \centering
  \subfloat[][]{\includegraphics[width=0.7\textwidth,trim={0cm 0cm 0cm 0cm},clip]{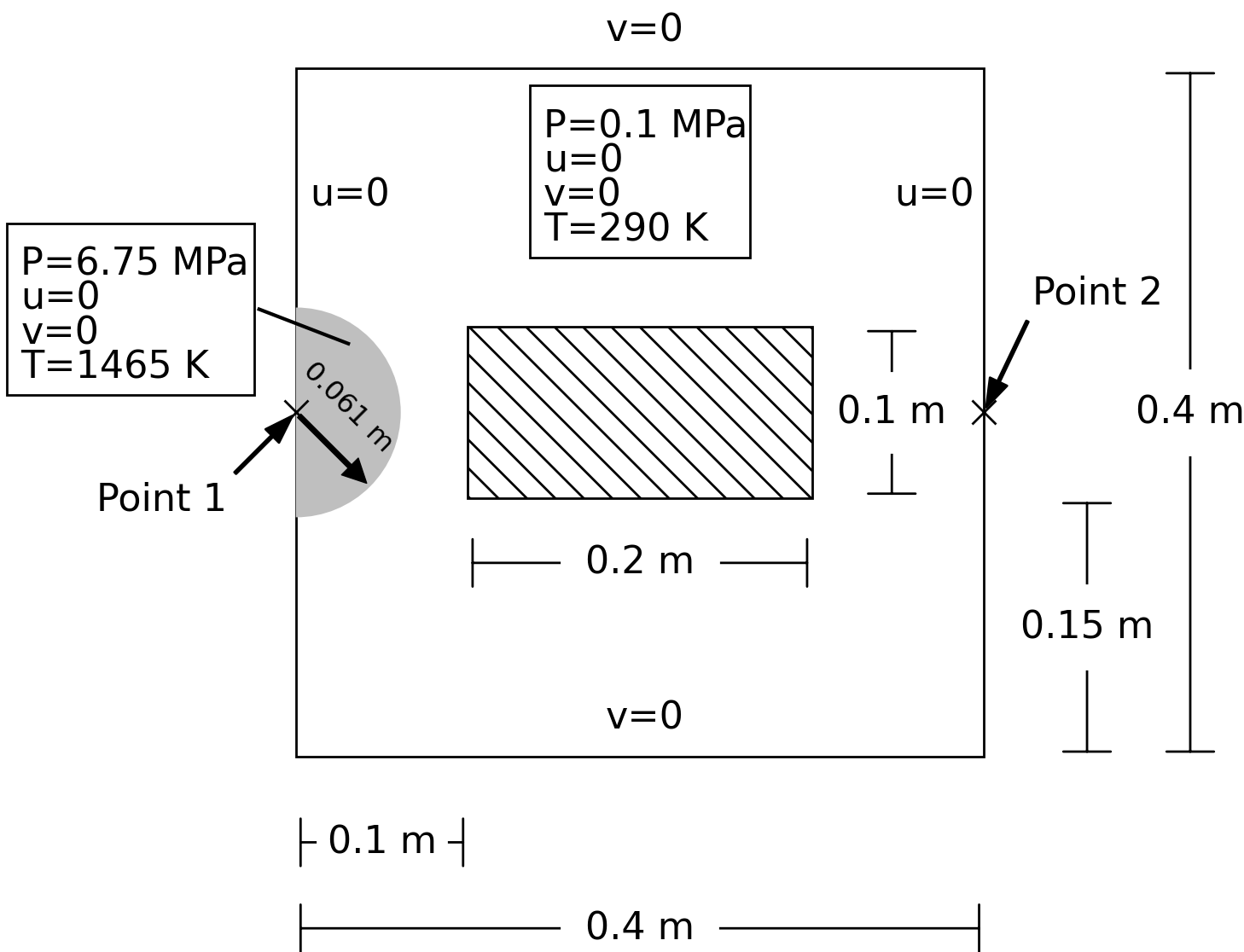}}
  \caption{Chamber detonation. Problem setup and geometry.}
  \label{fig:detonation_setup}
\end{figure*}

In this example, taken from~\cite{bazilevs2017new1}, a steel bar is subjected to a detonation load inside a closed chamber. The problem description is shown in \cref{fig:detonation_setup}, where a $0.2 \, {\rm m} \, \times \, 0.1 \, {\rm m}$ block is located at the center of a closed chamber with dimensions $0.4 \, {\rm m} \, \times \, 0.4 \, {\rm m}$.
The bar thickness is set to $3.5 \, {\rm mm}$. Isotropic linear hardening rule is used for the solid material with Young's modulus ${\rm E}=200 \, {\rm GPa}$, Poisson's ratio $\nu=0.29$, yield stress $\sigma_Y=0.4 \, {\rm GPa}$, hardening modulus $H = 0.1 \, {\rm GPa}$, and initial density $\rho^s=7870 \, {\rm kg}/{\rm m}^3$. Initially, the air in the chamber is at rest with $p = 0.1 \, {\rm MPa}$ and $T = 290 \, {\rm K}$. The detonation is initiated by setting higher-than-ambient values of the pressure and temperature, i.e., $p = 6.75 \, {\rm MPa}$ and $T = 1465 \, {\rm K}$, in a semi-circular area with a radius of $6.1 \, {\rm mm}$, centered on the left wall. No-penetration and free-slip boundary conditions are prescribed at the chamber walls. 

We carry out the weakly coupled simulations using four non-dimensional penalty constants, $\beta = 1/3, 1.0, 3.0, 9.0$, and compare with the strongly coupled case. Three different discretizations are considered for each case: coarse - Fluid: $20 \times 20$ elements; Solid: $30 \times 15$ elements (PD nodes); medium - Fluid: $40 \times 40$ elements; Solid: $60 \times 30$ elements (PD nodes); fine - Fluid: $80 \times 80$ elements; Solid: $120 \times 60$ elements (PD nodes). In the PD case, each foreground element is replaced by a meshfree node at its centroid with an equivalent volume. The time step size used for the coarse, medium, and fine strongly coupled cases is $1 \, {\rm ms}$, $0.5 \, {\rm ms}$, and $0.25 \, {\rm ms}$, respectively. The time step size used for all the weakly coupled cases is taken to be eight times smaller than the corresponding values for the strongly coupled cases. This factor is chosen such that the simulations remain stable for the largest value of $\beta = 9.0$. 

\begin{figure*}[!hbpt]
  \centering
  \subfloat{\includegraphics[width=0.35\textwidth,trim={0cm 0cm 0cm 4cm},clip]{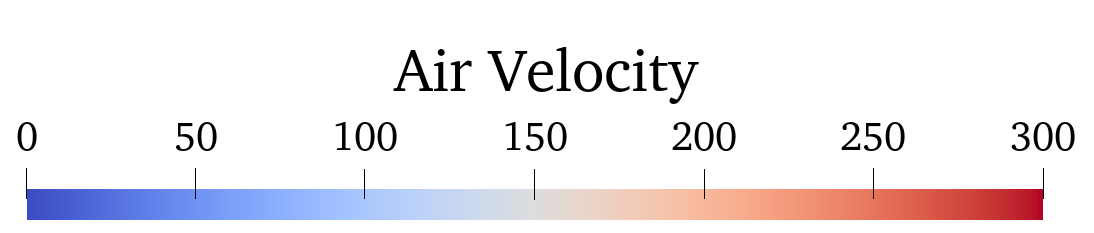}}
  \hspace{2cm}
  \subfloat{\includegraphics[width=0.35\textwidth,trim={0cm 0cm 0cm 4cm},clip]{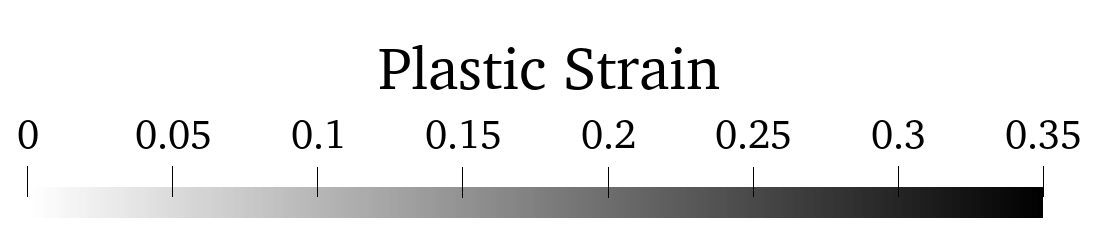}}
  
  \subfloat{\includegraphics[width=0.21\textwidth,trim={0cm 0cm 0cm 0cm},clip]{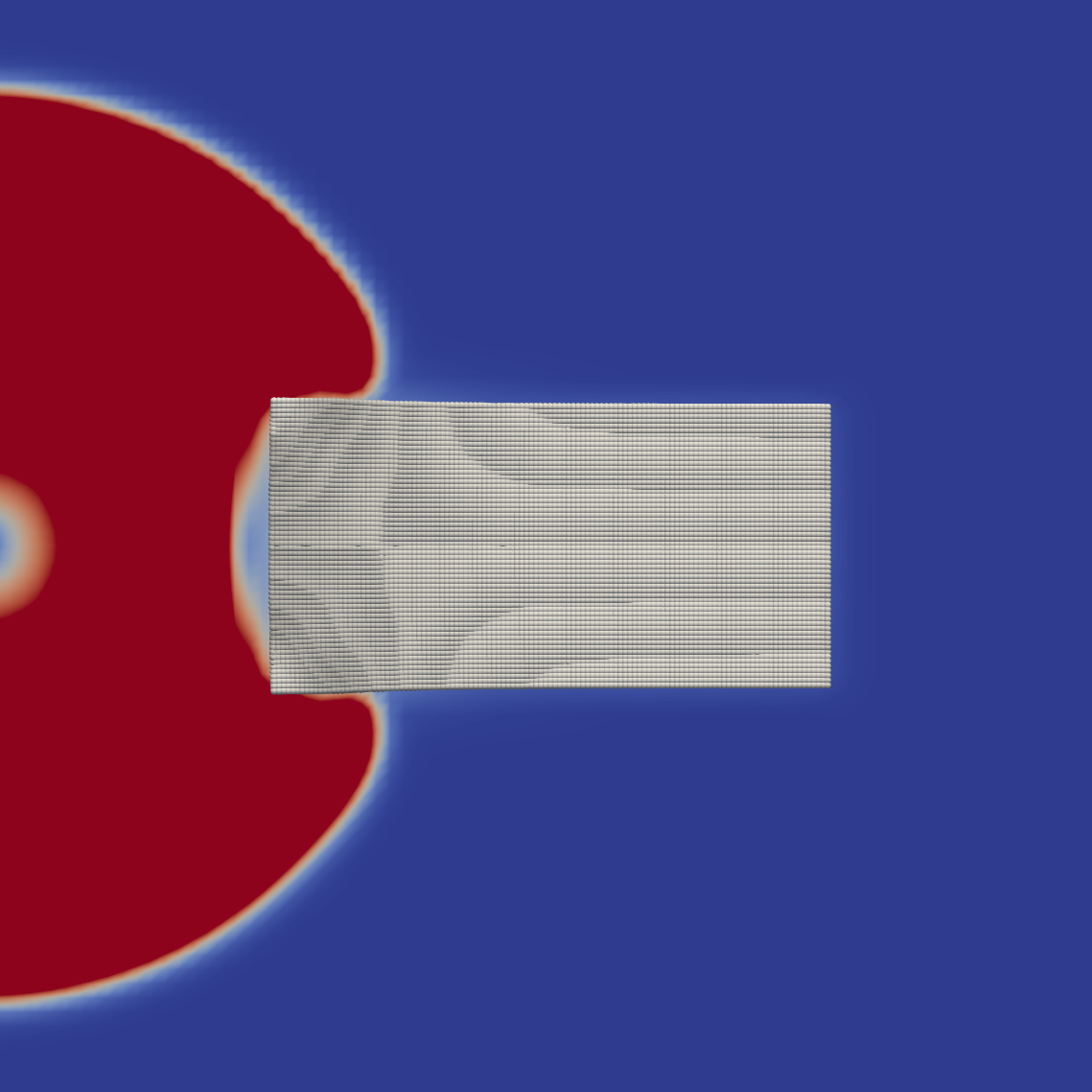}}
  \hspace{7pt}
  \subfloat{\includegraphics[width=0.21\textwidth,trim={0cm 0cm 0cm 0cm},clip]{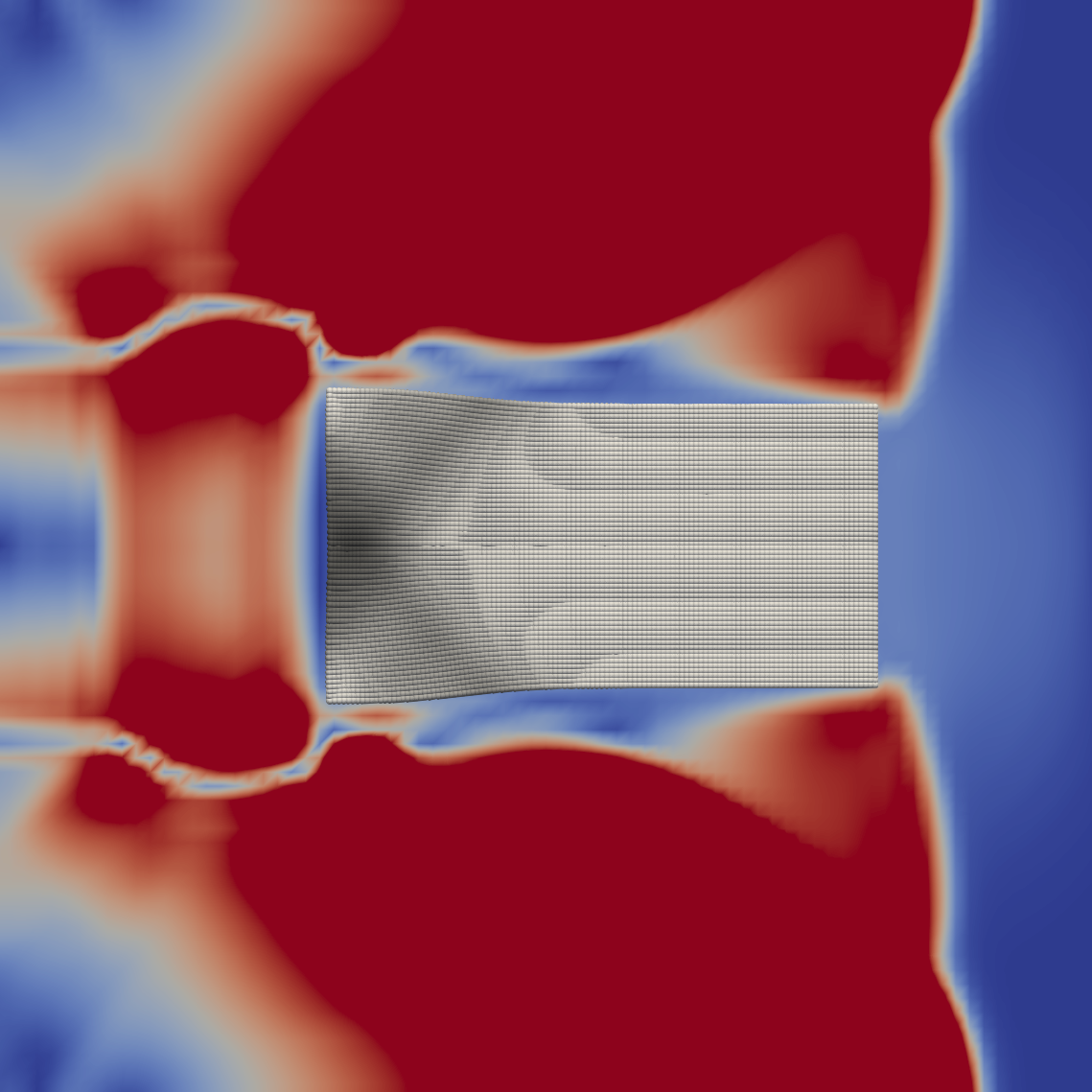}}
  \hspace{7pt}
  \subfloat{\includegraphics[width=0.21\textwidth,trim={0cm 0cm 0cm 0cm},clip]{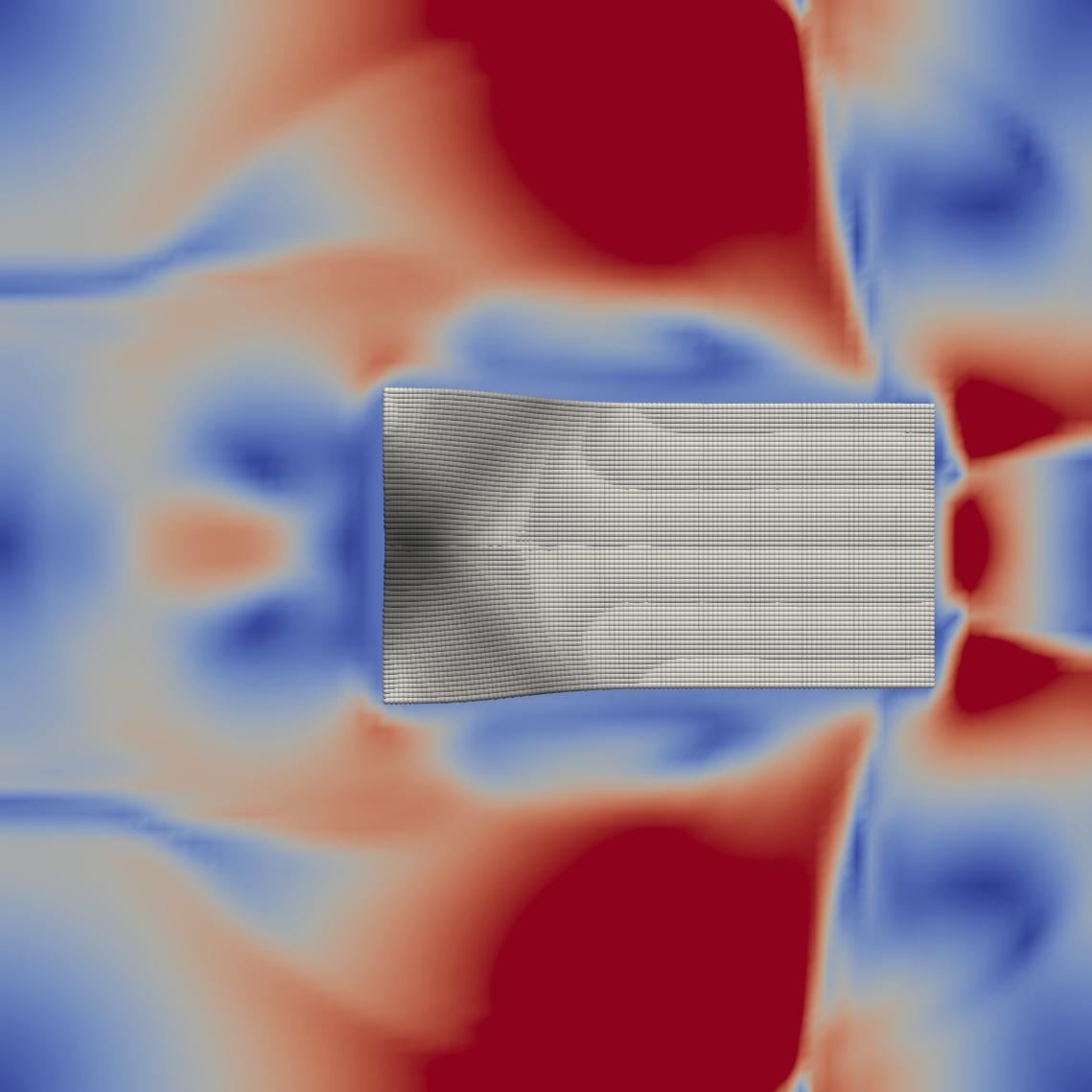}}
  \hspace{7pt}
  \subfloat{\includegraphics[width=0.26\textwidth,trim={0cm 0cm 0cm 0cm},clip]{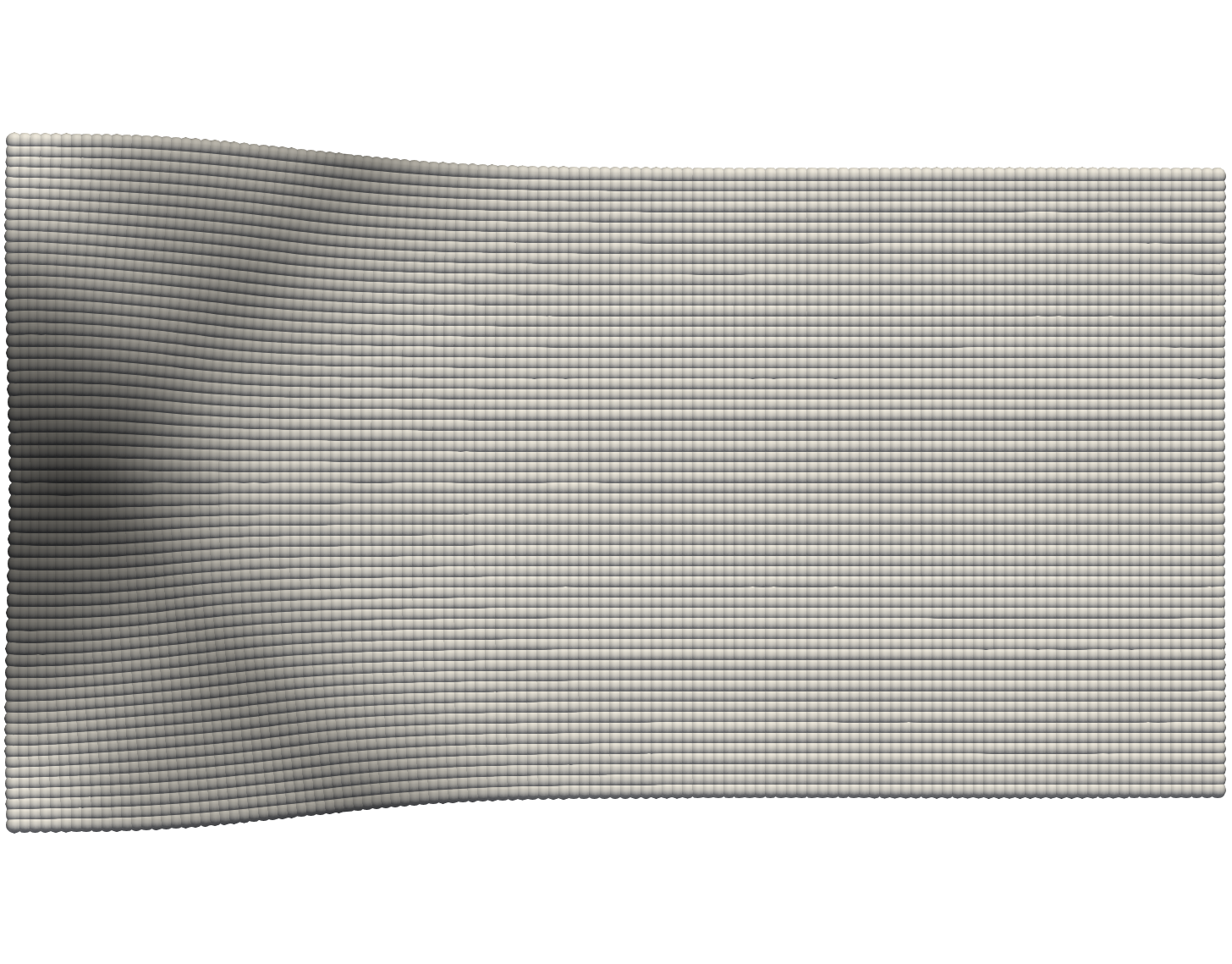}}
  
  \subfloat{\includegraphics[width=0.21\textwidth,trim={0cm 0cm 0cm 0cm},clip]{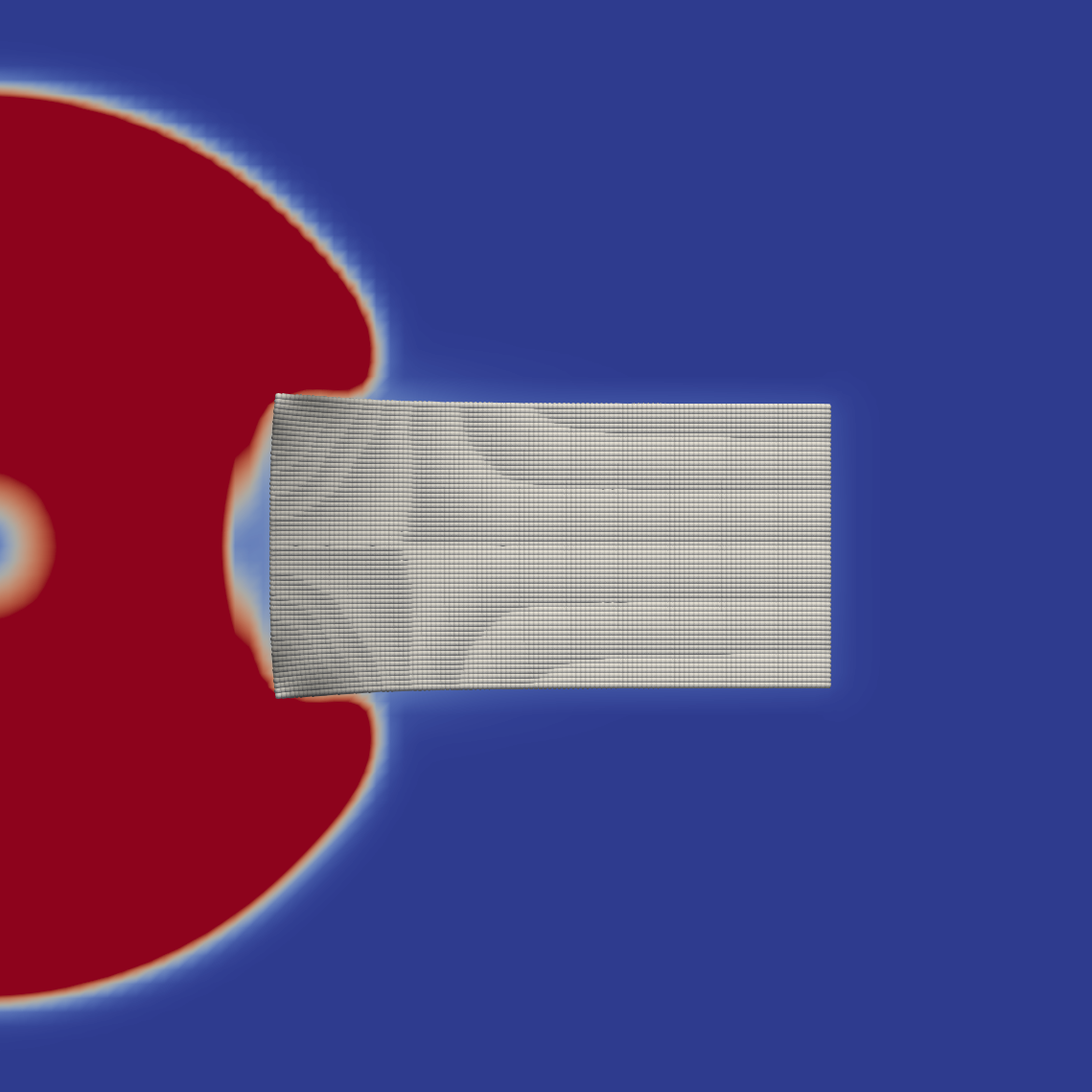}}
  \hspace{7pt}
  \subfloat{\includegraphics[width=0.21\textwidth,trim={0cm 0cm 0cm 0cm},clip]{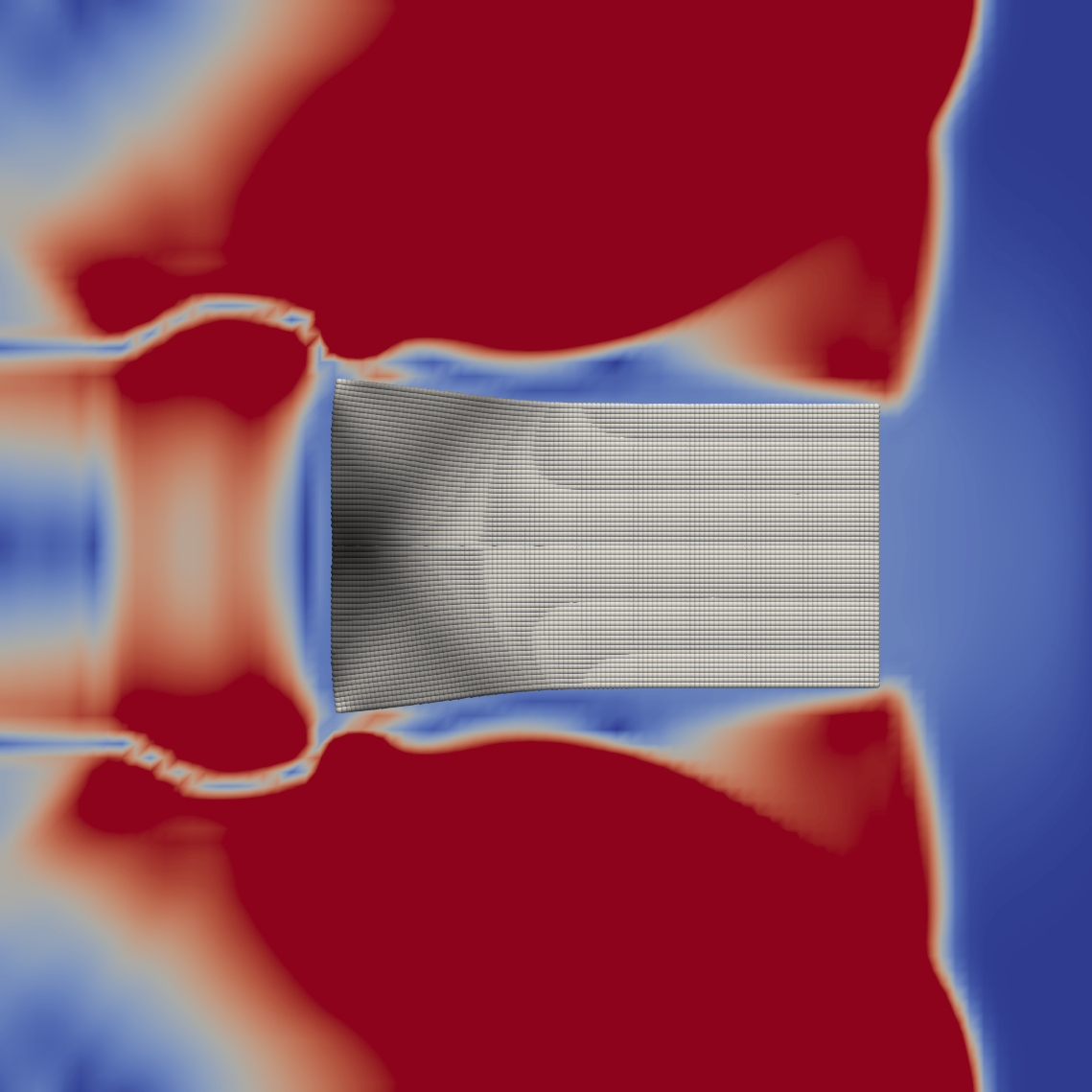}}
  \hspace{7pt}
  \subfloat{\includegraphics[width=0.21\textwidth,trim={0cm 0cm 0cm 0cm},clip]{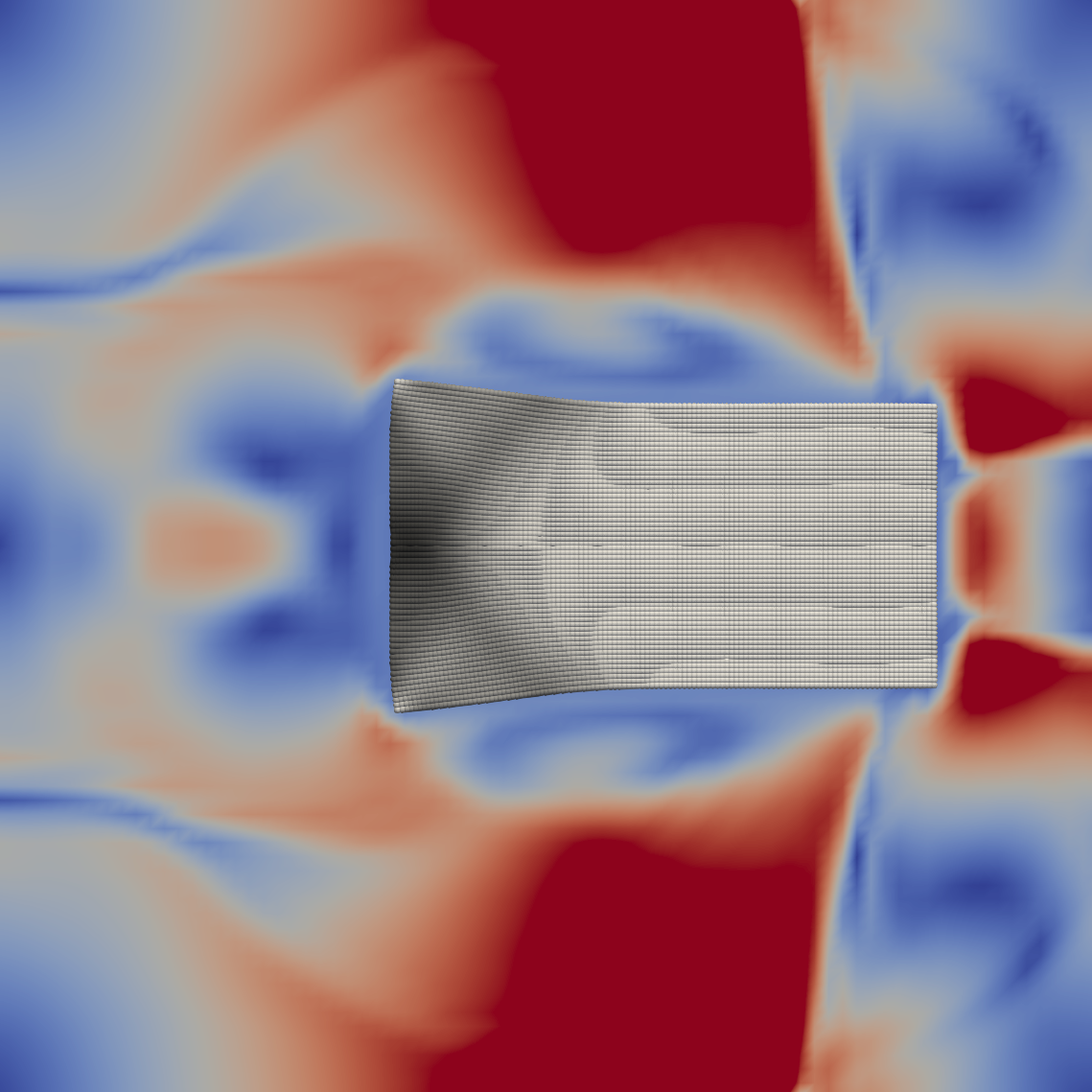}}
  \hspace{7pt}
  \subfloat{\includegraphics[width=0.26\textwidth,trim={0cm 0cm 0cm 0cm},clip]{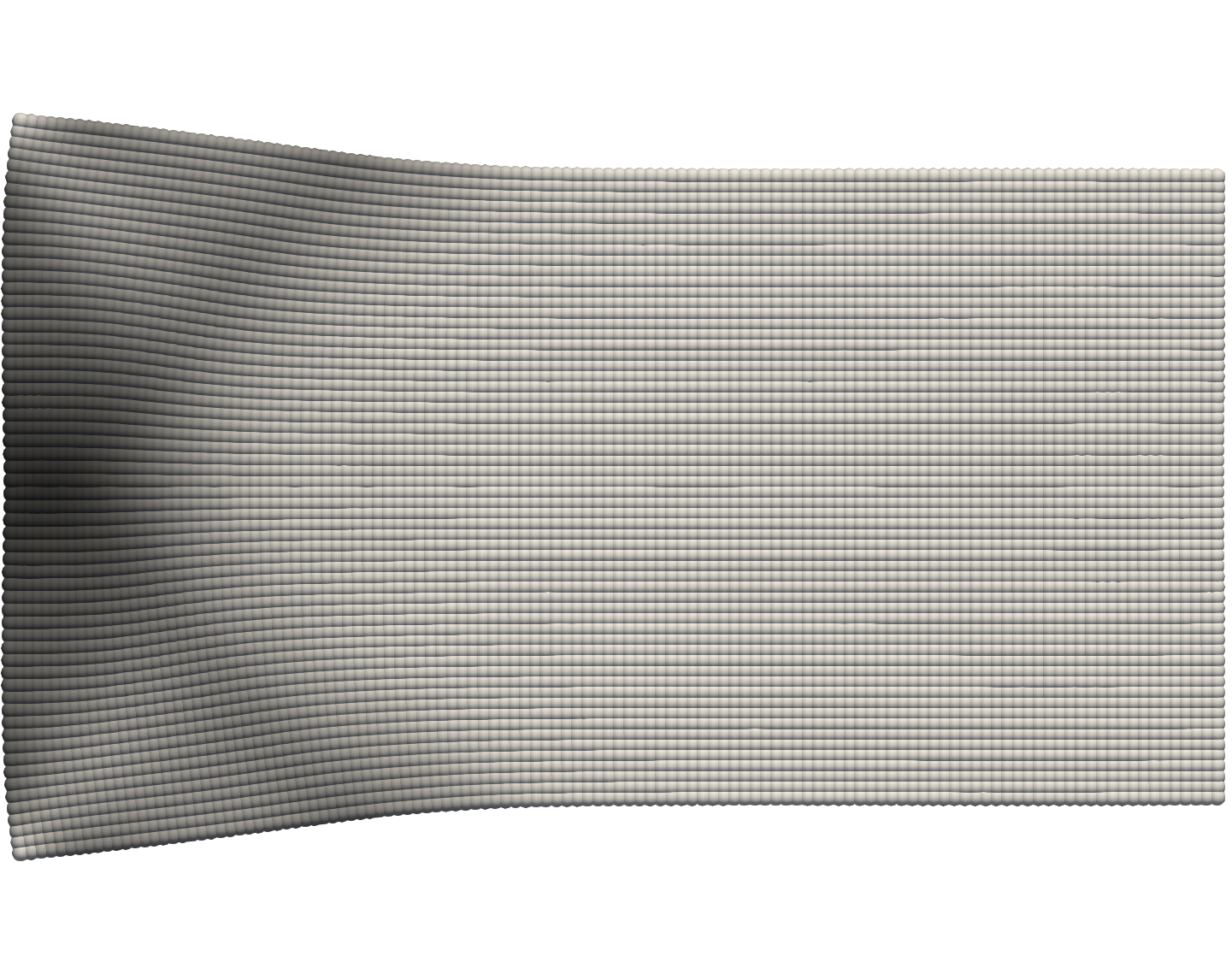}}
  
  \subfloat{\includegraphics[width=0.21\textwidth,trim={0cm 0cm 0cm 0cm},clip]{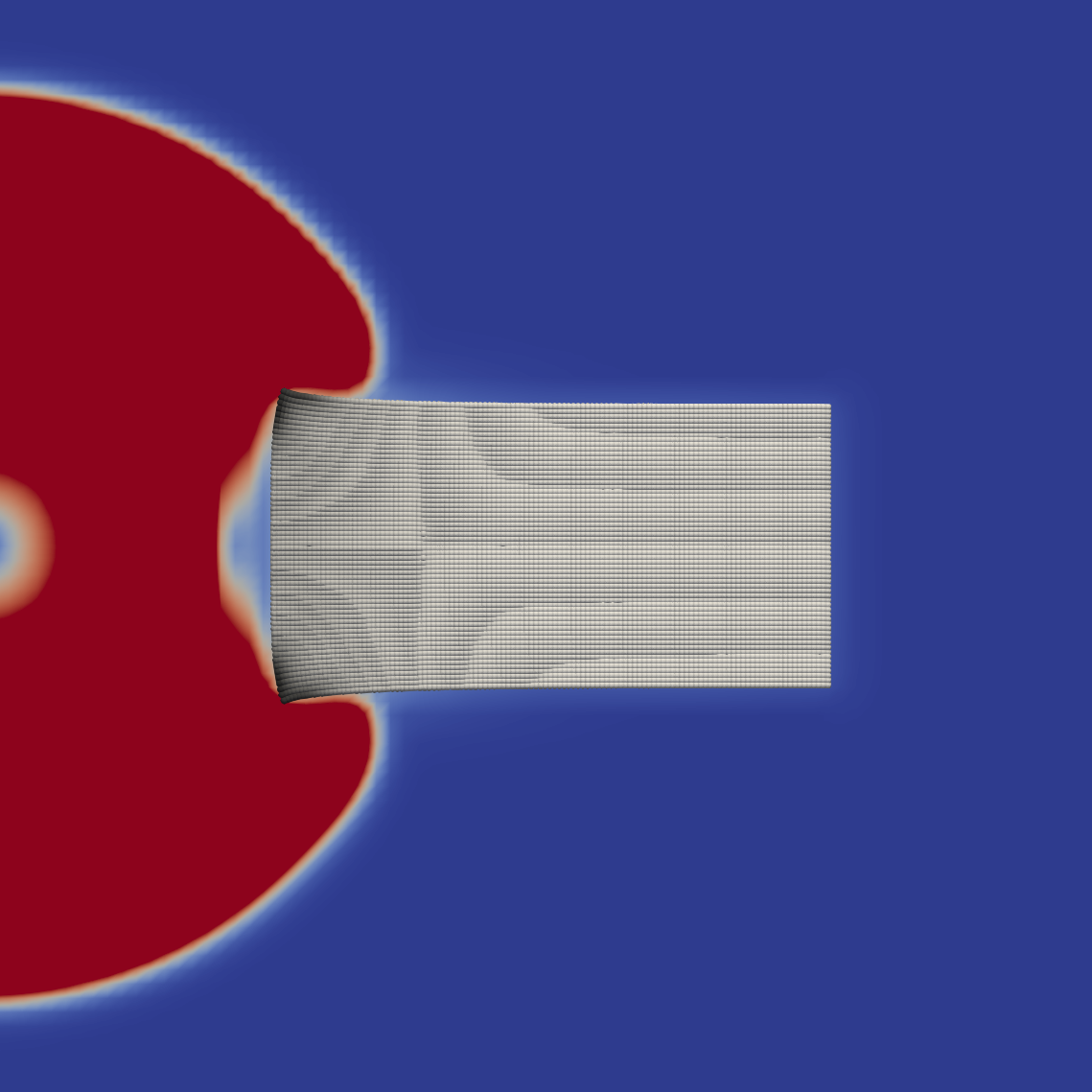}}
  \hspace{7pt}
  \subfloat{\includegraphics[width=0.21\textwidth,trim={0cm 0cm 0cm 0cm},clip]{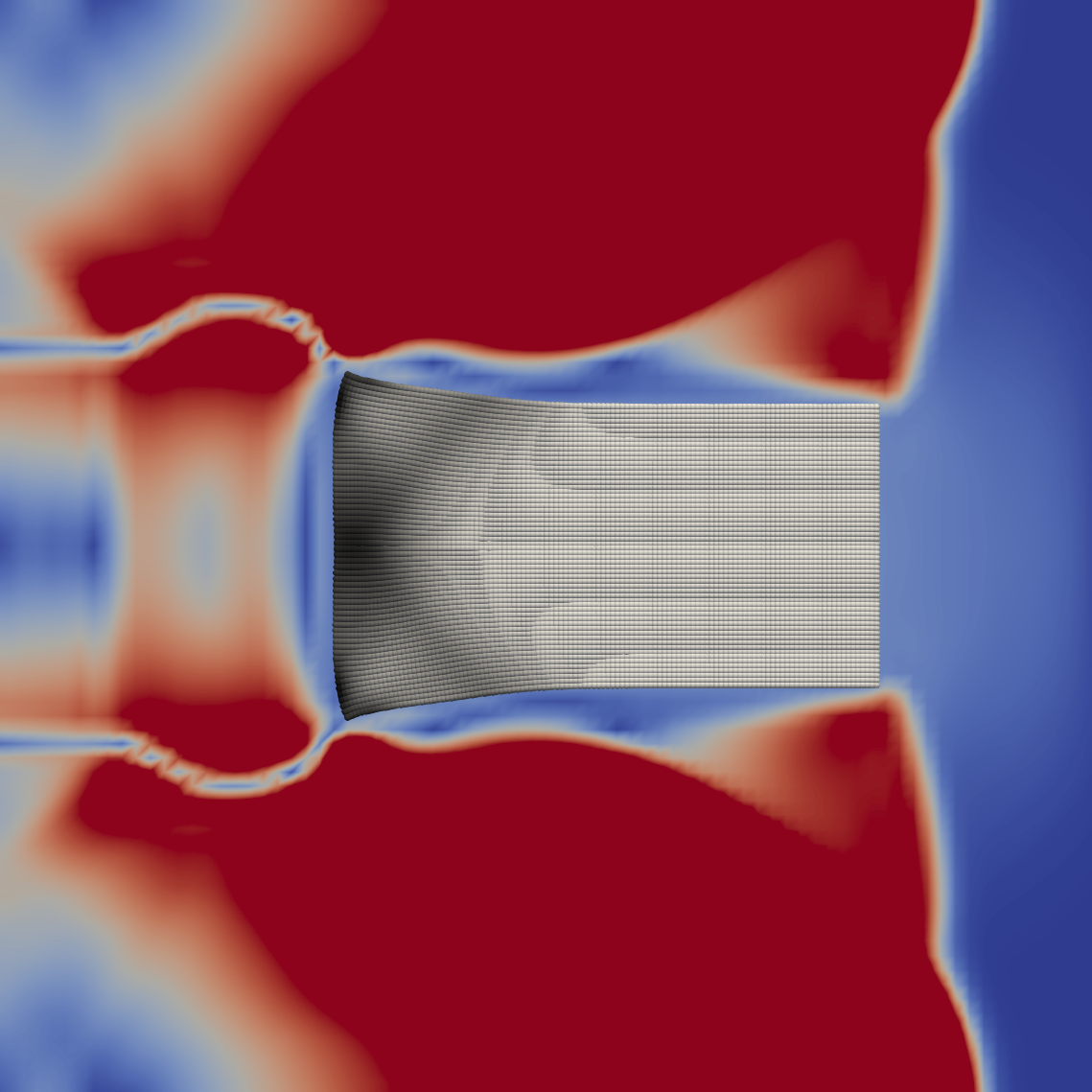}}
  \hspace{7pt}
  \subfloat{\includegraphics[width=0.21\textwidth,trim={0cm 0cm 0cm 0cm},clip]{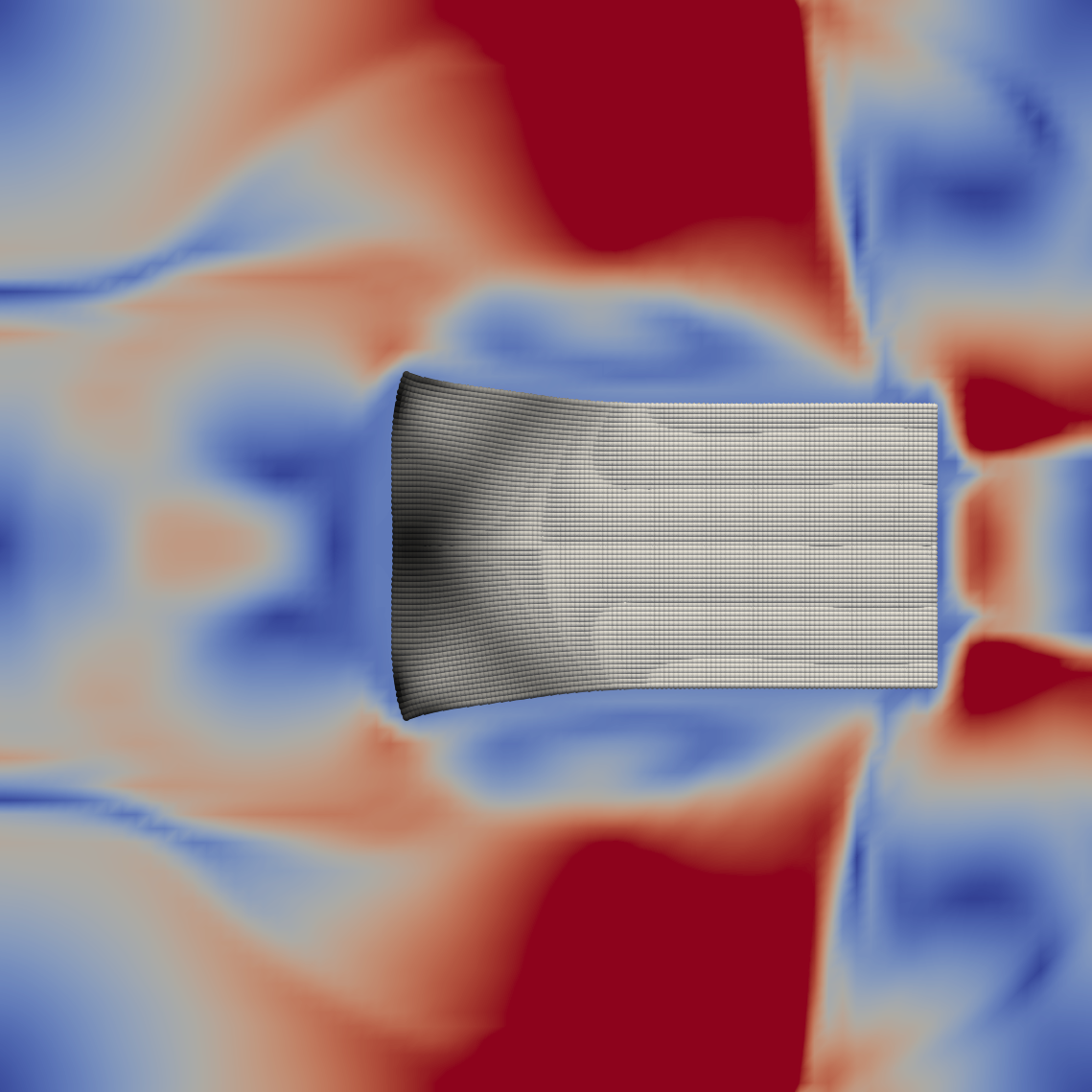}}
  \hspace{7pt}
  \subfloat{\includegraphics[width=0.26\textwidth,trim={0cm 0cm 0cm 0cm},clip]{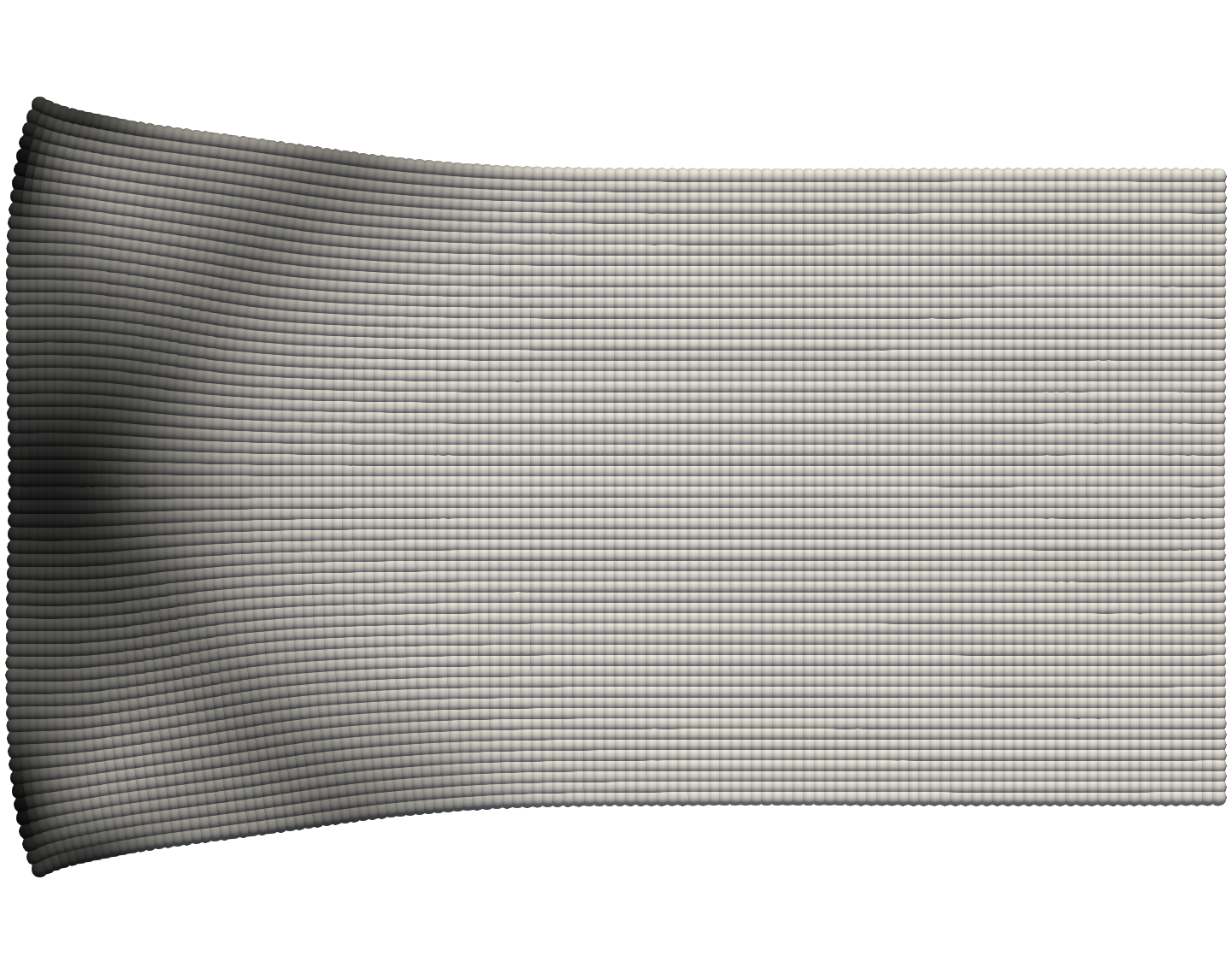}}
  
  \subfloat{\includegraphics[width=0.21\textwidth,trim={0cm 0cm 0cm 0cm},clip]{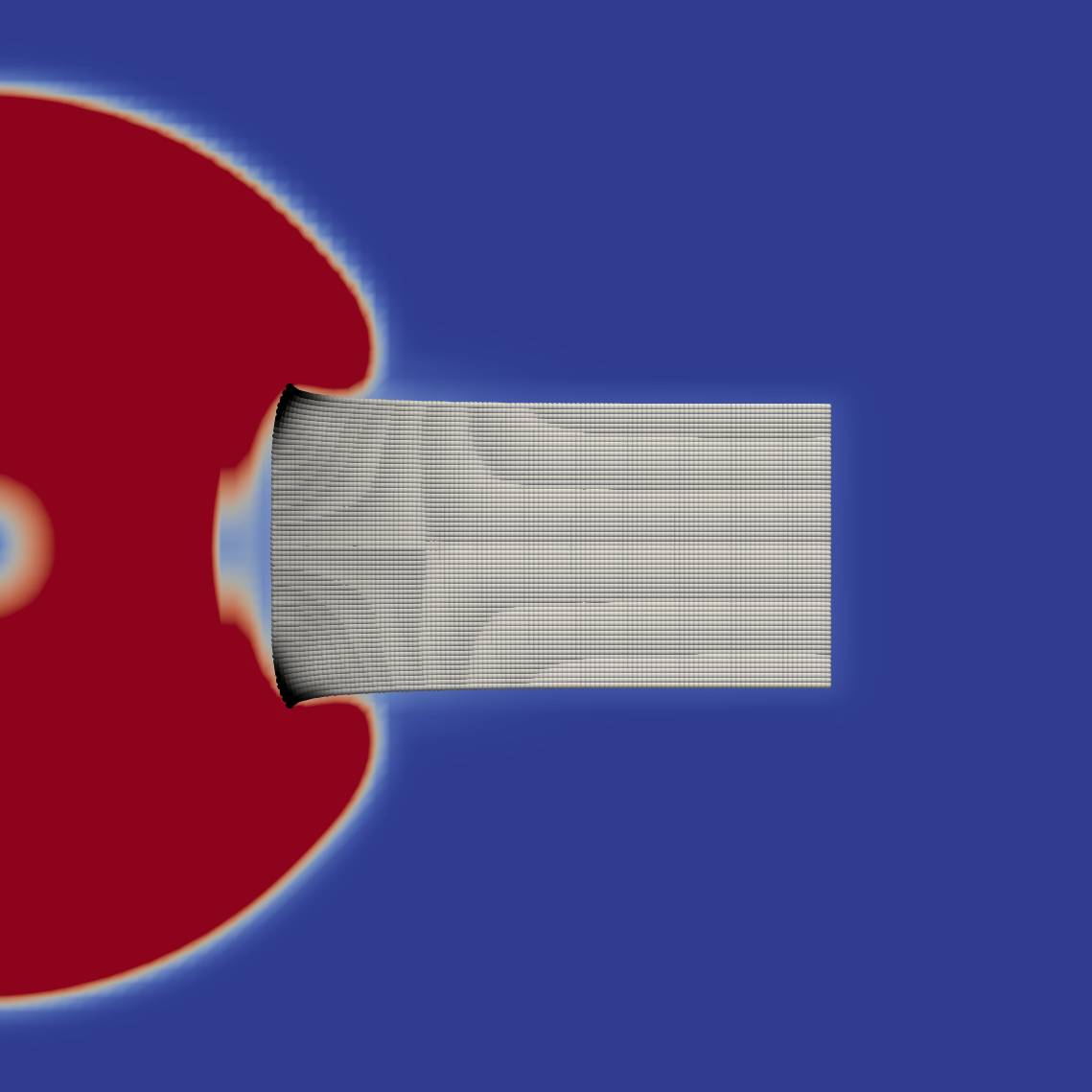}}
  \hspace{7pt}
  \subfloat{\includegraphics[width=0.21\textwidth,trim={0cm 0cm 0cm 0cm},clip]{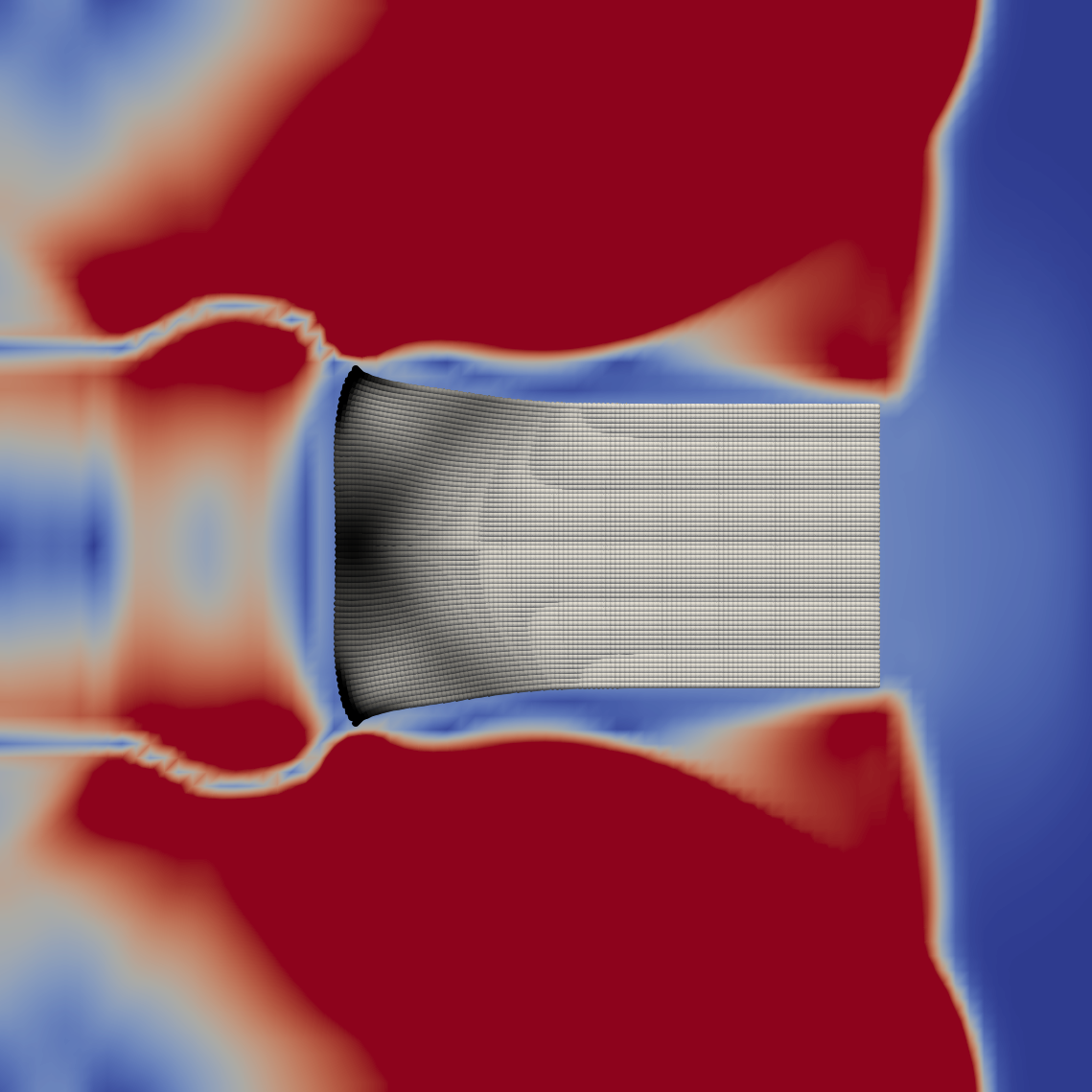}}
  \hspace{7pt}
  \subfloat{\includegraphics[width=0.21\textwidth,trim={0cm 0cm 0cm 0cm},clip]{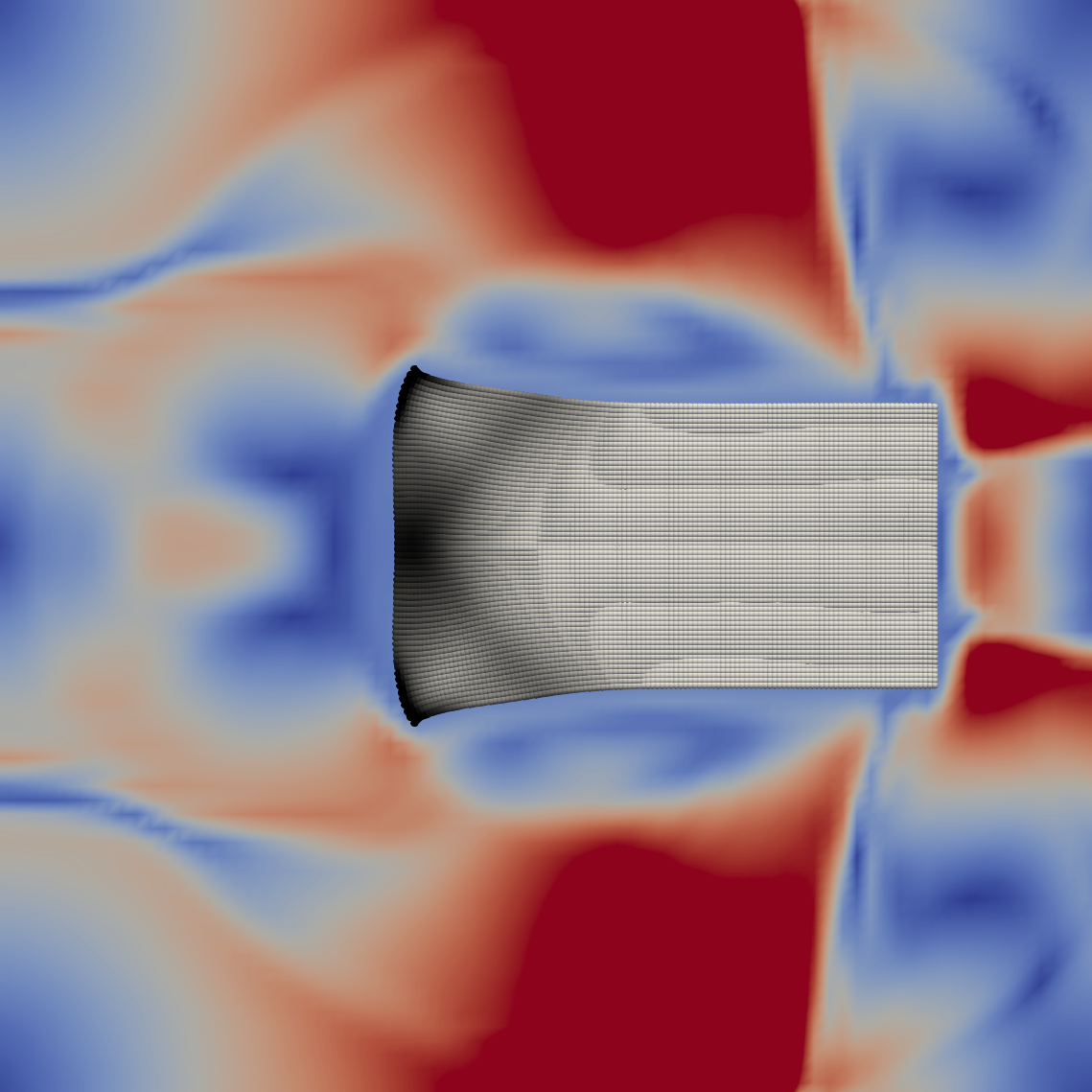}}
  \hspace{7pt}
  \subfloat{\includegraphics[width=0.26\textwidth,trim={0cm 0cm 0cm 0cm},clip]{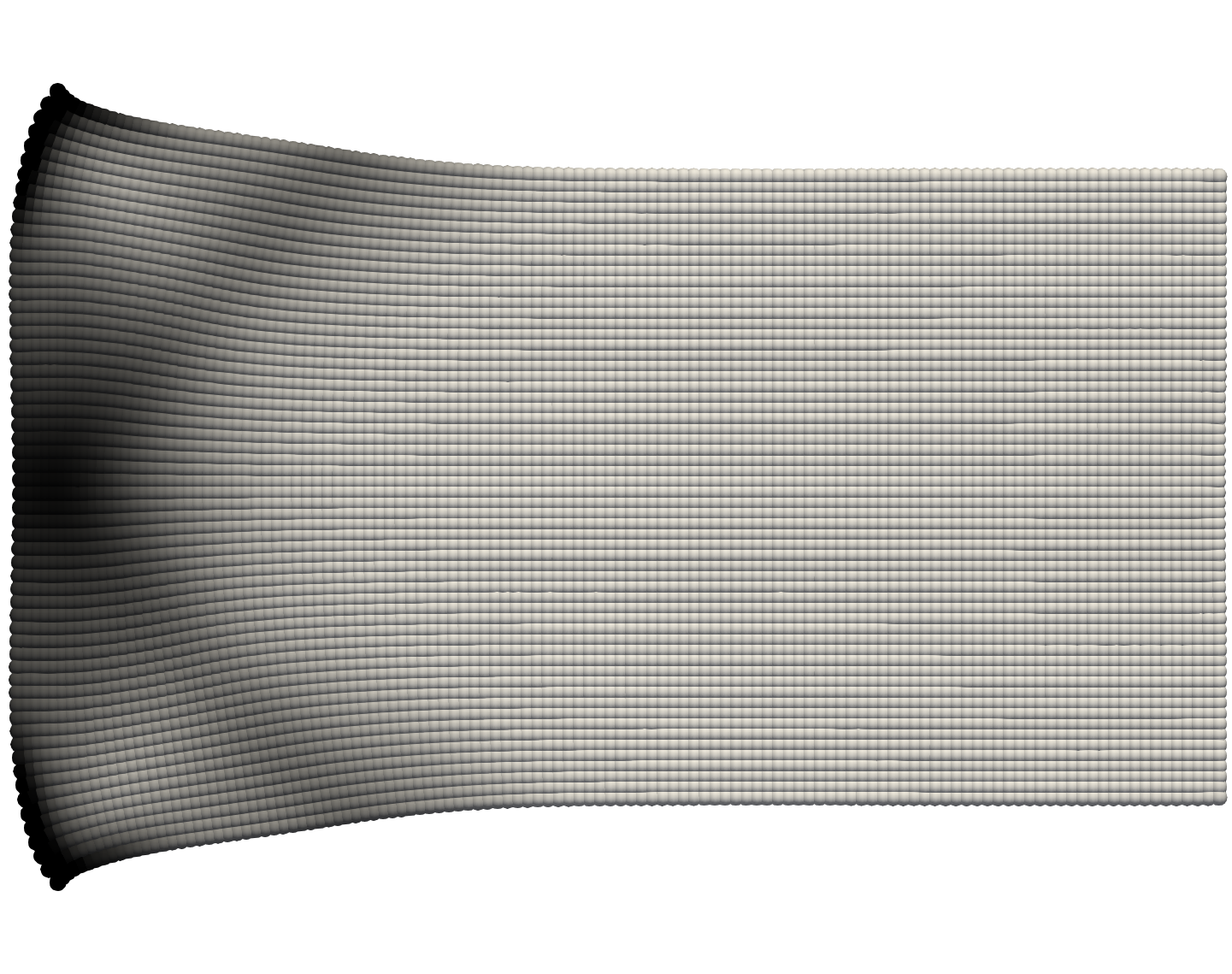}}
  
  \setcounter{subfigure}{0}

  \subfloat[][$t = \SI{0.1}{\milli s}$]{\includegraphics[width=0.21\textwidth,trim={0cm 0cm 0cm 0cm},clip]{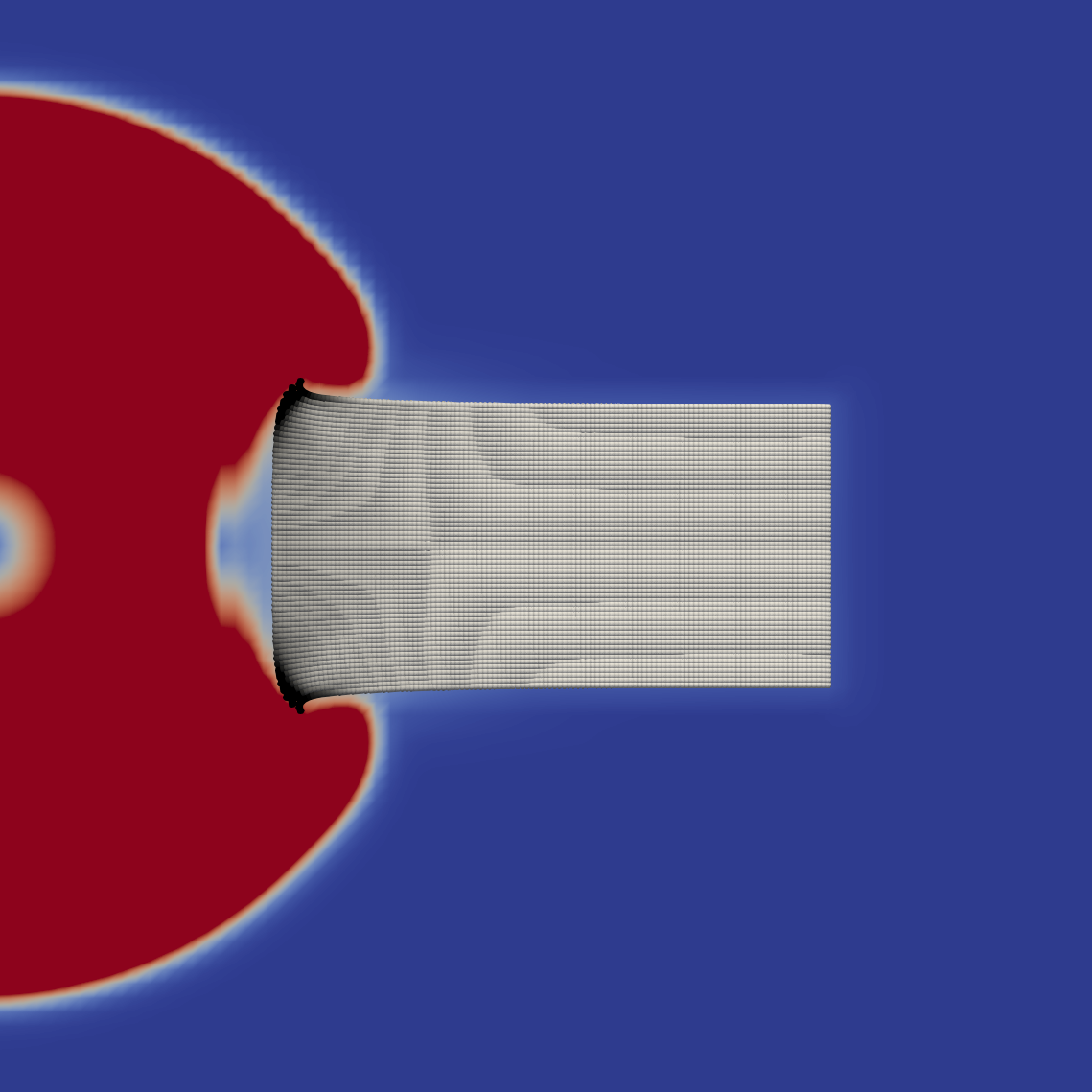}}
  \hspace{7pt}
  \subfloat[][$t = \SI{0.4}{\milli s}$]{\includegraphics[width=0.21\textwidth,trim={0cm 0cm 0cm 0cm},clip]{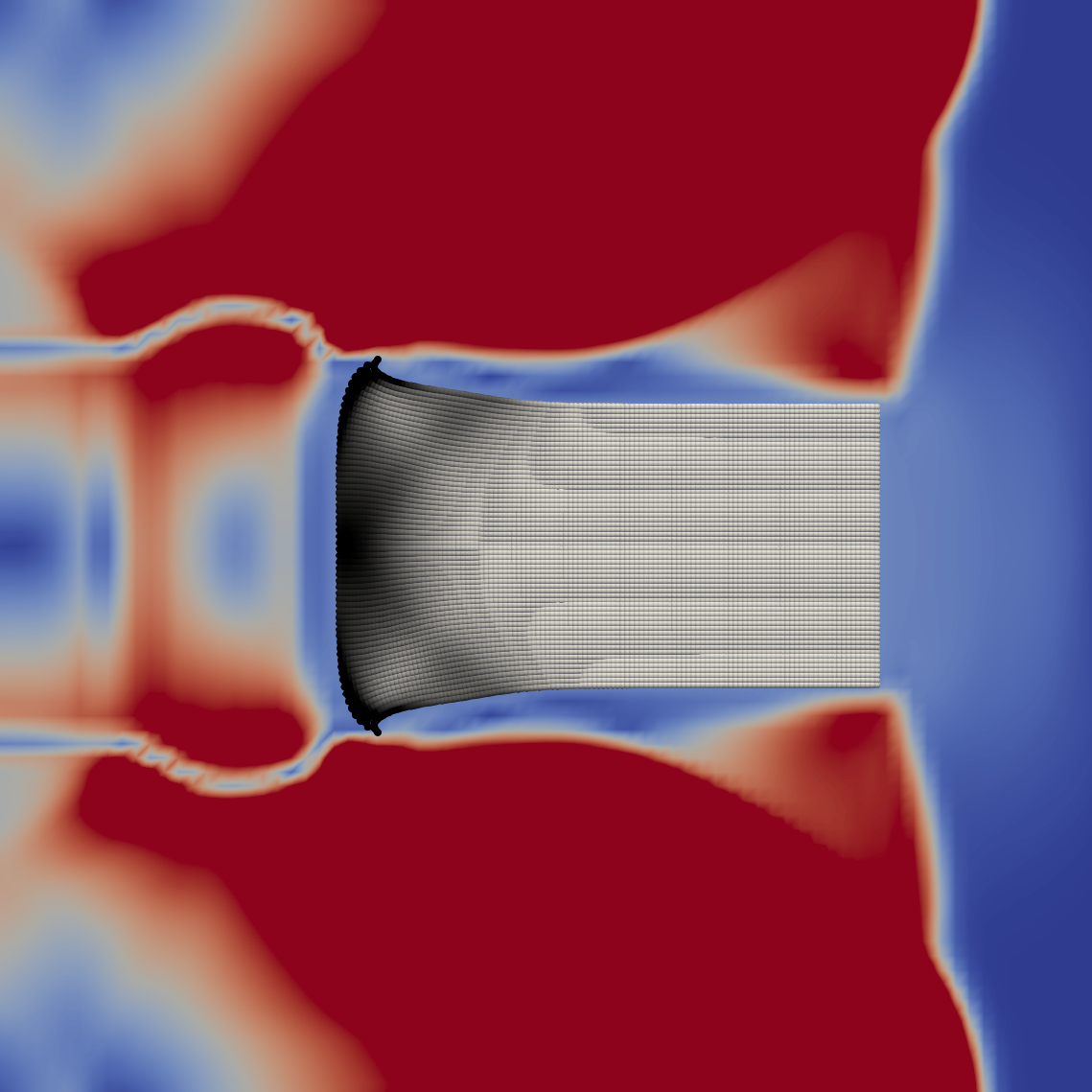}}
  \hspace{7pt}
  \subfloat[][$t = \SI{0.7}{\milli s}$]{\includegraphics[width=0.21\textwidth,trim={0cm 0cm 0cm 0cm},clip]{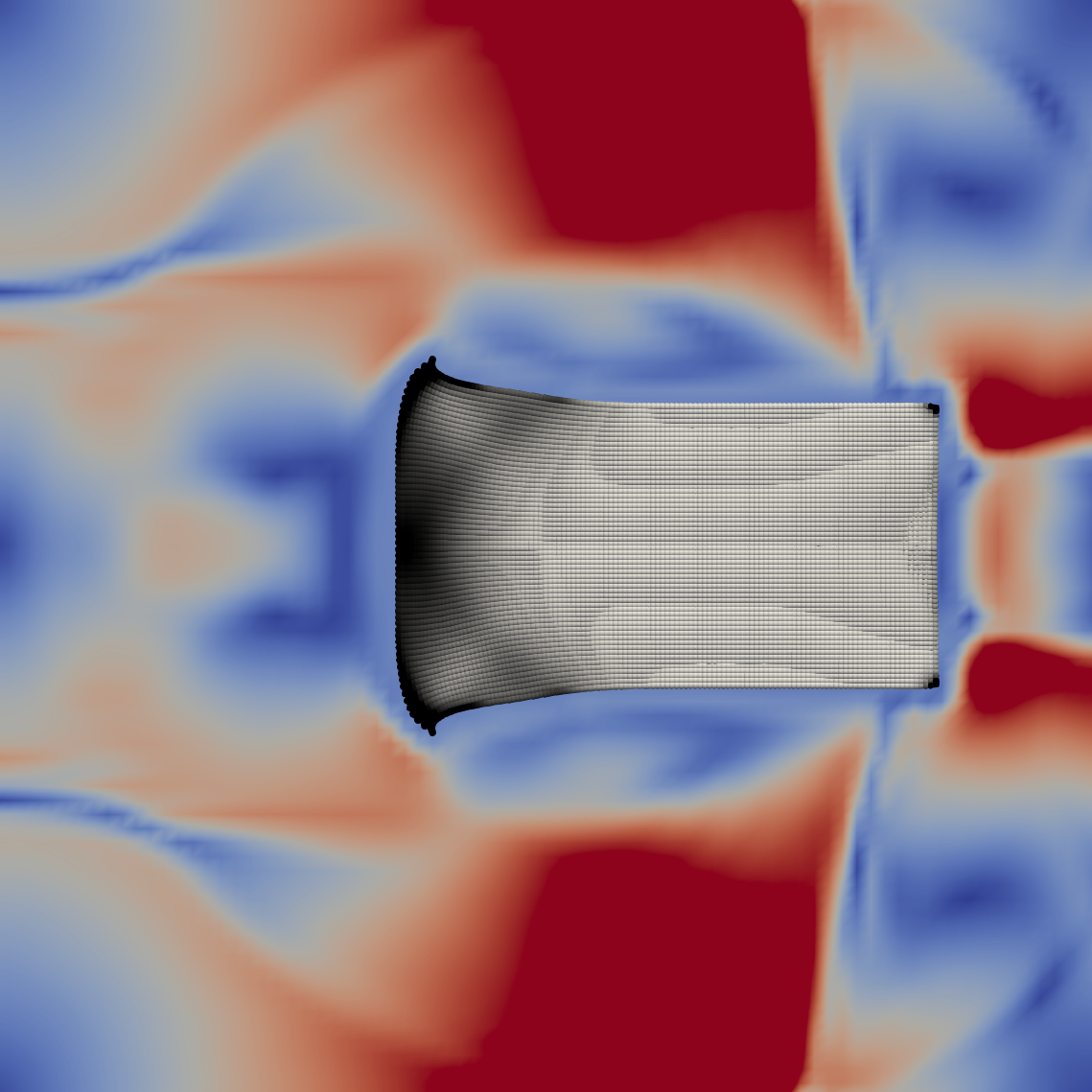}}
  \hspace{7pt}
  \subfloat[][$t = \SI{0.7}{\milli s}$]{\includegraphics[width=0.26\textwidth,trim={0cm 0cm 0cm 0cm},clip]{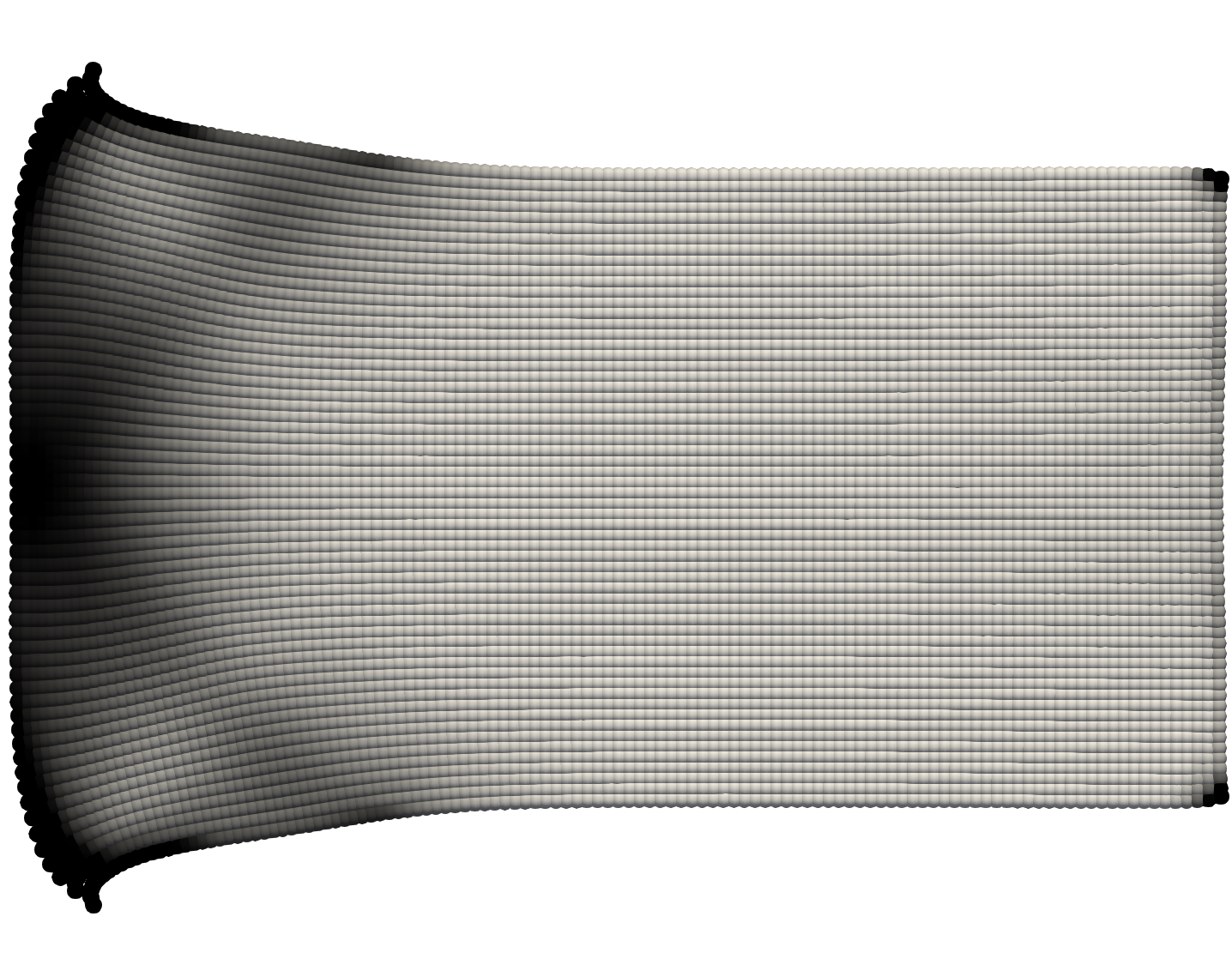}}

  \caption{Chamber detonation problem. (a--c): Snapshots of air speed (in m/s) and solid plastic strain in the current configuration at different times for the finest mesh. From top to bottom, the individual rows correspond to $\beta = 1/3$, $\beta = 1.0$, $\beta = 3.0$, $\beta = 9.0$, and strong coupling, respectively. (d): Solid plastic strain at $t = \SI{0.7}{\milli s}$.}
  \label{fig:plastic_contours}
\end{figure*}

Air speed and solid plastic strain contours at several time instants are compared between the different cases computed on the finest mesh in \cref{fig:plastic_contours}. The fluid response appears to be very similar in all cases. As the penalty constant increases, plastic contours in the solid domain become more pronounced, more so near the domain boundaries and, especially, at the corners. The strong coupling produces excessive deformation and plastic straining at the domain corners because the discrete solution tries to accommodate high velocity gradients near the corners of the solid. Switching to weak coupling allows for some mismatch in the kinematics, and the domain corners experience little distortion, which is a real advantage of the proposed weak coupling.

\begin{figure*}[!hbpt]
  \centering
   \subfloat[][]{\includegraphics[width=0.47\textwidth,trim={0cm 0cm 0cm 0cm},clip]{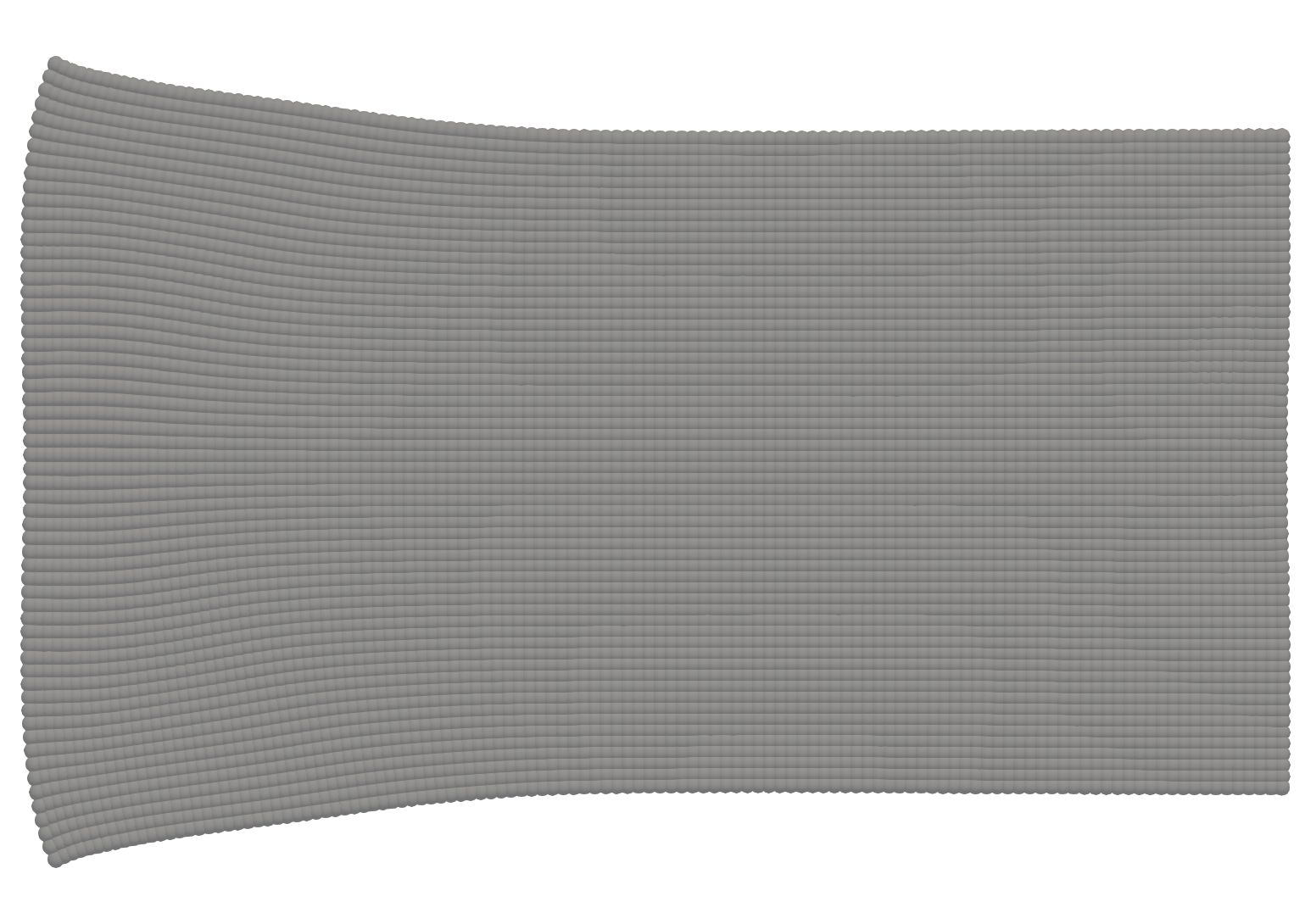}}
  \hfill
  \subfloat[][]{\includegraphics[width=0.47\textwidth,trim={0cm 0cm 0cm 0cm},clip]{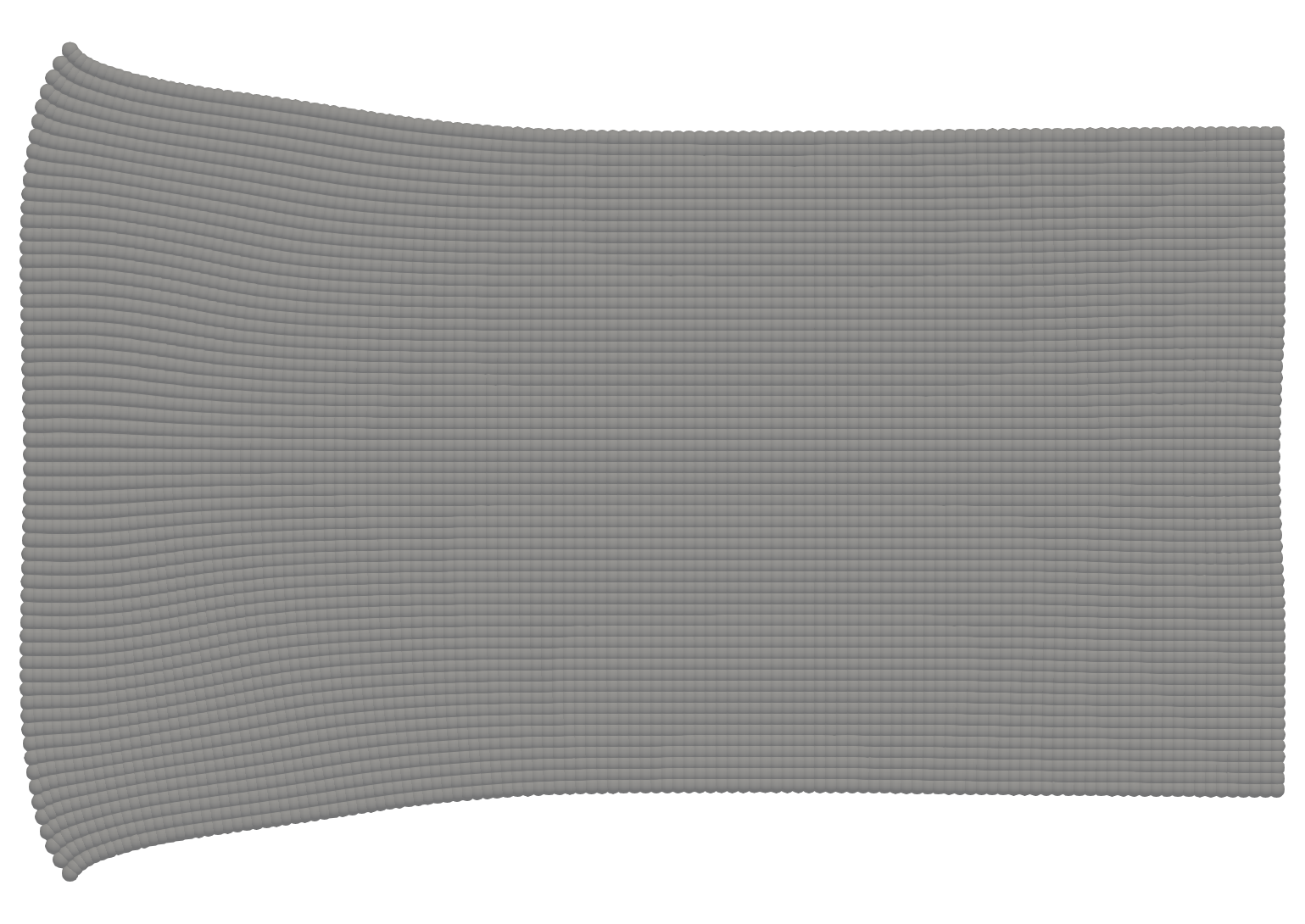}}
  
  \subfloat[][]{\includegraphics[width=0.47\textwidth,trim={0cm 0cm 0cm 0cm},clip]{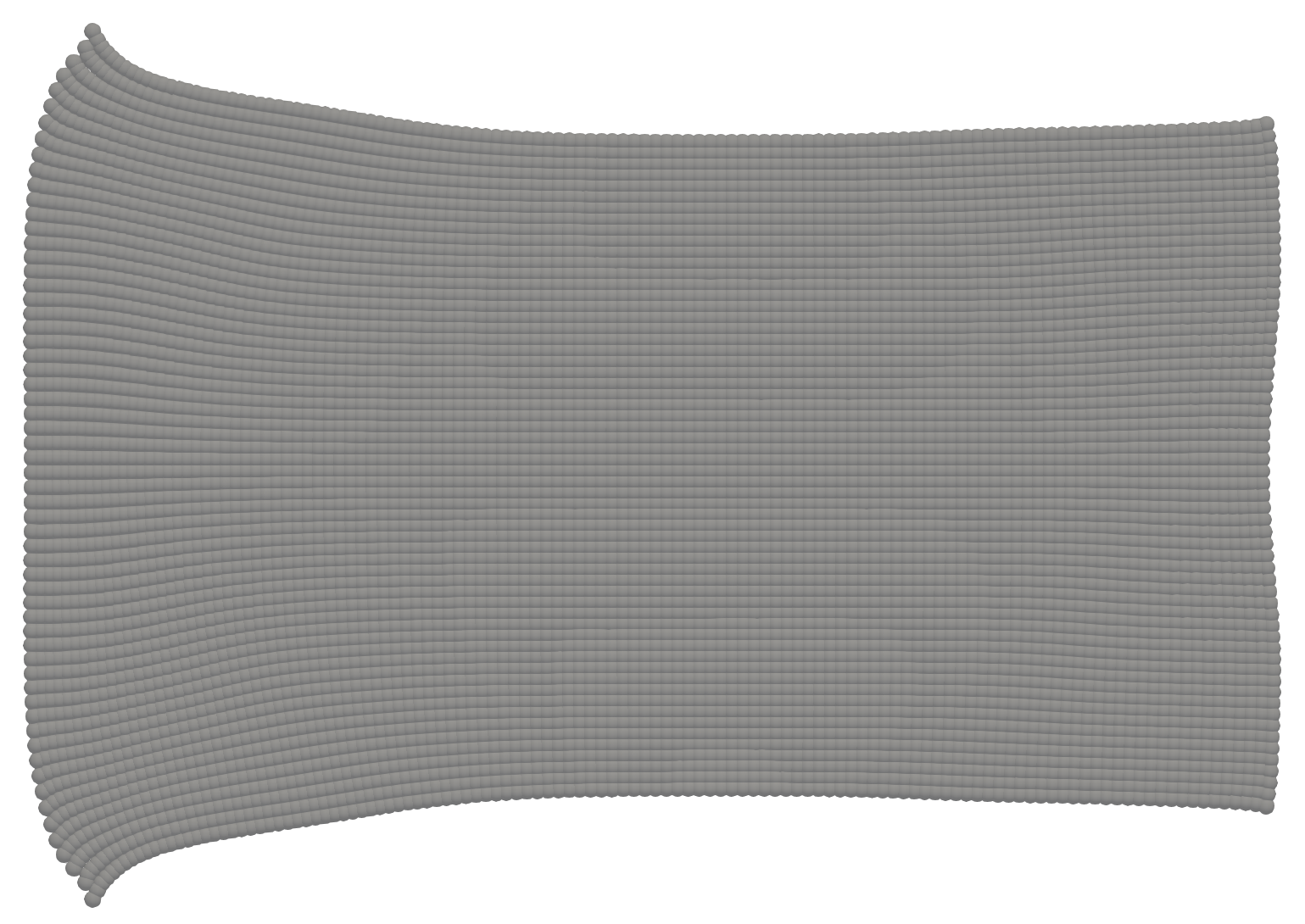}}
  \hfill
  \subfloat[][]{\includegraphics[width=0.47\textwidth,trim={14.5cm 12cm 4.5cm 12cm},clip]{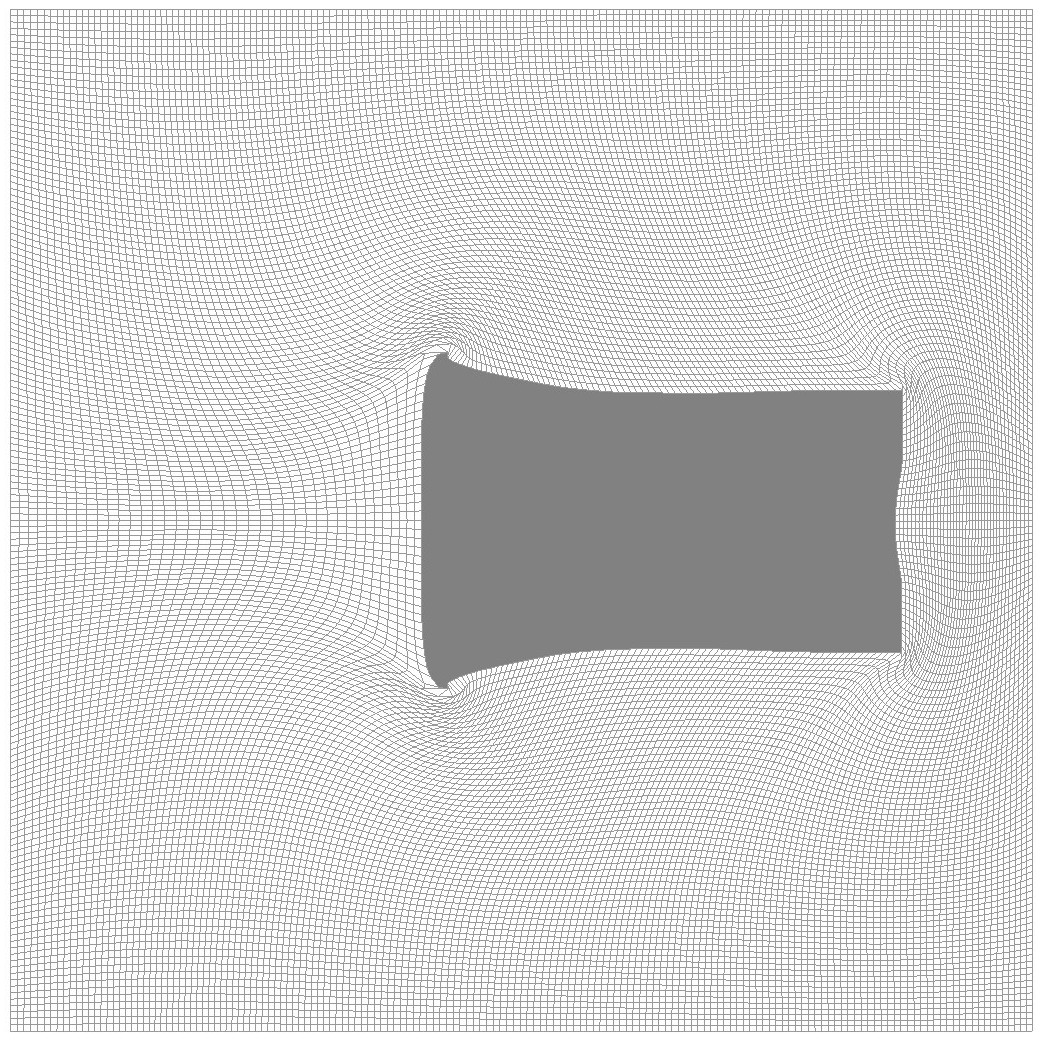}}
  \caption{Chamber detonation problem. The solid deformed shape at $t = \SI{1.5}{\milli s}$. (a): $\beta = 3.0$. (b): $\beta = 9.0$. (c): Strong coupling. (d): Conforming-mesh ALE result from~\cite{bazilevs2017new1}.}
  \label{fig:plastic_solid}
\end{figure*}

\begin{figure*}[!hbpt]
  \centering
  \subfloat{\includegraphics[width=0.36\textwidth,height=0.27\textwidth,trim={0cm 1.885cm 0cm 0cm},clip]{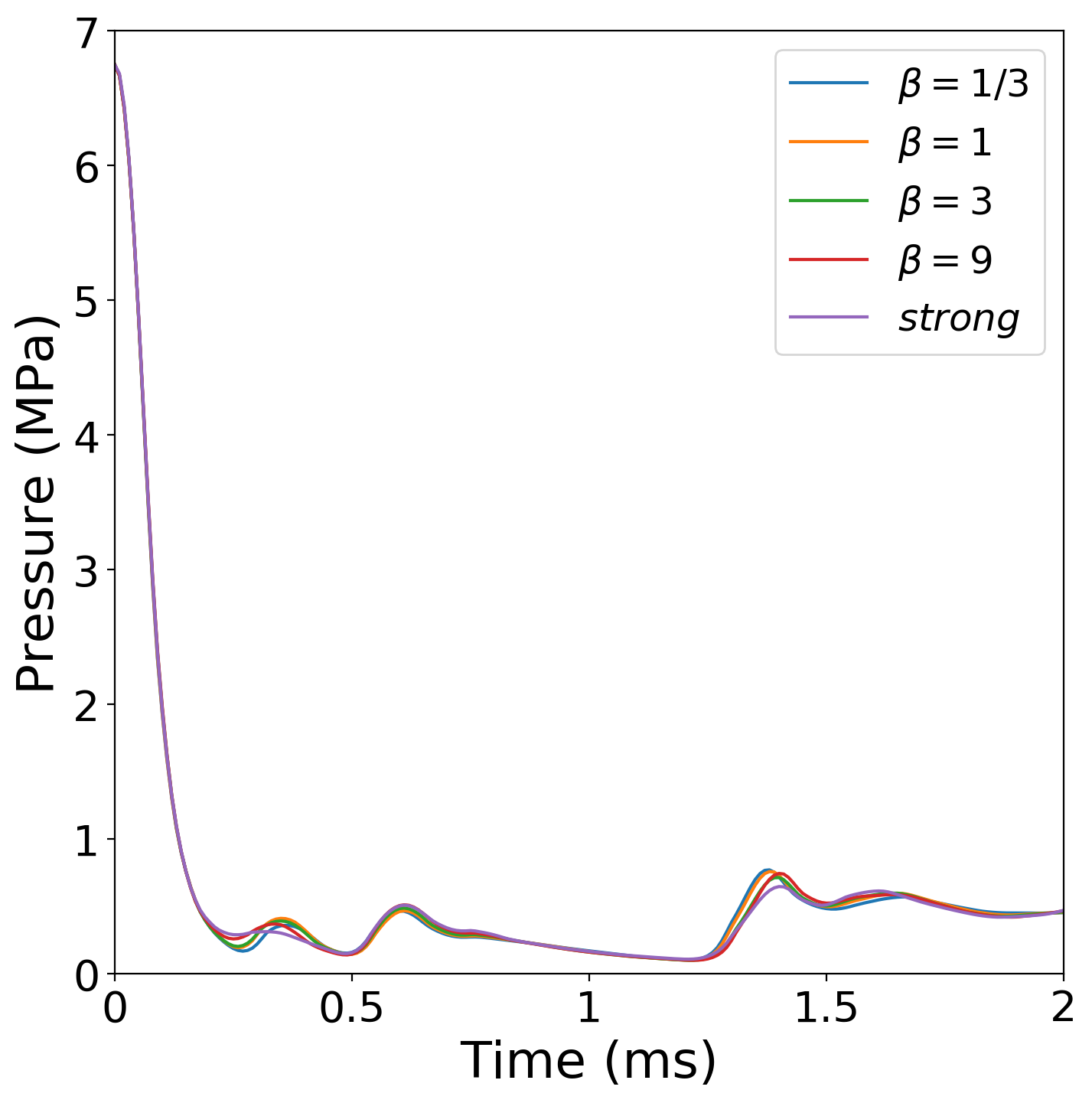}}
  \hspace{4pt}
  \subfloat{\includegraphics[width=0.303\textwidth,height=0.27\textwidth,trim={1.7cm 1.885cm 0cm 0cm},clip]{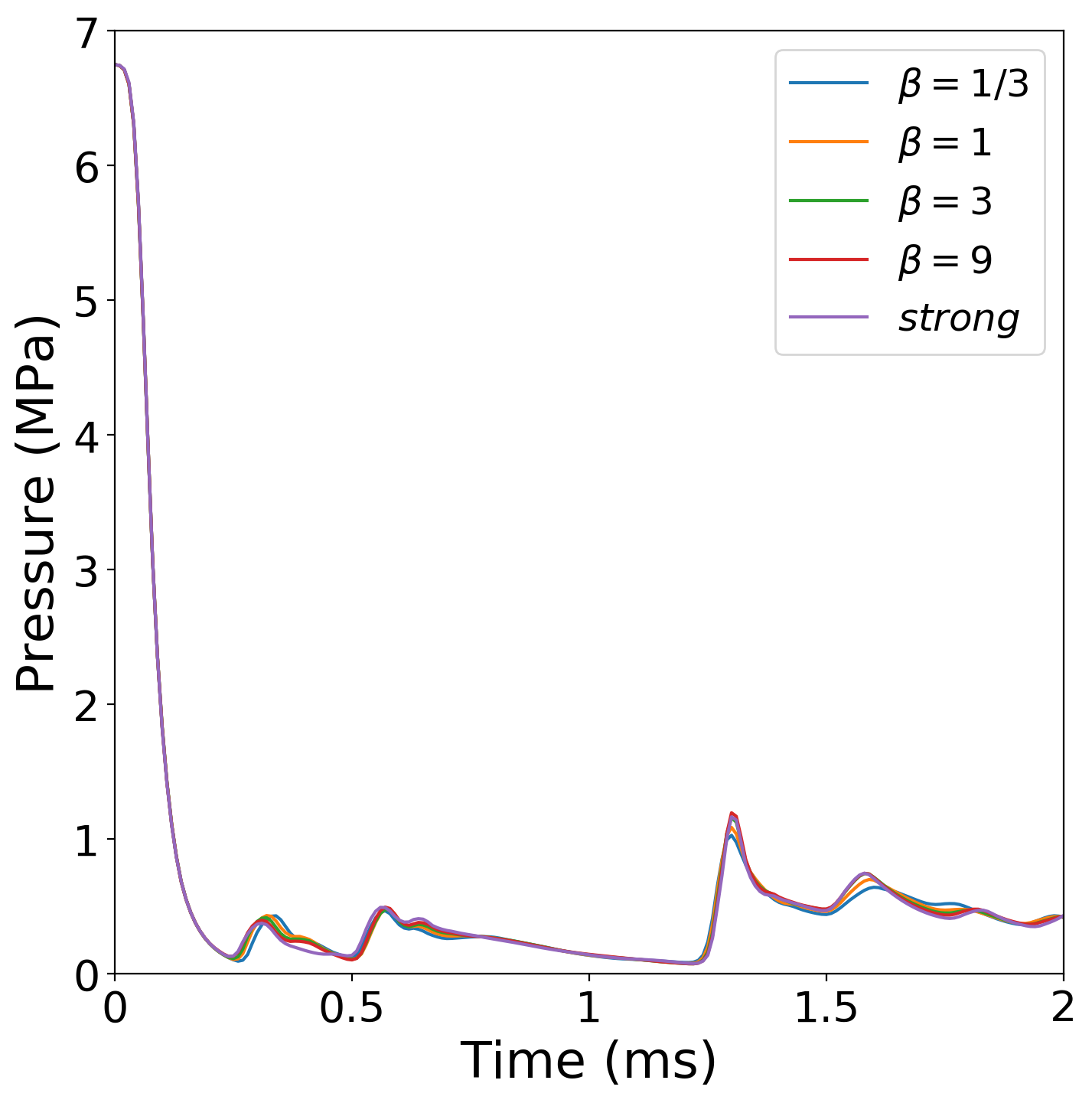}}
  \hspace{4pt}
  \subfloat{\includegraphics[width=0.303\textwidth,height=0.27\textwidth,trim={1.7cm 1.885cm 0cm 0cm},clip]{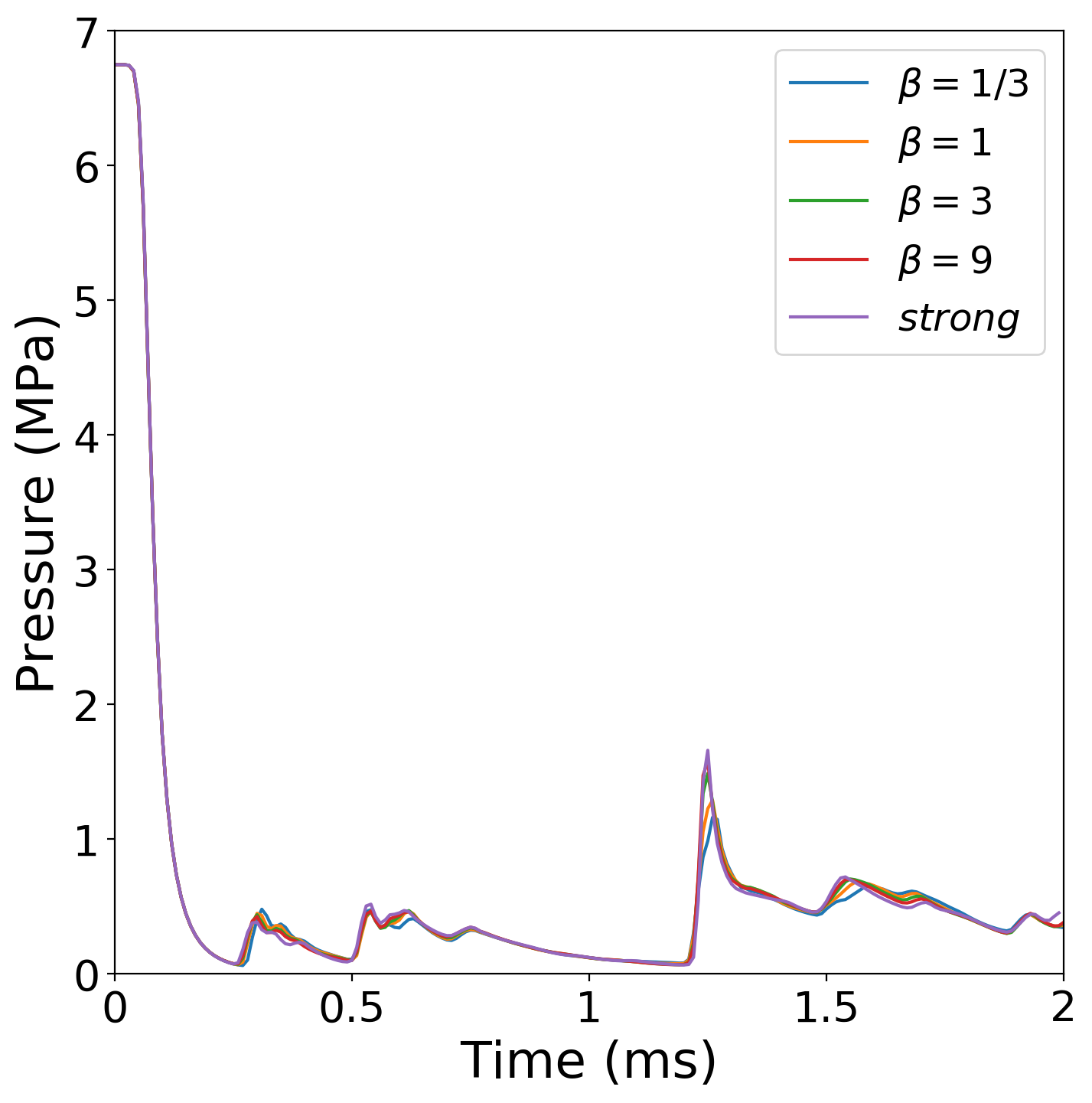}}
  
  \subfloat{\includegraphics[width=0.36\textwidth,height=0.27\textwidth,trim={0cm 1.885cm 0cm 0cm},clip]{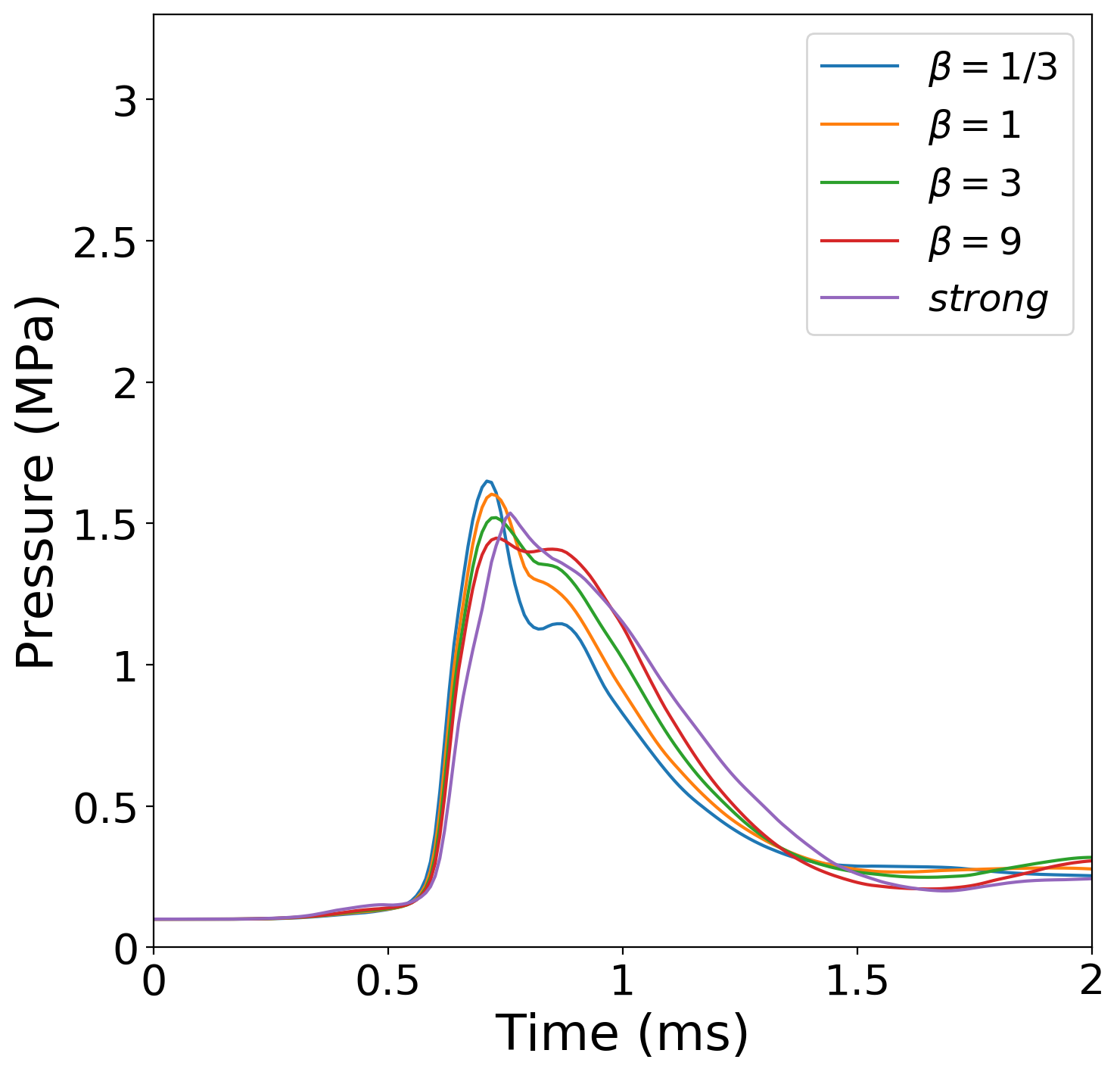}}
  \hspace{4pt}
  \subfloat{\includegraphics[width=0.303\textwidth,height=0.27\textwidth,trim={2.45cm 1.885cm 0cm 0cm},clip]{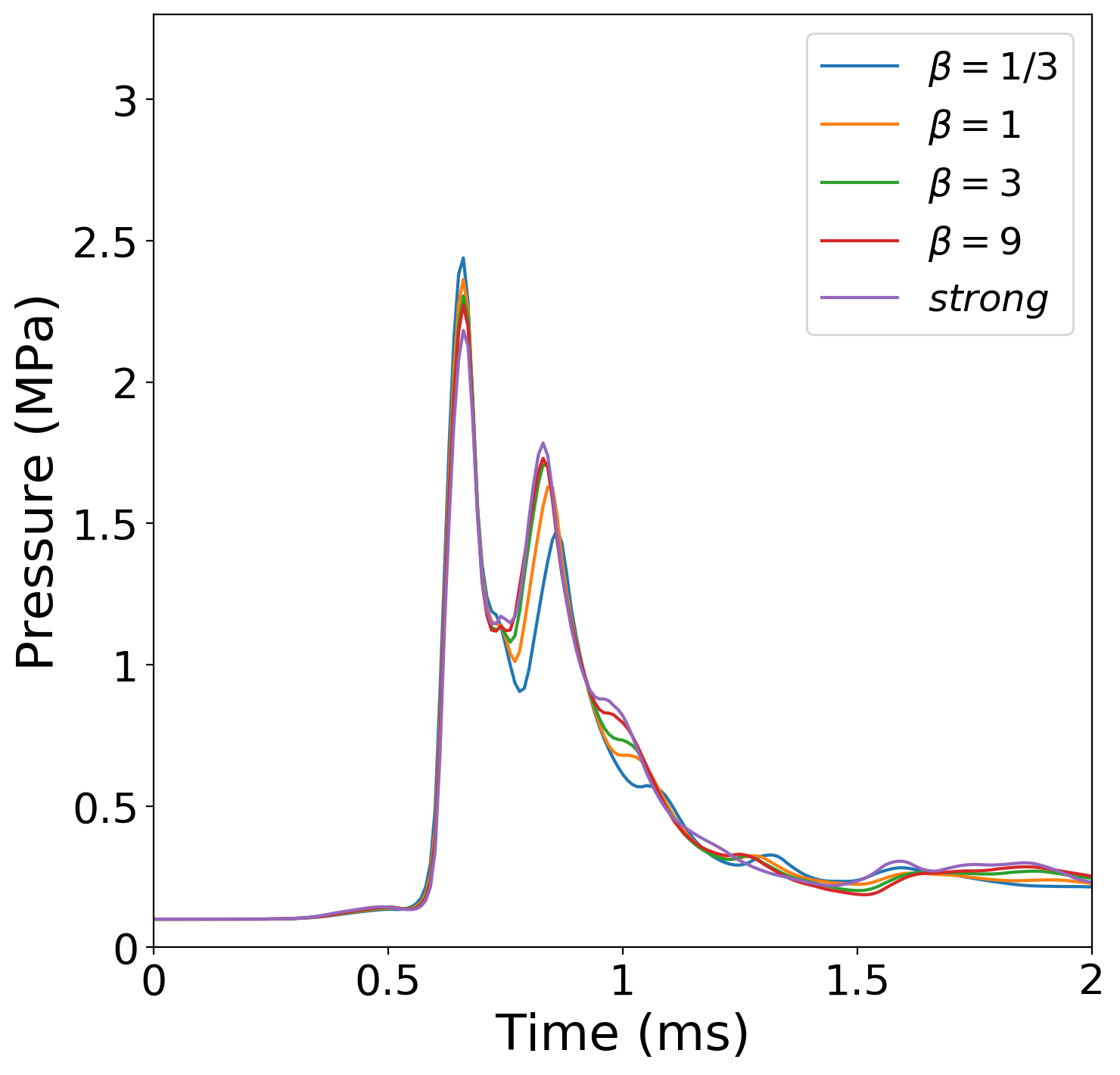}}
  \hspace{4pt}
  \subfloat{\includegraphics[width=0.303\textwidth,height=0.27\textwidthwidth=0.303\textwidth,height=0.28\textwidth,trim={2.45cm 1.885cm 0cm 0cm},clip]{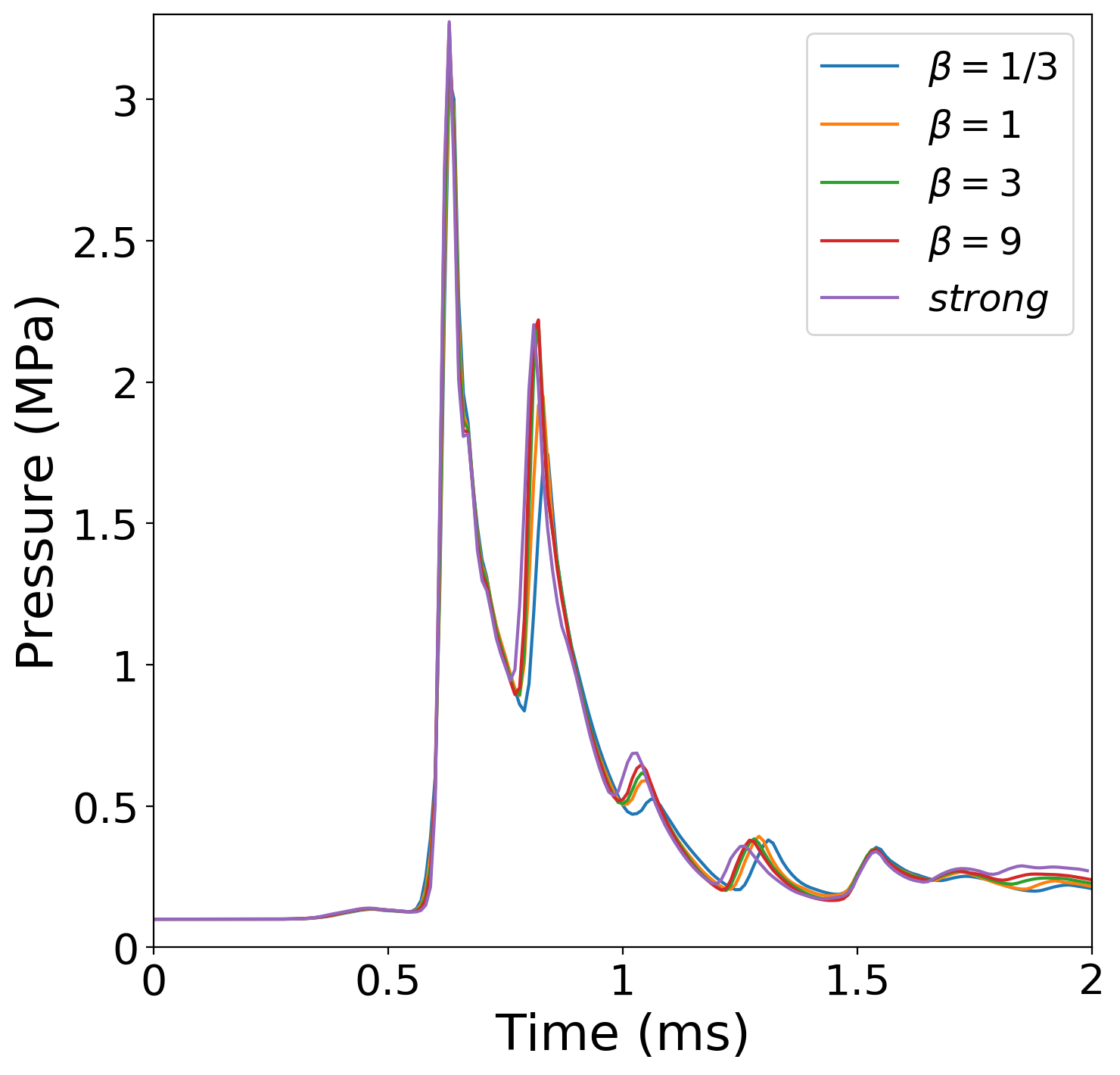}}
  
  \subfloat{\includegraphics[width=0.36\textwidth,height=0.27\textwidth,trim={0cm 1.885cm 0cm 0cm},clip]{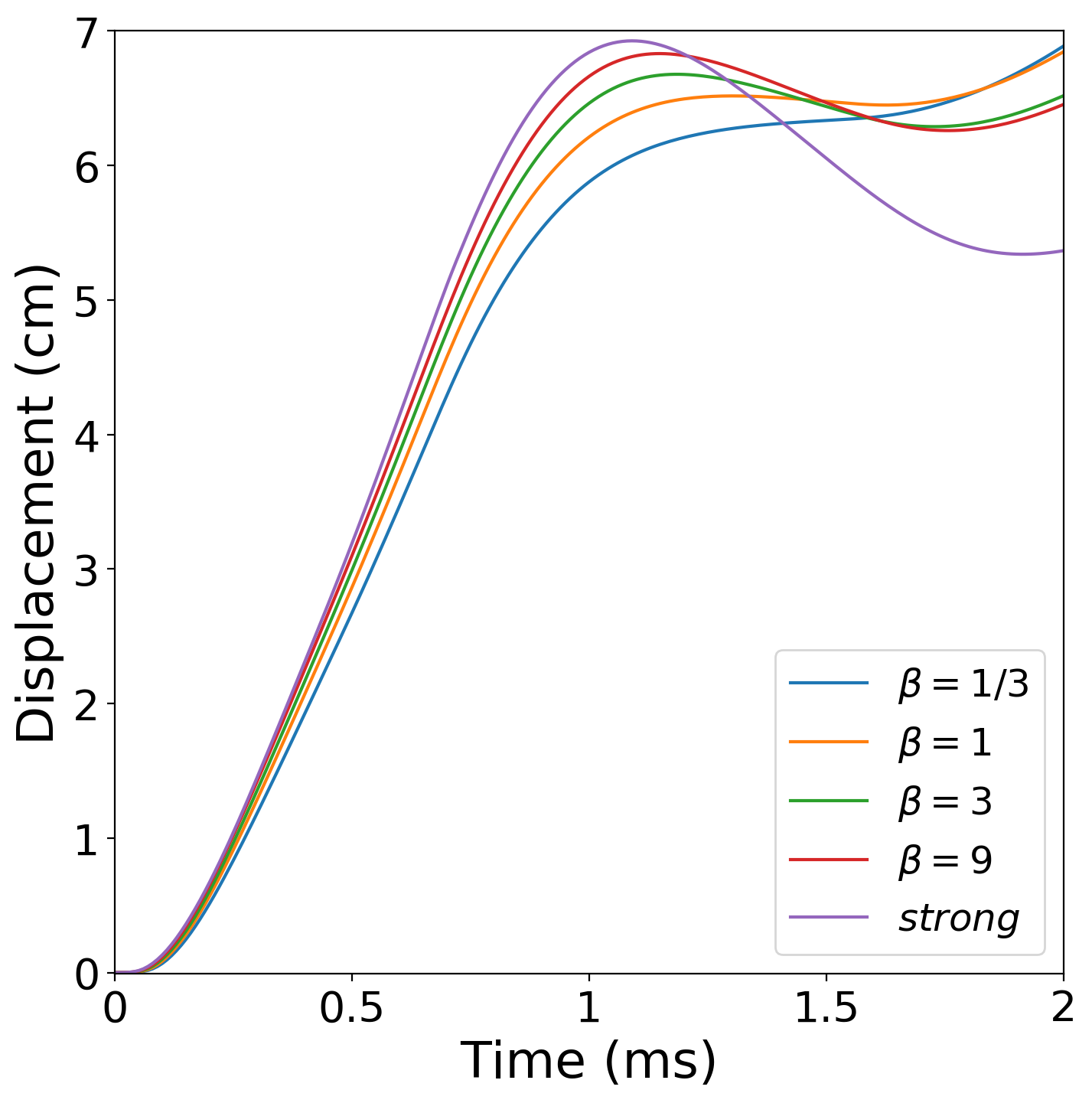}}
  \hspace{4pt}
  \subfloat{\includegraphics[width=0.303\textwidth,height=0.27\textwidth,trim={1.7cm 1.885cm 0cm 0cm},clip]{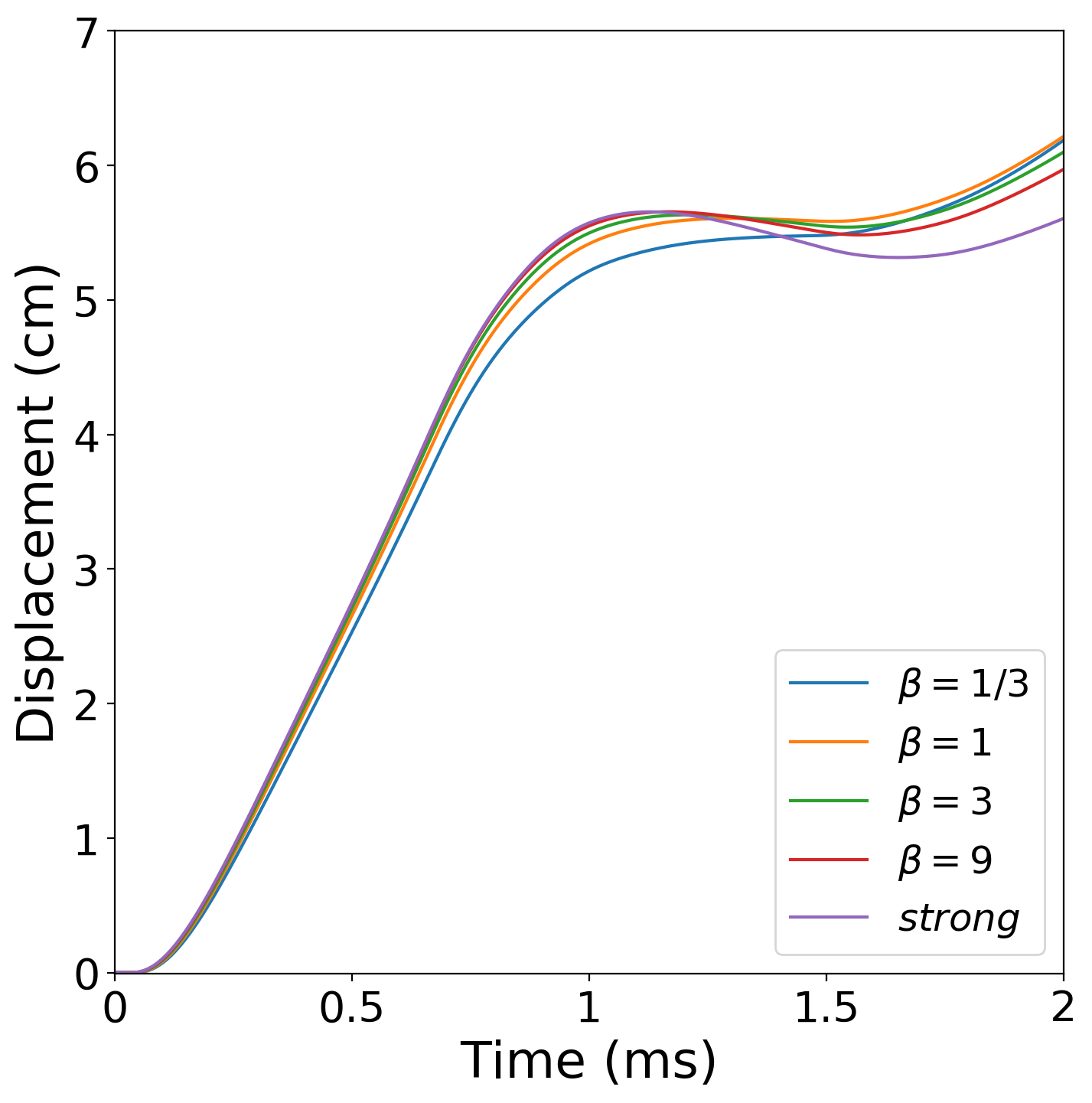}}
  \hspace{4pt}
  \subfloat{\includegraphics[width=0.303\textwidth,height=0.27\textwidth,trim={1.7cm 1.885cm 0cm 0cm},clip]{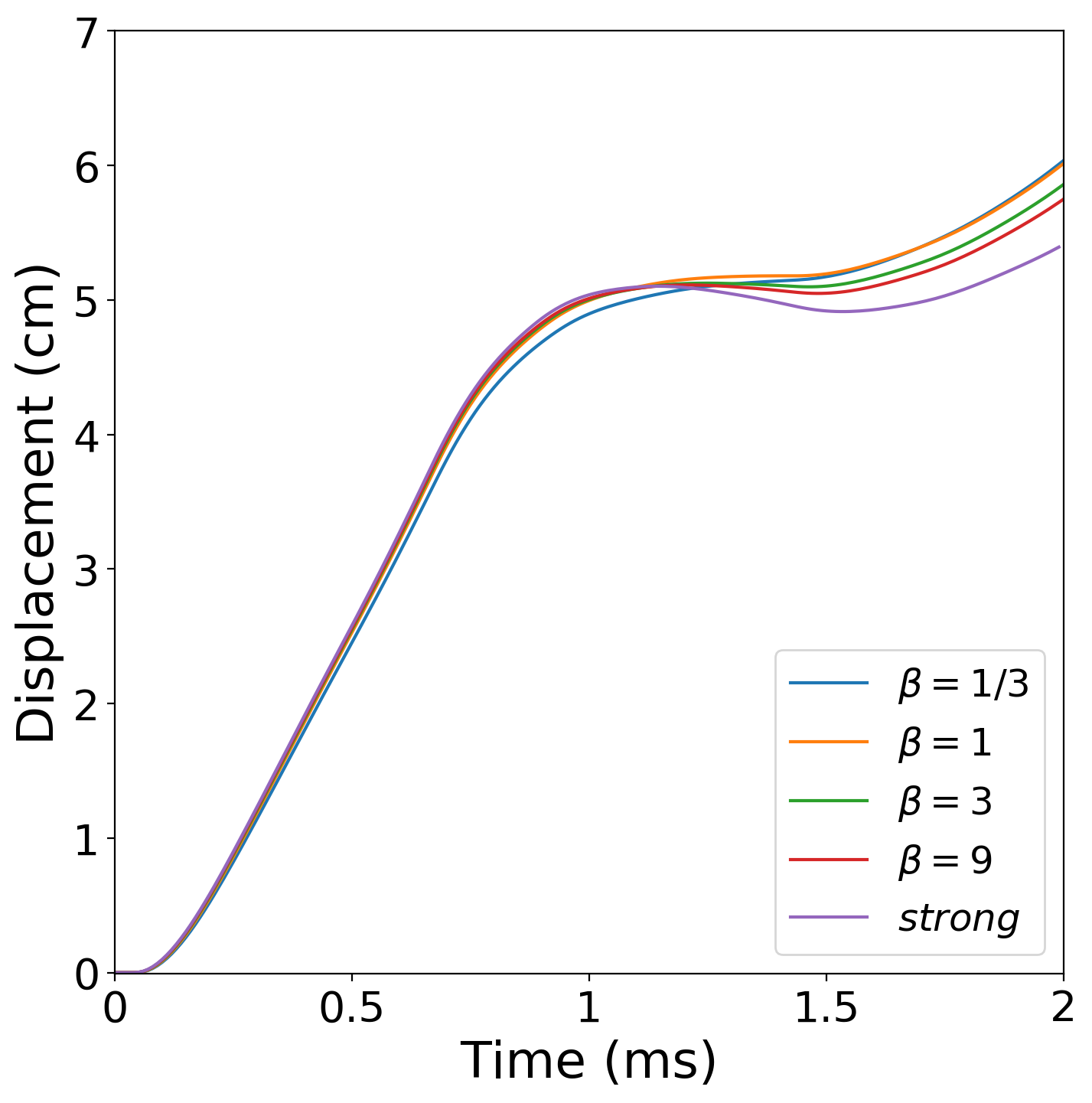}}
  
  \setcounter{subfigure}{0}
  \subfloat[][]{\includegraphics[width=0.36\textwidth,height=0.27\textwidth,trim={0cm 0cm 0cm 0cm},clip]{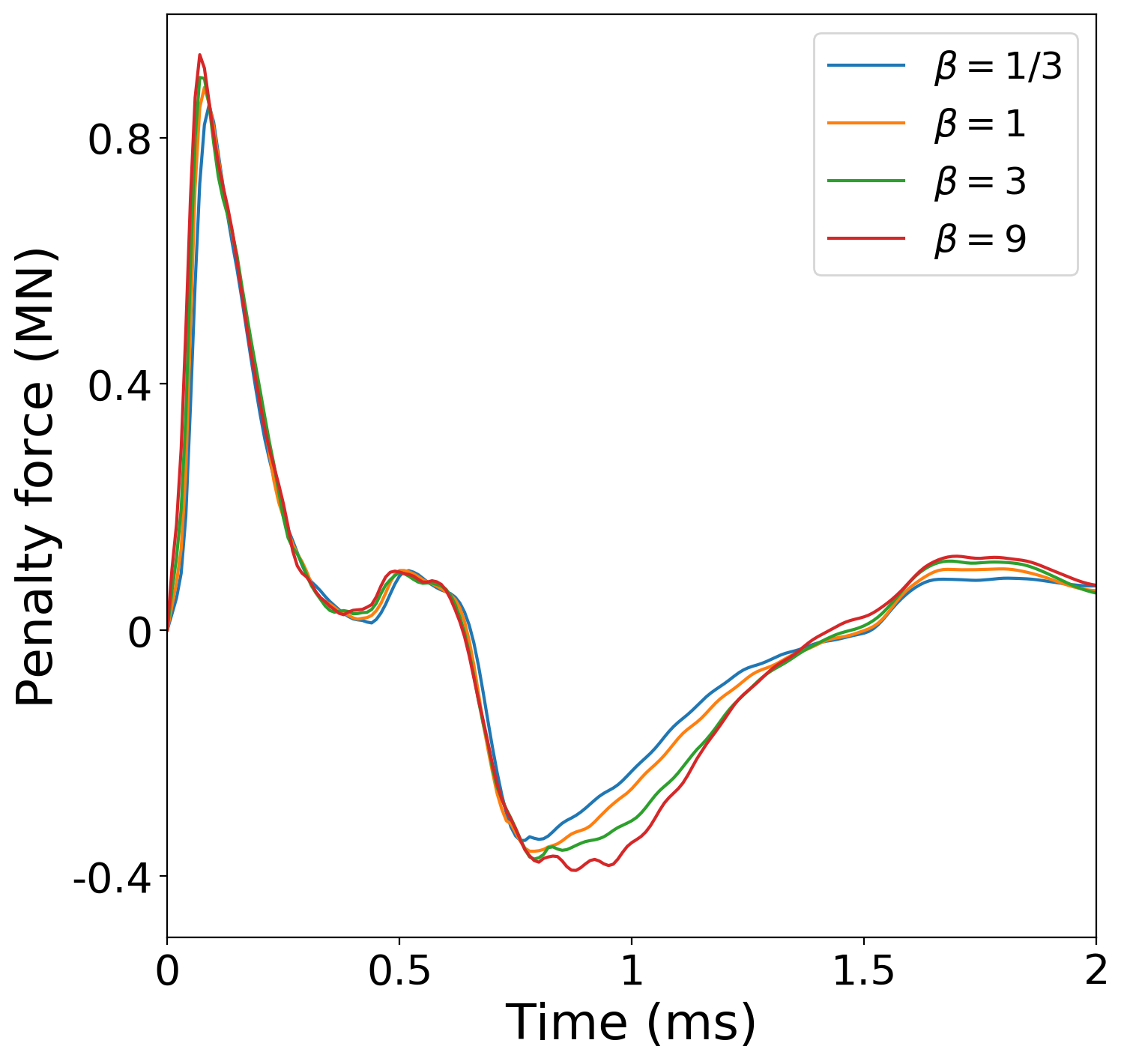}}
  \hspace{4pt}
  \subfloat[][]{\includegraphics[width=0.303\textwidth,height=0.27\textwidth,trim={2.6cm 0cm 0cm 0cm},clip]{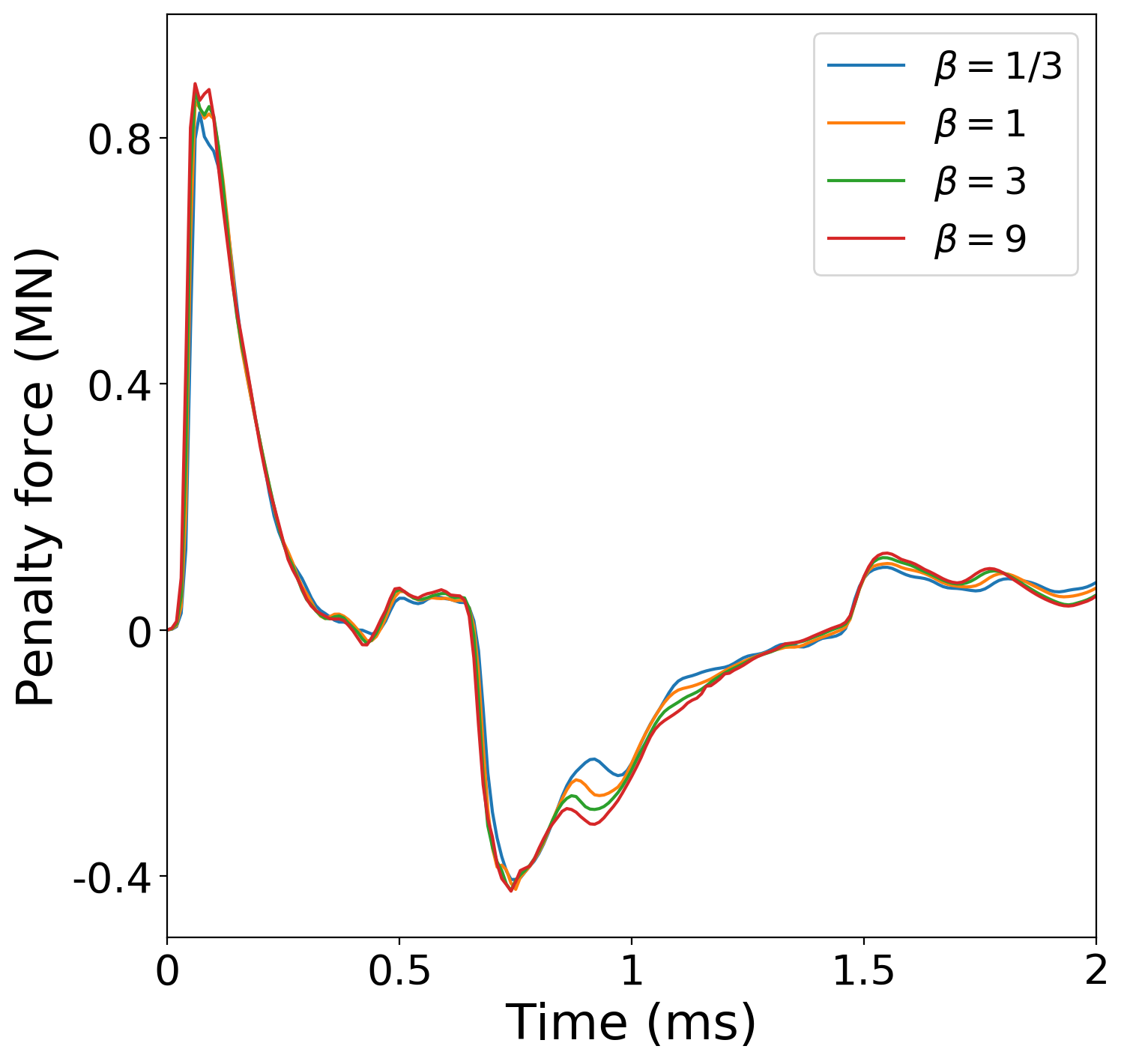}}
  \hspace{4pt}
  \subfloat[][]{\includegraphics[width=0.303\textwidth,height=0.27\textwidth,trim={2.6cm 0cm 0cm 0cm},clip]{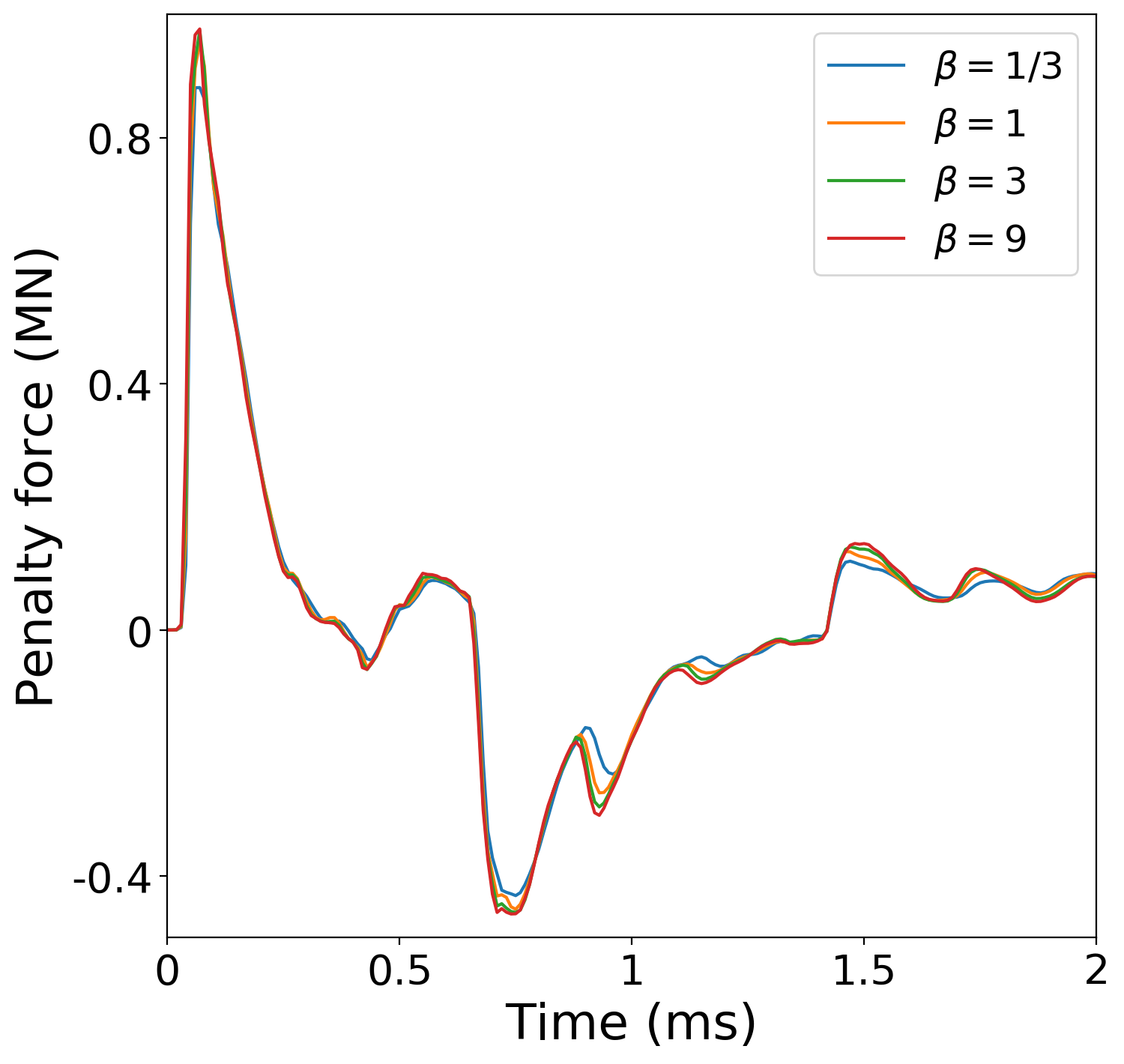}}
  \caption{Chamber detonation problem. Comparison of the results for different discretization levels and coupling approaches. (a): Coarse mesh. (b): Medium mesh. (c): Fine mesh. From top to bottom, individual rows correspond to: 1) Pressure at the center of detonation; 2) Pressure at the center of the right wall; 3) Horizontal displacement of the bar center of mass; and 4) Horizontal component of the integrated penalty force.}
  \label{fig:plastic_results}
\end{figure*}

\begin{figure*}[!hbpt]
  \centering
  \subfloat[][]{\includegraphics[width=0.24\textwidth,height=0.24\textwidth,trim={0cm 0cm 0cm 0cm},clip]{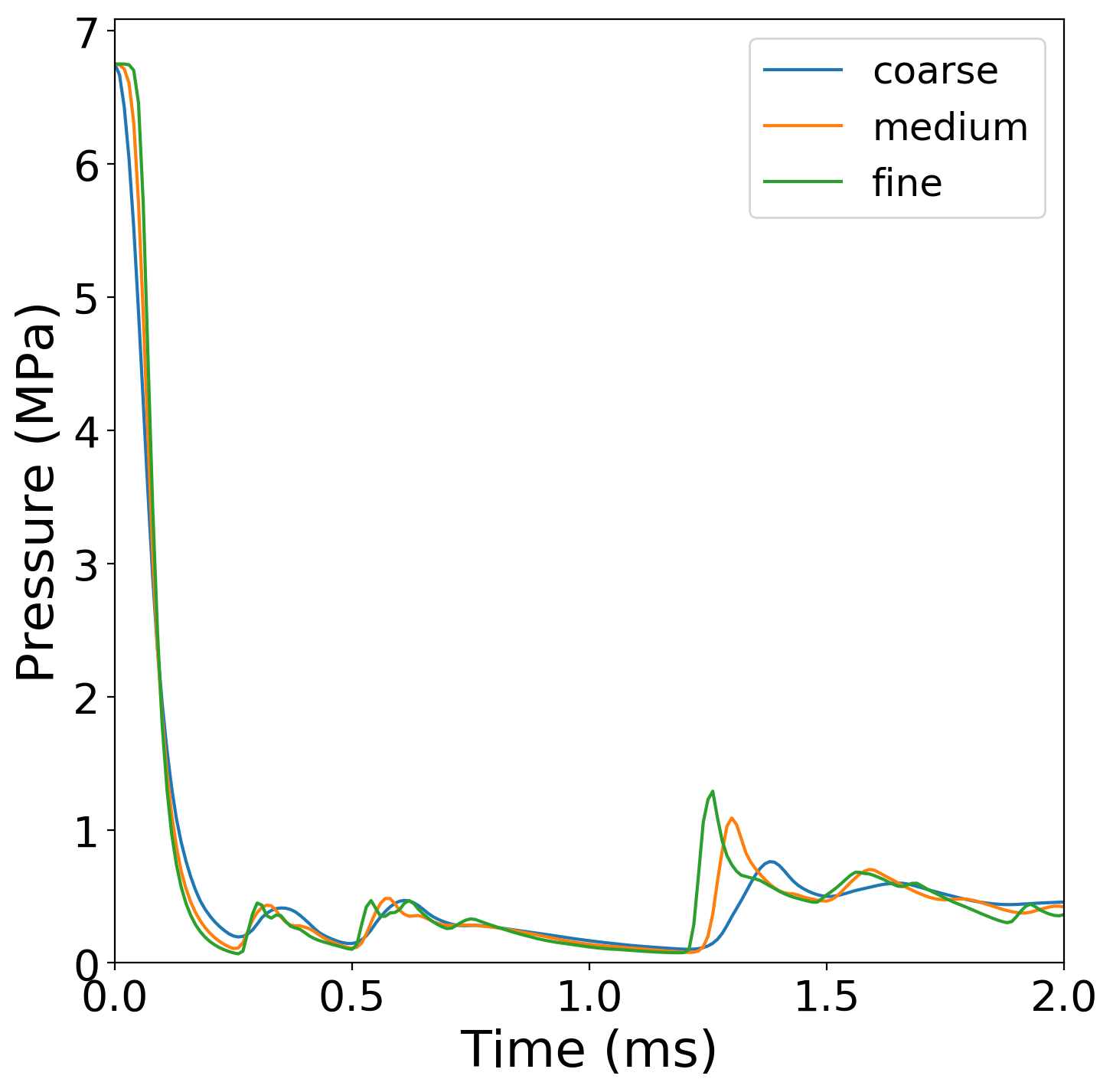}}
  \hspace{2pt}
  \subfloat[][]{\includegraphics[width=0.24\textwidth,height=0.24\textwidth,trim={0cm 0cm 0cm 0cm},clip]{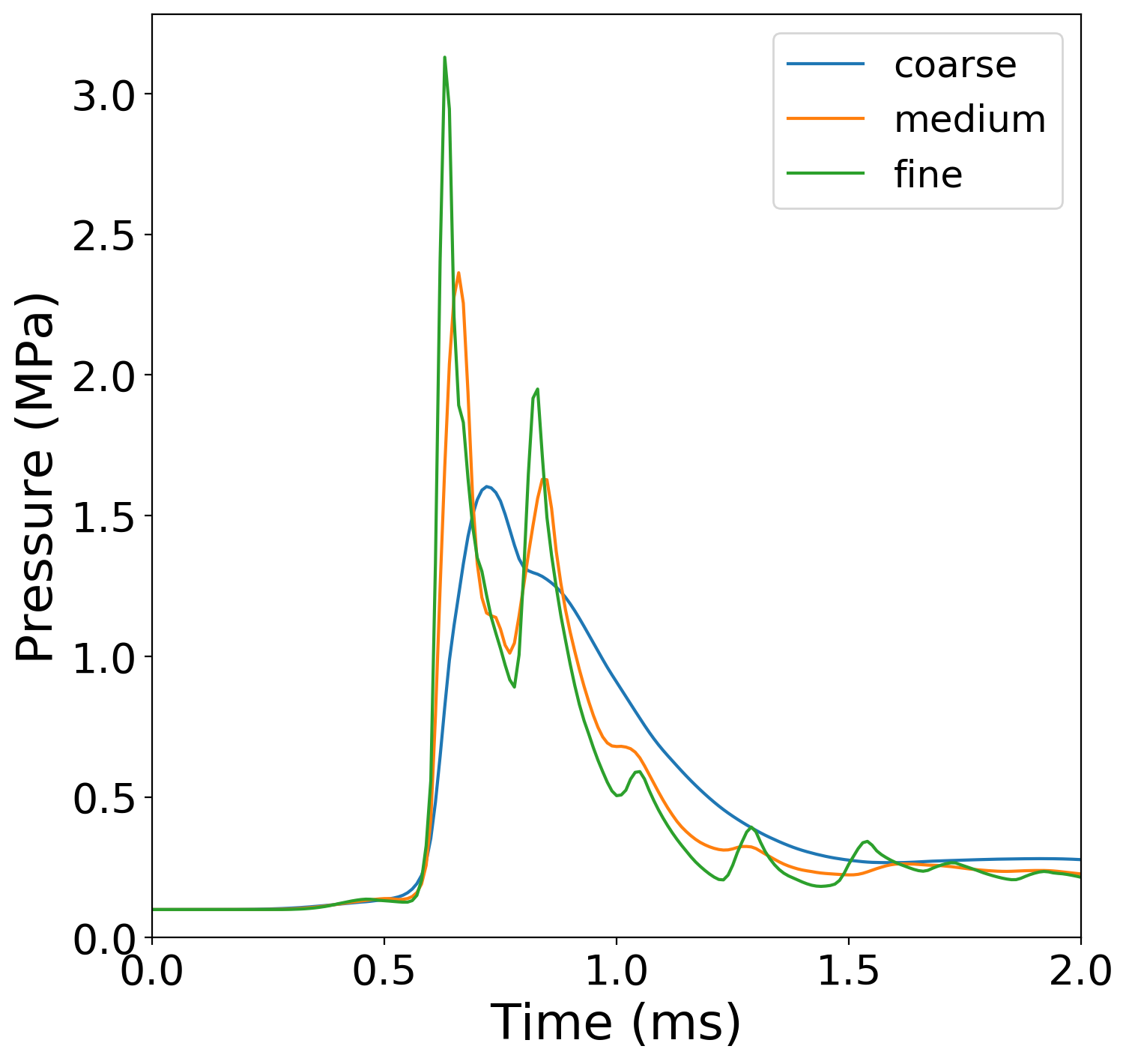}}
  \hspace{2pt}
  \subfloat[][]{\includegraphics[width=0.24\textwidth,height=0.24\textwidth,trim={0cm 0cm 0cm 0cm},clip]{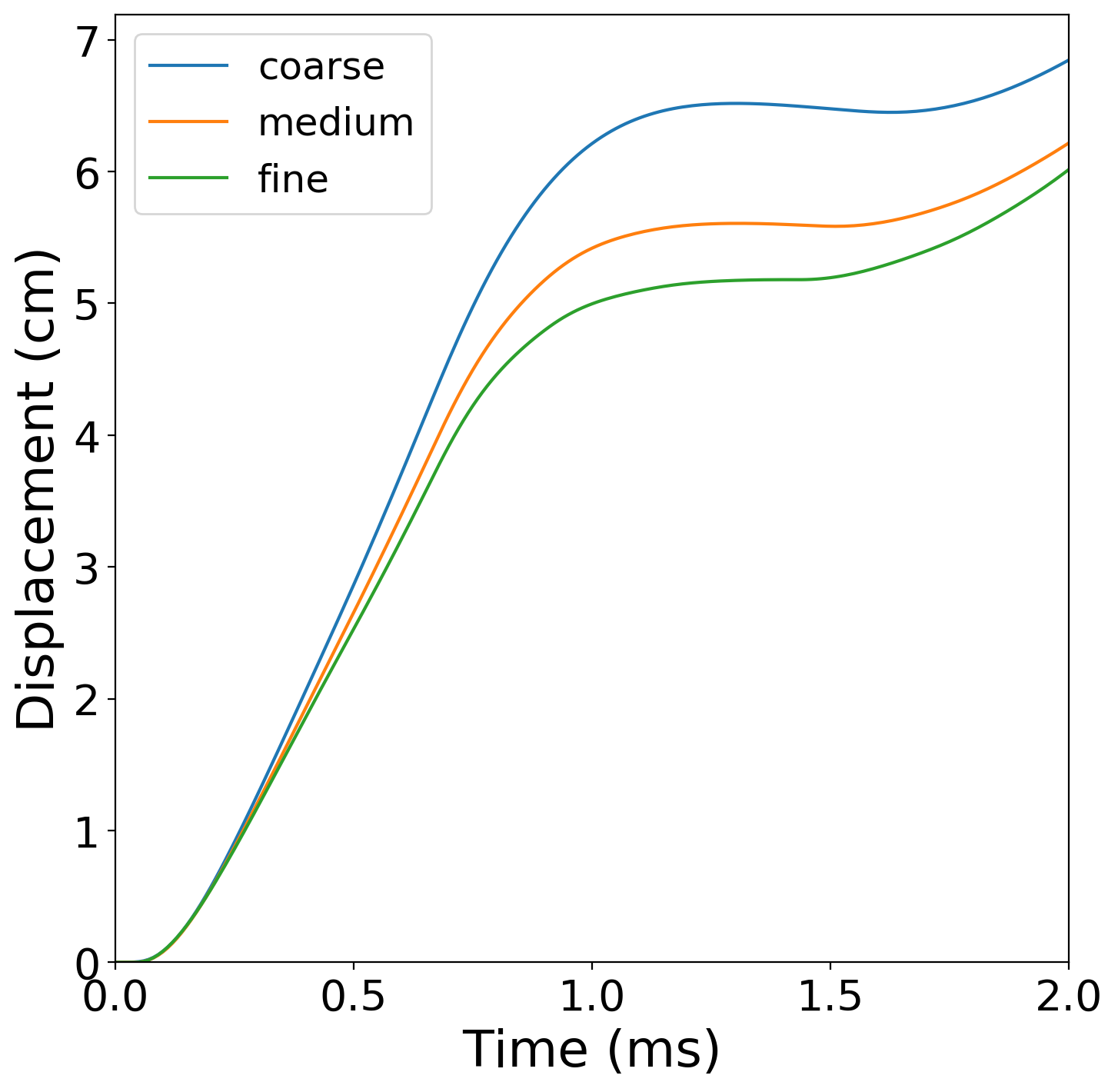}}
  \hspace{2pt}
  \subfloat[][]{\includegraphics[width=0.24\textwidth,height=0.24\textwidth,trim={0cm 0cm 0cm 0cm},clip]{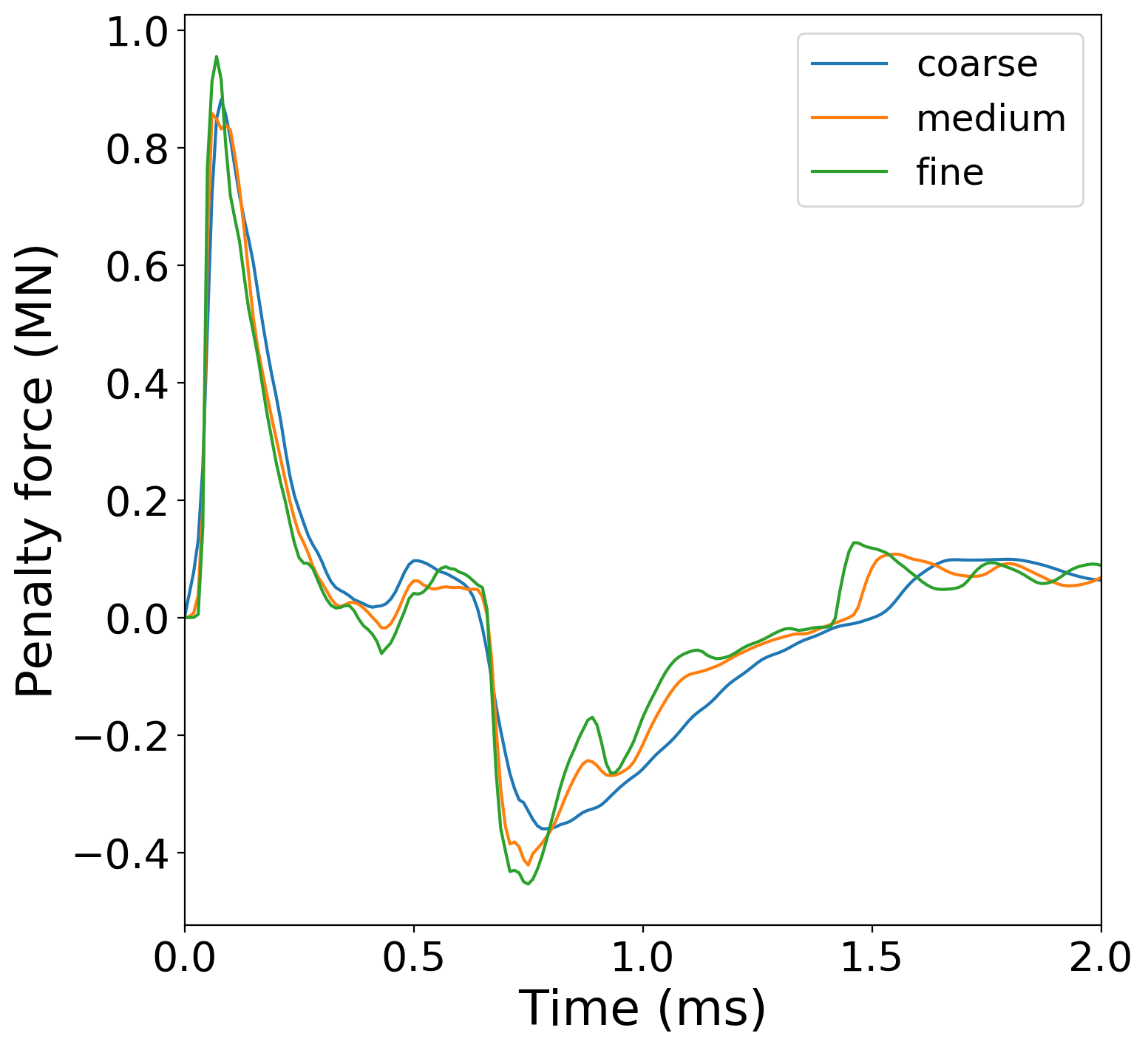}}
  \caption{Chamber detonation problem. Comparison of the results for different discretization levels of $\beta = 1$. 1) Pressure at the center of detonation; 2) Pressure at the center of the right wall; 3) Horizontal displacement of the bar center of mass; and 4) Horizontal component of the integrated penalty force.}
  \label{fig:plastic_c1}
\end{figure*}

In \cref{fig:plastic_solid}, the final shape of the bar (for the finest discretization) is compared for the weak and strong immersed FSI coupling approaches and the conforming-discretization Arbitrary Lagrangian–Eulerian (ALE) simulation results taken from~\cite{bazilevs2017new1}. The agreement with the ALE results is quite good, especially for the higher penalty-constant case. Note the fluid mesh distortion in the vicinity of the corners for the ALE case. %
\cref{fig:plastic_results} shows the time history of the pressure, the solid-object center-of-mass displacement, and the integrated penalty force on the solid in the $x-$direction for the three discretizations employed. Convergence with mesh refinement for all quantities may be inferred. In addition, the Figure shows that the penalty force is not a strong function of the penalty constant value, which is an important observation on the robustness of the approach. \cref{fig:plastic_c1} shows the same quantities for the $\beta = 1.0$ case and for the three discretizations employed. We note that all the quantities clearly converge under mesh refinement, as well as the fact that the global penalty force is likewise not a strong function of the mesh size. This suggests that even using relatively coarse background meshes may result in an accurate total force the structure feels from the surrounding fluid.   

\subsection{Ductile Solid Subjected to Internal Explosion}
\label{sec:ductile}

\begin{figure*}[!hbpt]
  \centering
  \subfloat[][]{\includegraphics[width=0.7\textwidth,trim={0cm 0cm 0cm 0cm},clip]{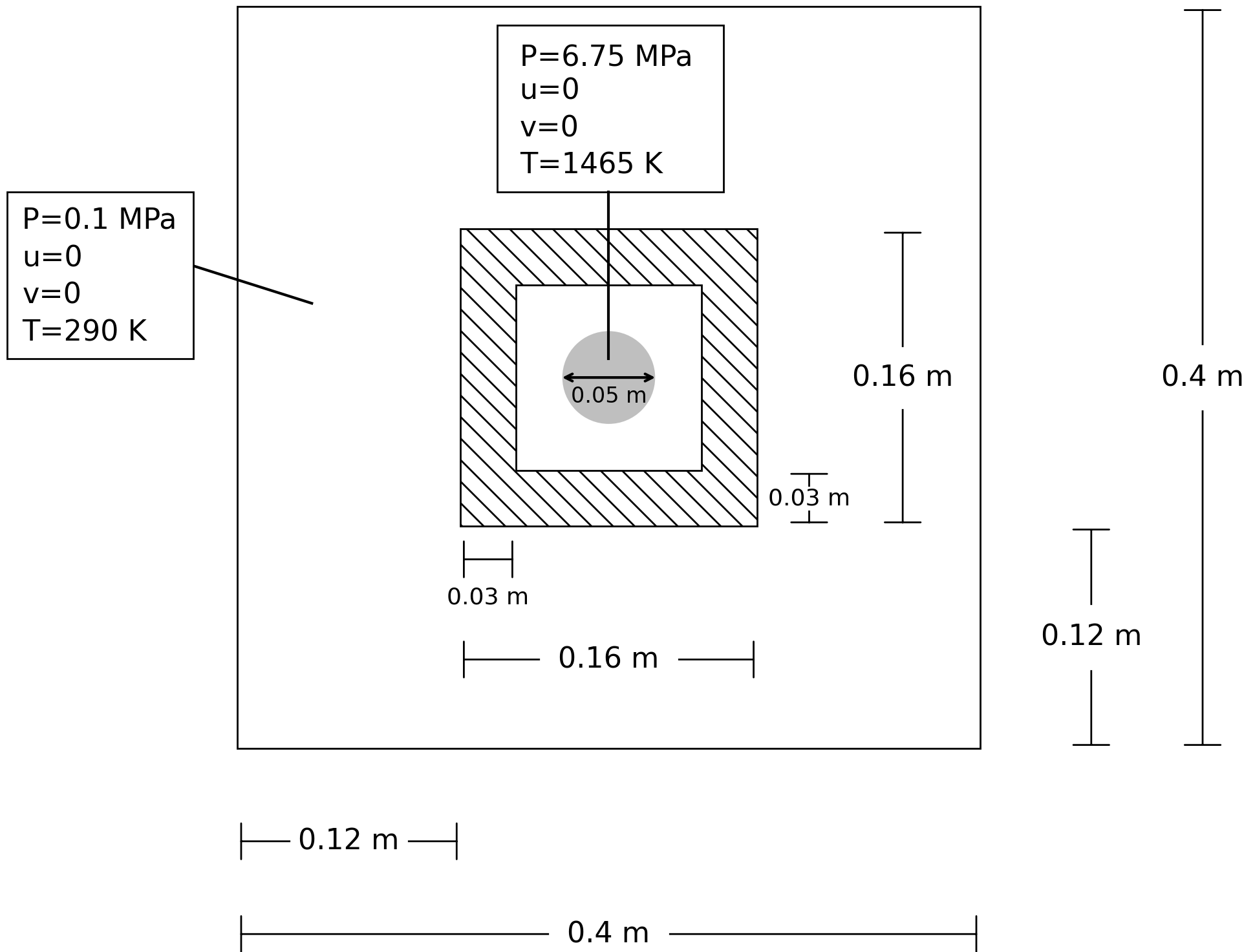}}
  \caption{Ductile fracture problem. A hollow square block of ductile material subjected to internal detonation. Problem setup and geometry.}
  \label{fig:ductile_setup}
\end{figure*}

The problem consists of a hollow square block of ductile material subjected to internal explosion. As shown in \cref{fig:ductile_setup}, the detonation is initiated at the center of the hollow square with inner dimension of $10 \, {\rm cm}$, outer dimension of $16 \, {\rm cm}$, and thickness of $3.5 \, {\rm mm}$. The background domain has dimensions of $40 \, {\rm cm} \, \times \, 40 \, {\rm cm}$. Isotropic linear hardening rule is used for the solid ductile material, which has Young's modulus ${\rm E}=200 \, {\rm GPa}$, Poisson's ratio $\nu=0.29$, yield stress $\sigma_Y=0.4 \, {\rm GPa}$, hardening modulus $H = 0.1 \, {\rm GPa}$, and initial density $\rho=7870 \, {\rm kg}/{\rm m}^3$. To simulate ductile fracture, a plasticity-driven failure approach described in~\cite{behzadinasab2020semi,behzadinasab2020peridynamic} is employed with $\bar{\epsilon}^P_{\rm th} = 0.18$ and $\bar{\epsilon}^P_{\rm cr} = 0.2$. The air is at rest initially with $p = 0.1 \, {\rm MPa}$ and $T = 290 \, {\rm K}$. The detonation condition is initiated by setting the pressure $p = 6.75 \, {\rm MPa}$ and temperature $T = 1465 \, {\rm K}$ in a circular area with a radius of $5 \, {\rm cm}$ centered inside the hollow square. 
In this example, the background domain and foreground solid are discretized uniformly. Three discretization levels are considered here, with the solid node spacing of $h = 2 \, {\rm mm}$, $1.5 \, {\rm mm}$, and $1 \, {\rm mm}$, respectively. In each case, the fluid mesh size is set to four times that of the solid node spacing. The time step size used for the coarse, medium, and fine strongly coupled cases is $\SI{0.4}{\micro s}$, $\SI{0.3}{\micro s}$, and $\SI{0.2}{\micro s}$, respectively. The time step used for the weakly coupled cases is four times smaller for the strongly coupled cases on the respective meshes. This choice is made to stably accommodate a penalty constant of $\beta = 1$ employed in the weakly coupled simulations. For the weakly coupled cases, we also examine the effect of adding the damage variable in the definition of the penalty parameter (see \cref{sec:discFSI}).

\begin{figure*}[!hbpt]
  \centering
  \subfloat{\includegraphics[width=0.6\textwidth,trim={0cm 0cm 0cm 4cm},clip]{air_vel_scale.png}}
  
  \subfloat{\includegraphics[width=0.235\textwidth,trim={0cm 0cm 0cm 0cm},clip]{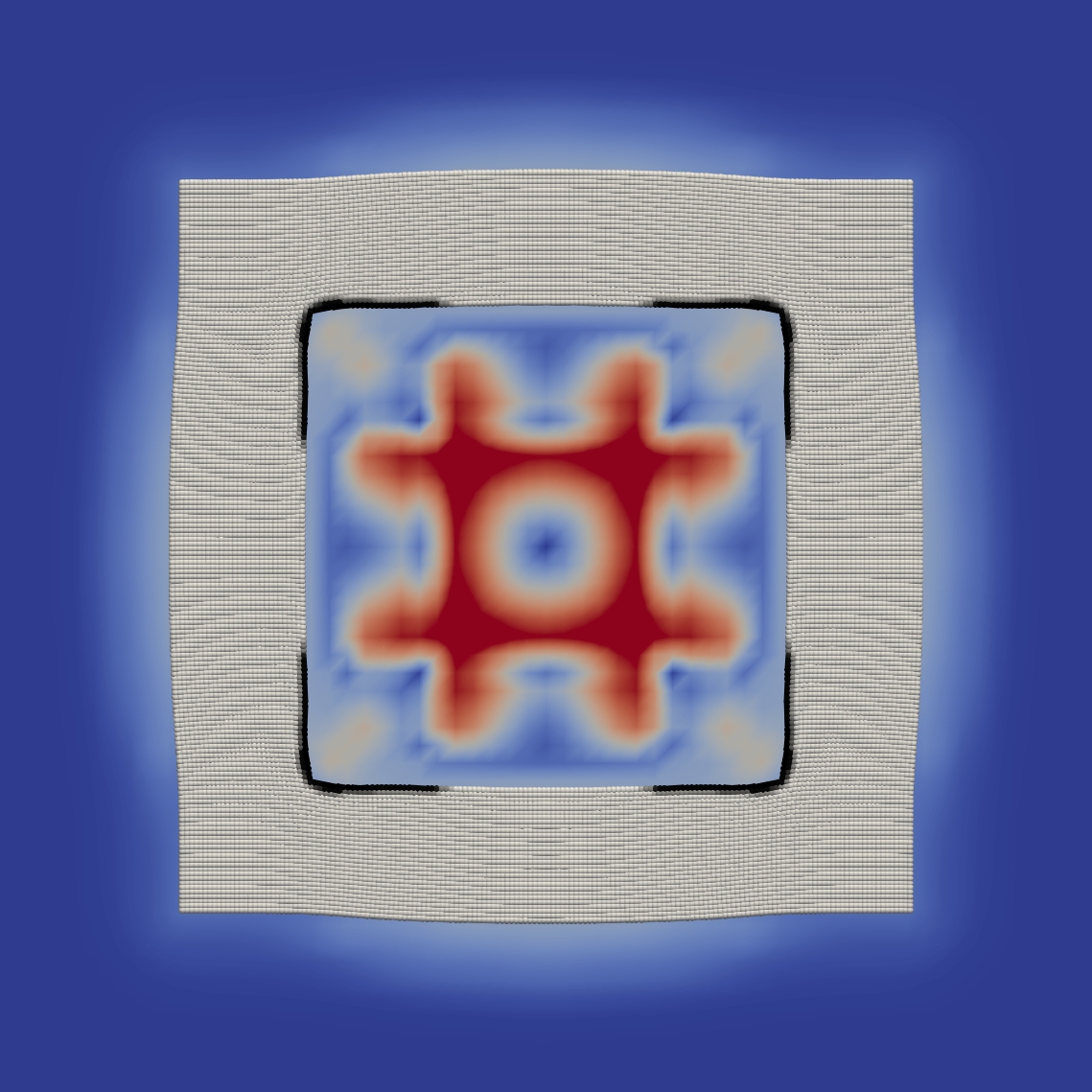}}
  \hspace{5pt}
  \subfloat{\includegraphics[width=0.235\textwidth,trim={0cm 0cm 0cm 0cm},clip]{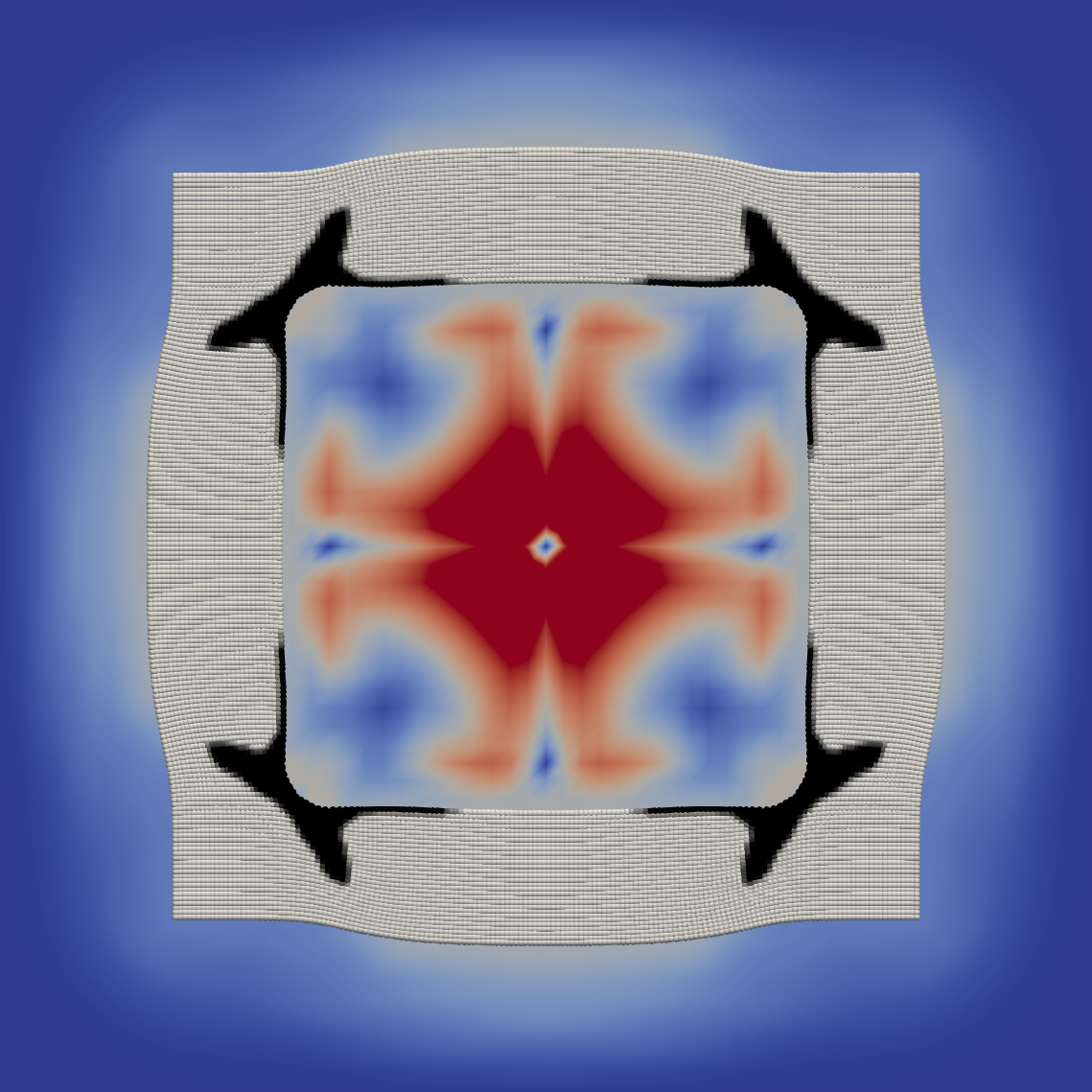}}
  \hspace{5pt}
  \subfloat{\includegraphics[width=0.235\textwidth,trim={0cm 0cm 0cm 0cm},clip]{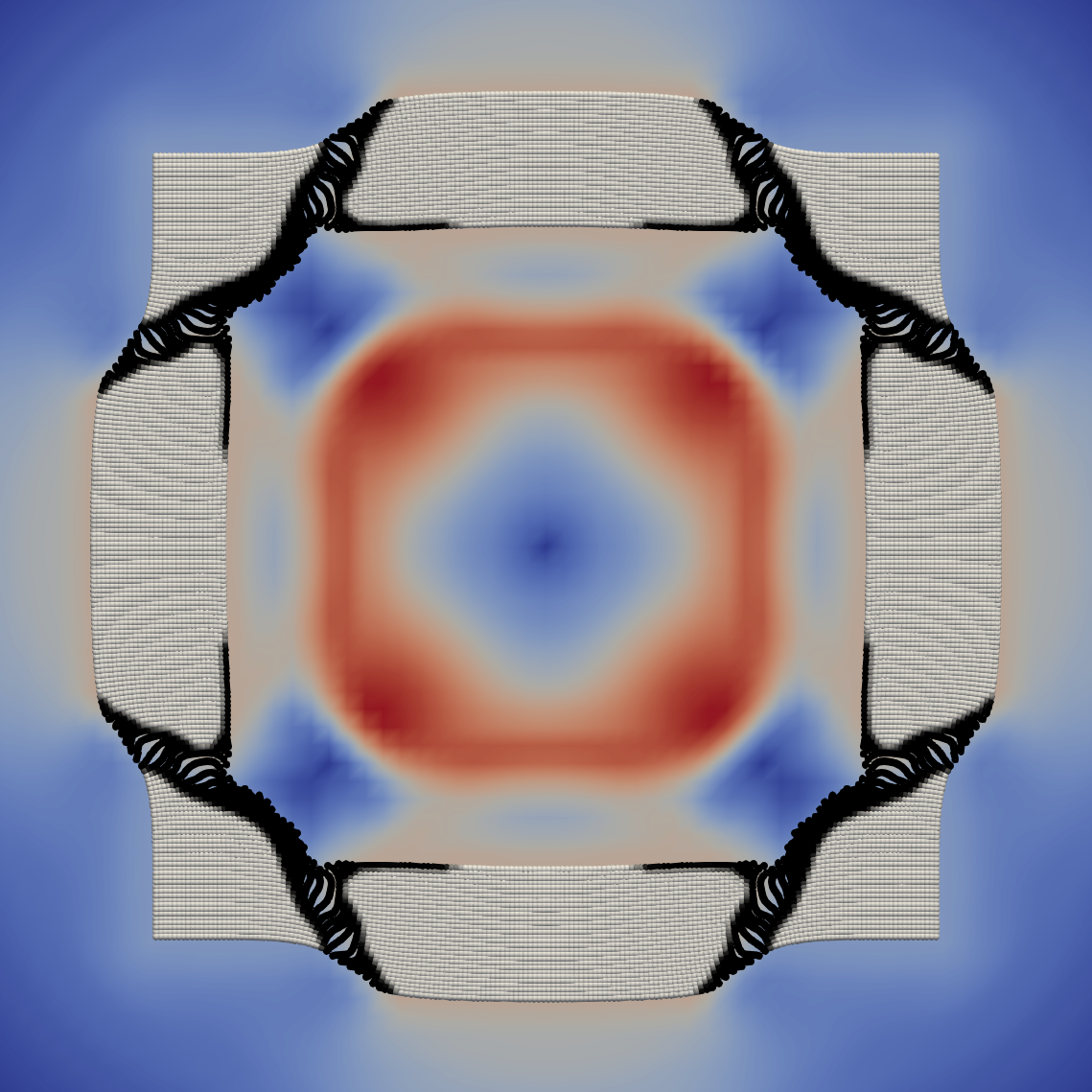}}
  \hspace{5pt}
  \subfloat{\includegraphics[width=0.235\textwidth,trim={0cm 0cm 0cm 0cm},clip]{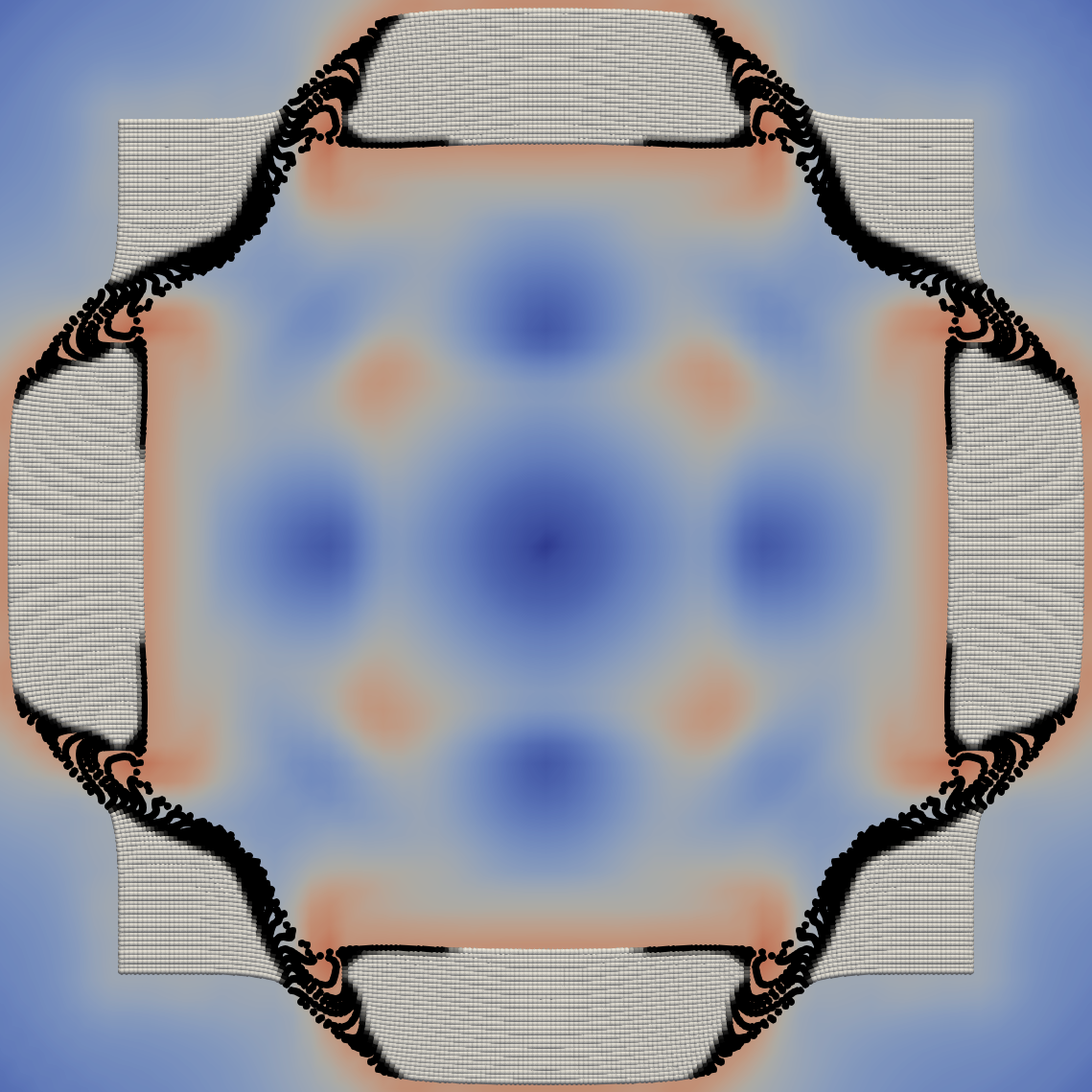}}
  
  \subfloat{\includegraphics[width=0.235\textwidth,trim={0cm 0cm 0cm 0cm},clip]{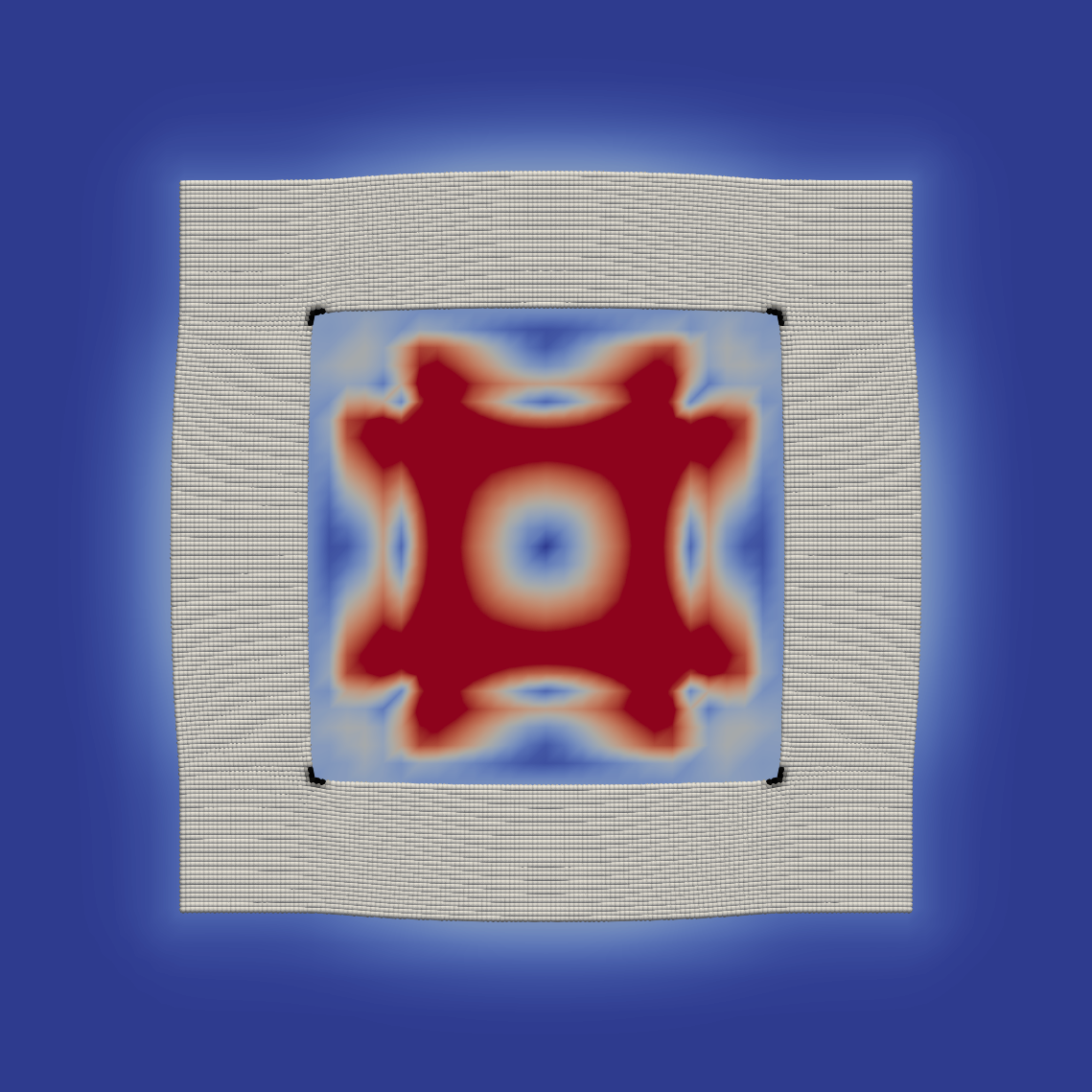}}
  \hspace{5pt}
  \subfloat{\includegraphics[width=0.235\textwidth,trim={0cm 0cm 0cm 0cm},clip]{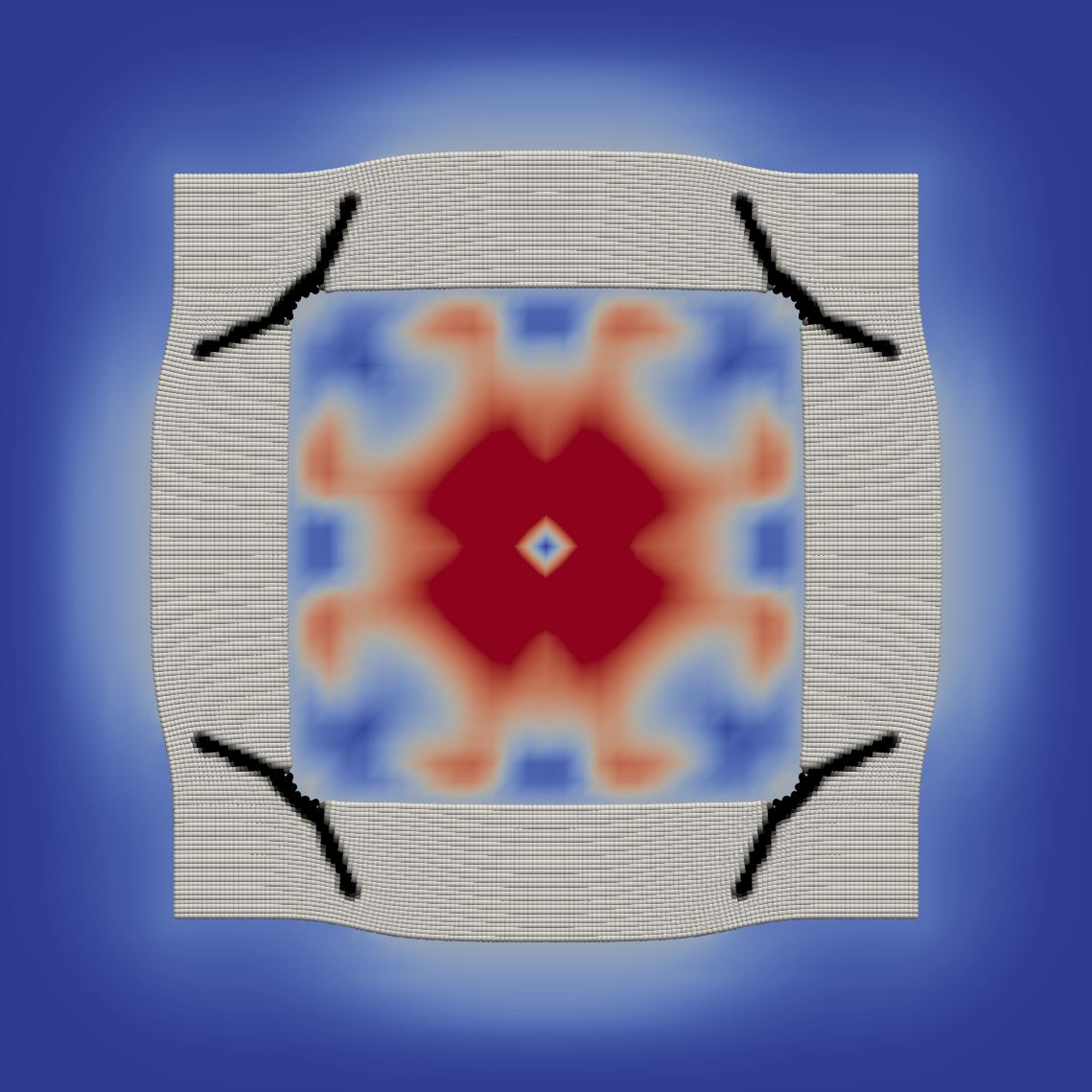}}
  \hspace{5pt}
  \subfloat{\includegraphics[width=0.235\textwidth,trim={0cm 0cm 0cm 0cm},clip]{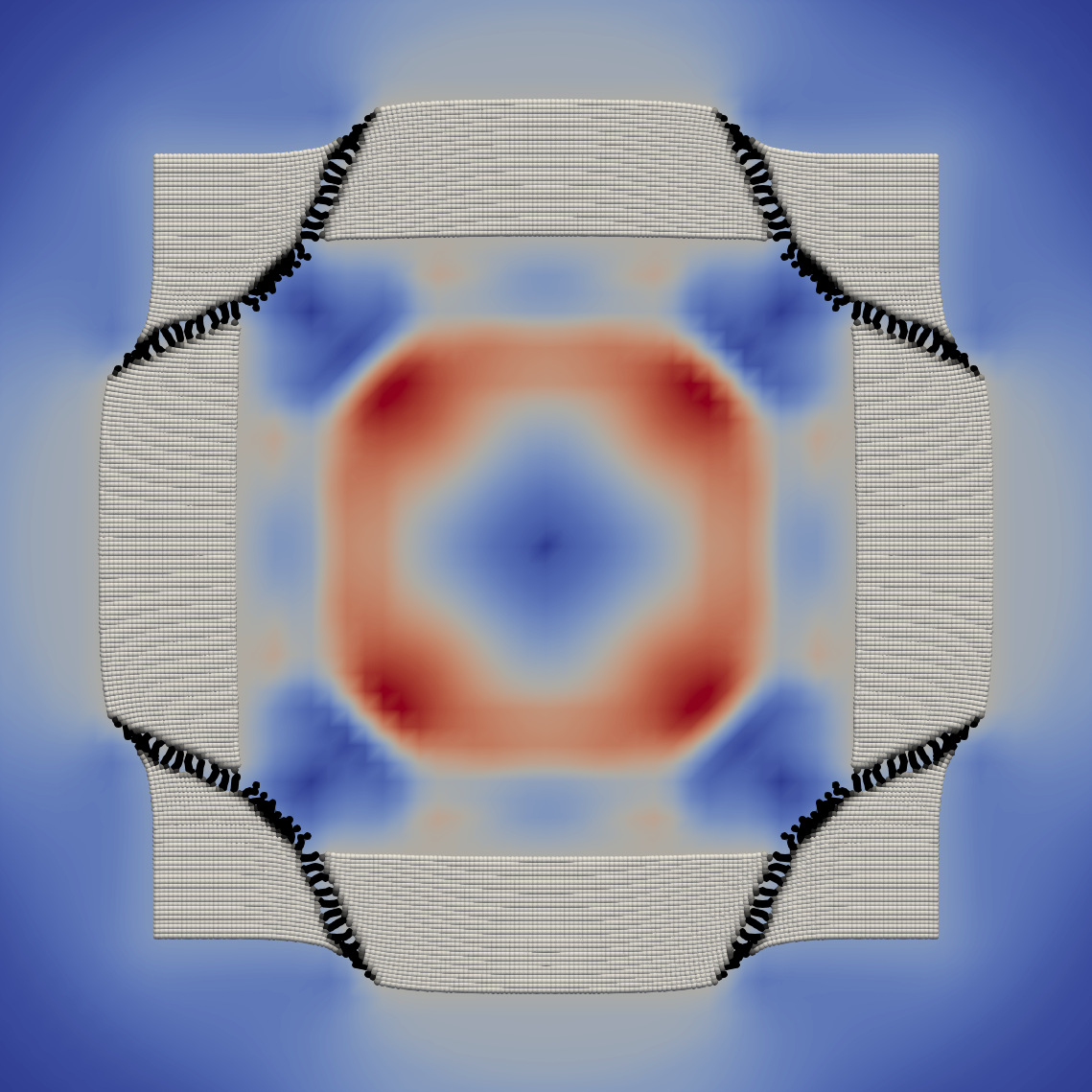}}
  \hspace{5pt}
  \subfloat{\includegraphics[width=0.235\textwidth,trim={0cm 0cm 0cm 0cm},clip]{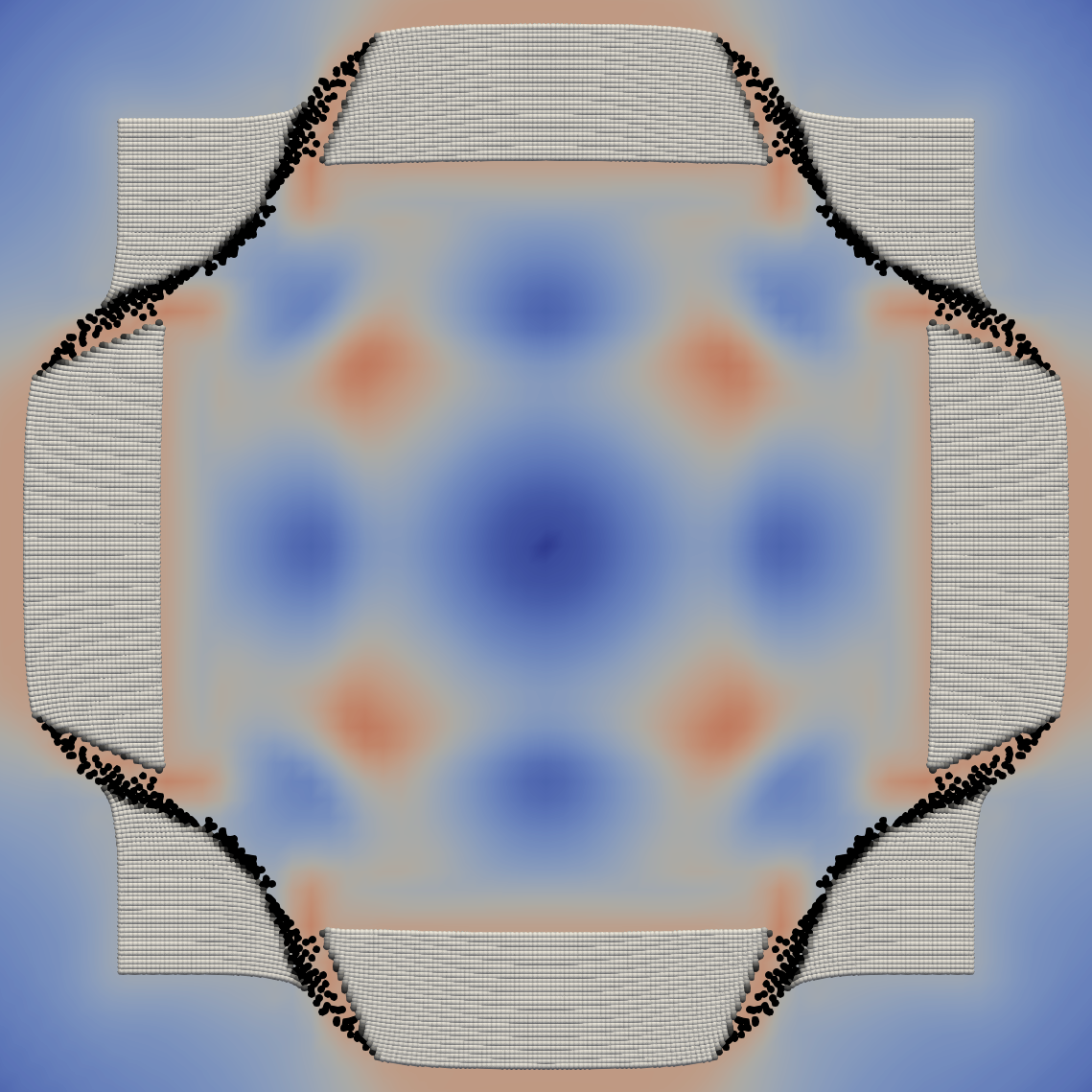}}

  \setcounter{subfigure}{0}
  \subfloat[][\SI{80}{\micro s}]{\includegraphics[width=0.235\textwidth,trim={0cm 0cm 0cm 0cm},clip]{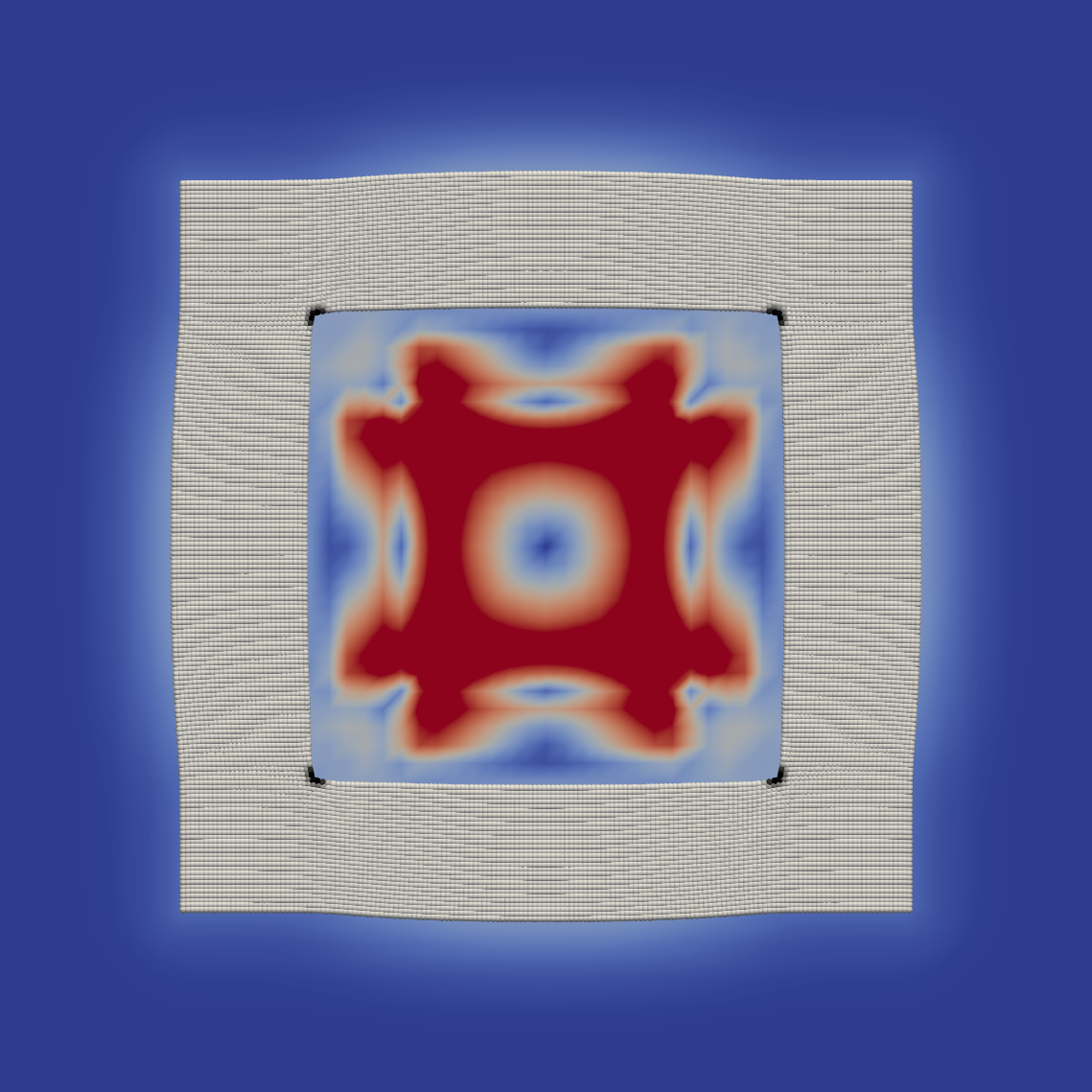}}
  \hspace{5pt}
  \subfloat[][\SI{120}{\micro s}]{\includegraphics[width=0.235\textwidth,trim={0cm 0cm 0cm 0cm},clip]{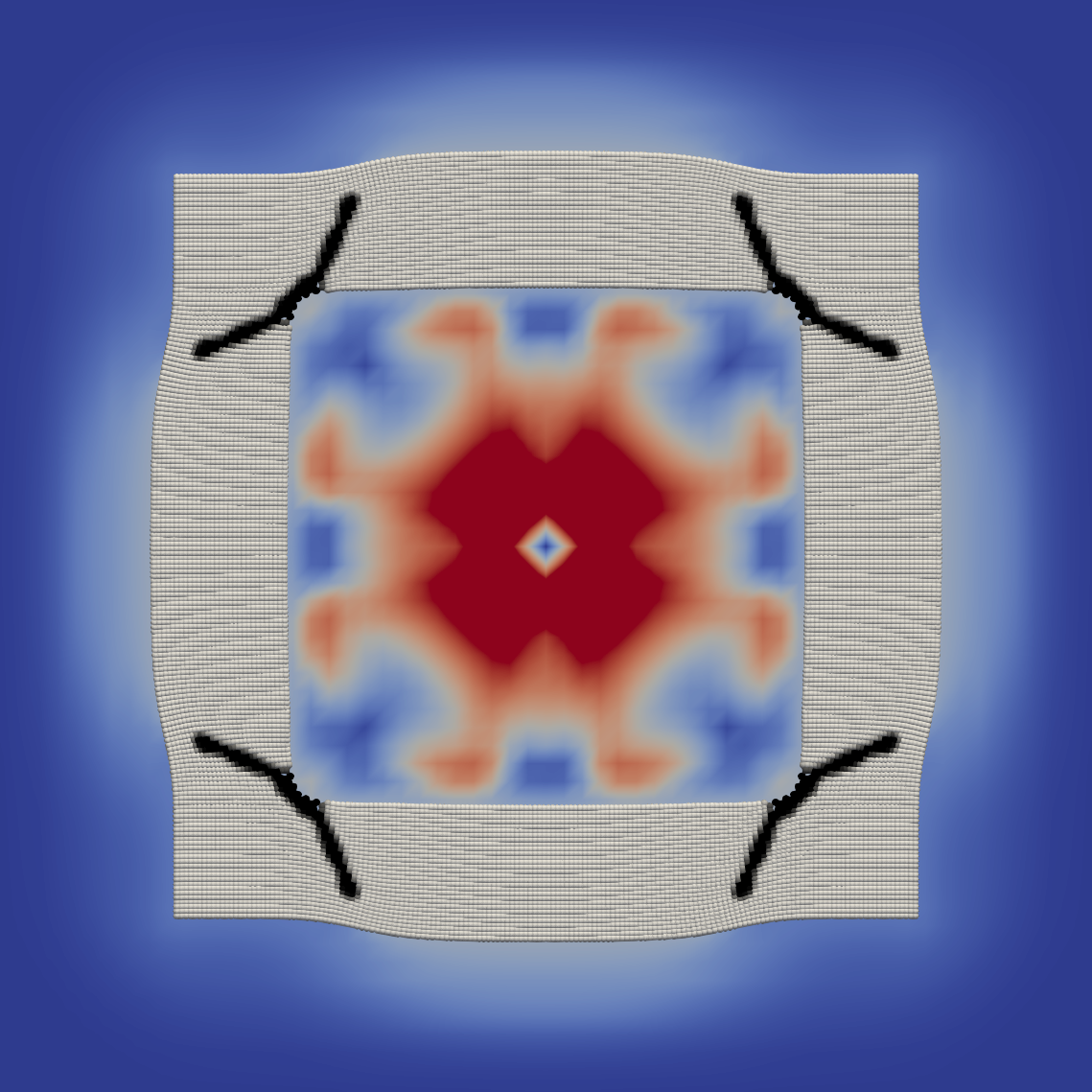}}
  \hspace{5pt}
  \subfloat[][\SI{200}{\micro s}]{\includegraphics[width=0.235\textwidth,trim={0cm 0cm 0cm 0cm},clip]{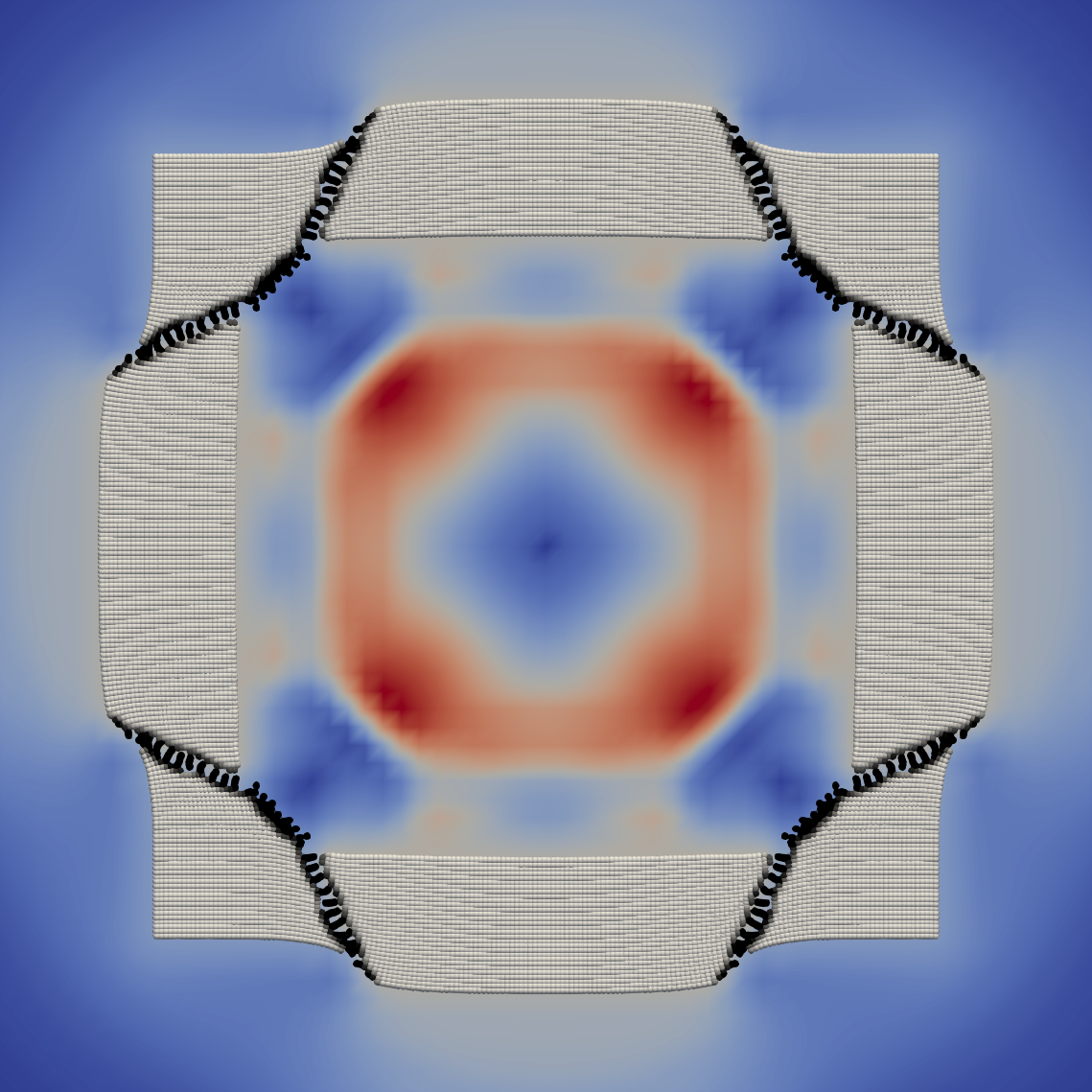}}
  \hspace{5pt}
  \subfloat[][\SI{300}{\micro s}]{\includegraphics[width=0.235\textwidth,trim={0cm 0cm 0cm 0cm},clip]{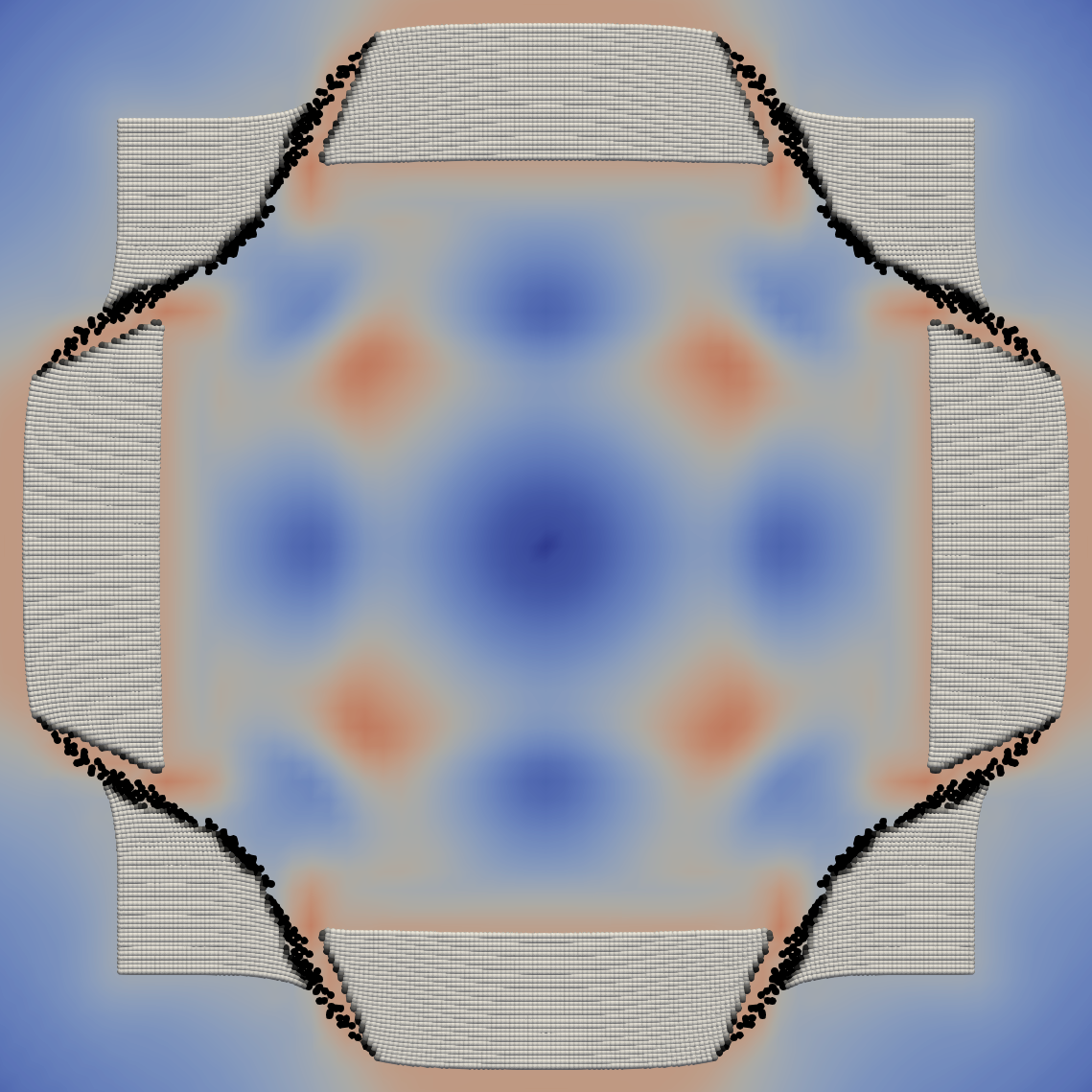}}

  \caption{Ductile fracture problem. Snapshots of air speed (in m/s) and solid damage in the current configuration at different stages during the simulation for the finest mesh. Top, middle, and bottom rows correspond to strong coupling, weak coupling without damage in the penalty stiffness, and weak coupling with damage in the penalty stiffness, respectively.}
  \label{fig:ductile_contours}
\end{figure*}

\begin{figure*}[!hbpt]
  \centering
  \subfloat{\includegraphics[width=0.32\textwidth,height=0.31\textwidth,trim={0cm 0cm 0cm 0cm},clip]{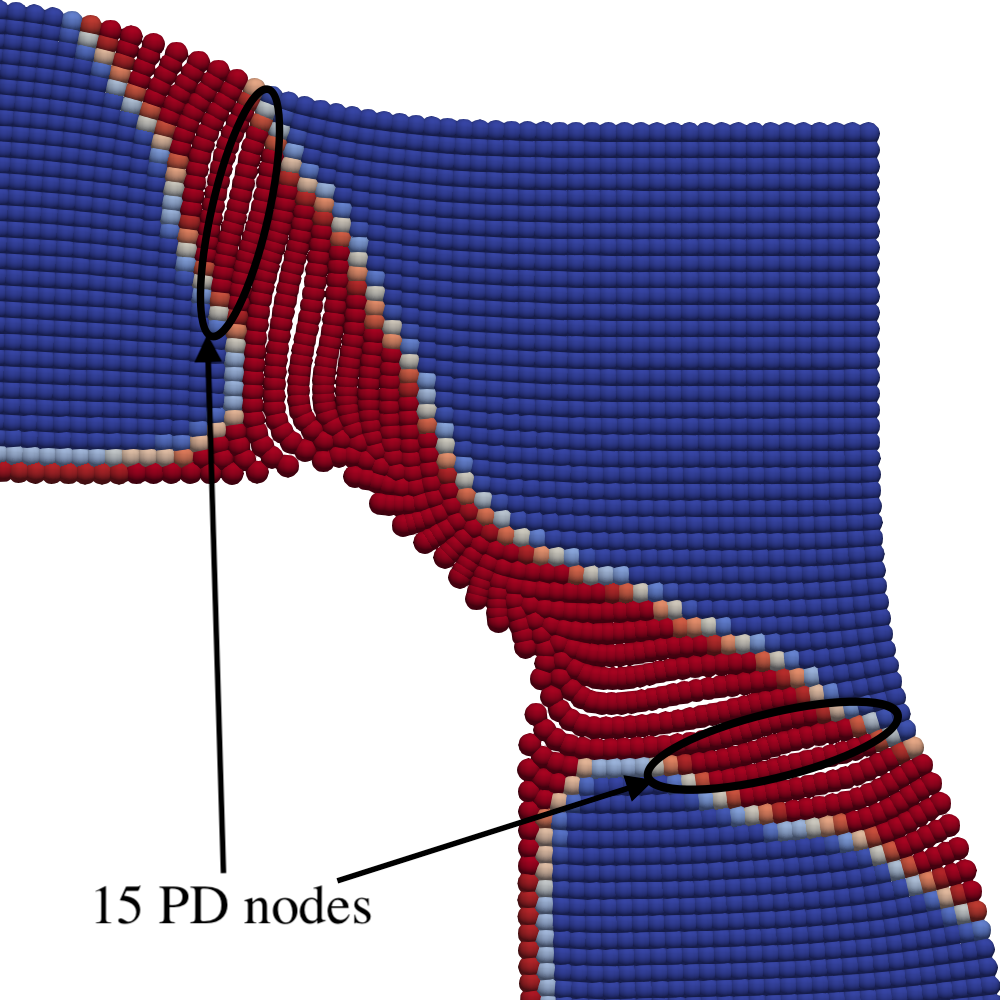}}
  \hspace{6pt}
  \subfloat{\includegraphics[width=0.32\textwidth,height=0.31\textwidth,trim={0cm 0cm 0cm 0cm},clip]{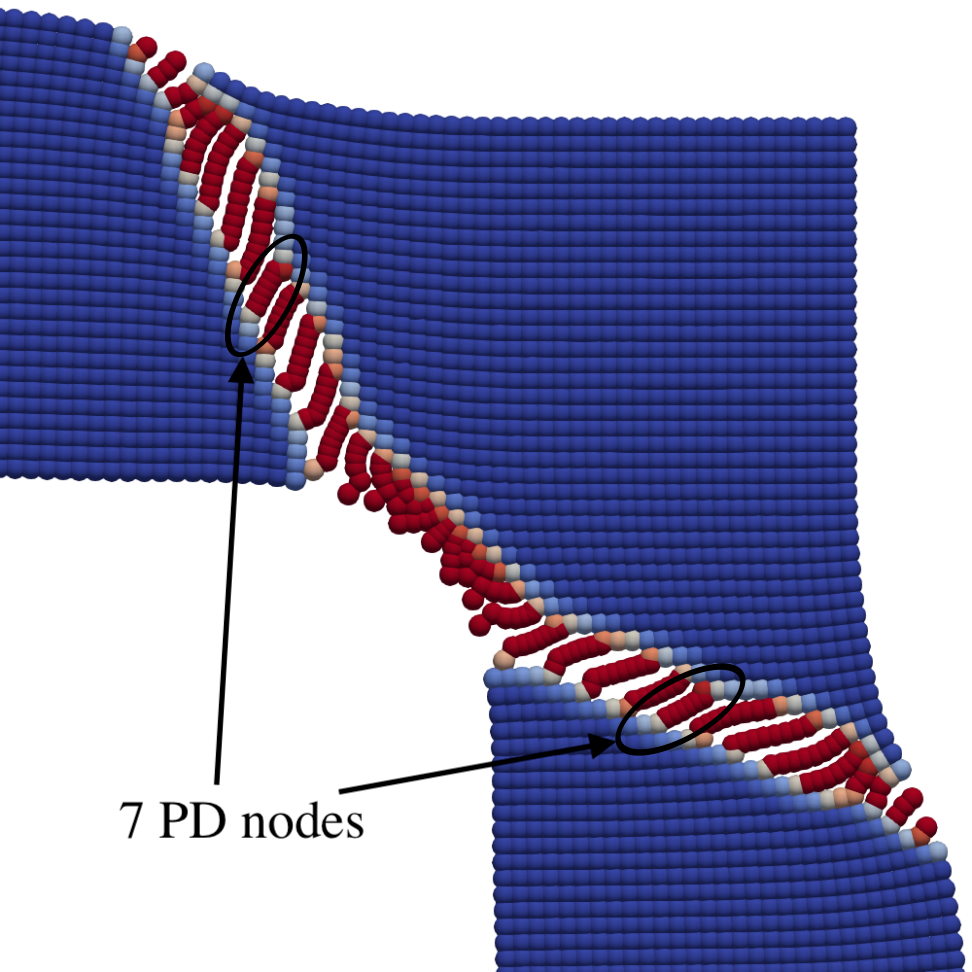}}
  \hspace{6pt}
  \subfloat{\includegraphics[width=0.32\textwidth,height=0.31\textwidth,trim={0cm 0cm 0cm 0cm},clip]{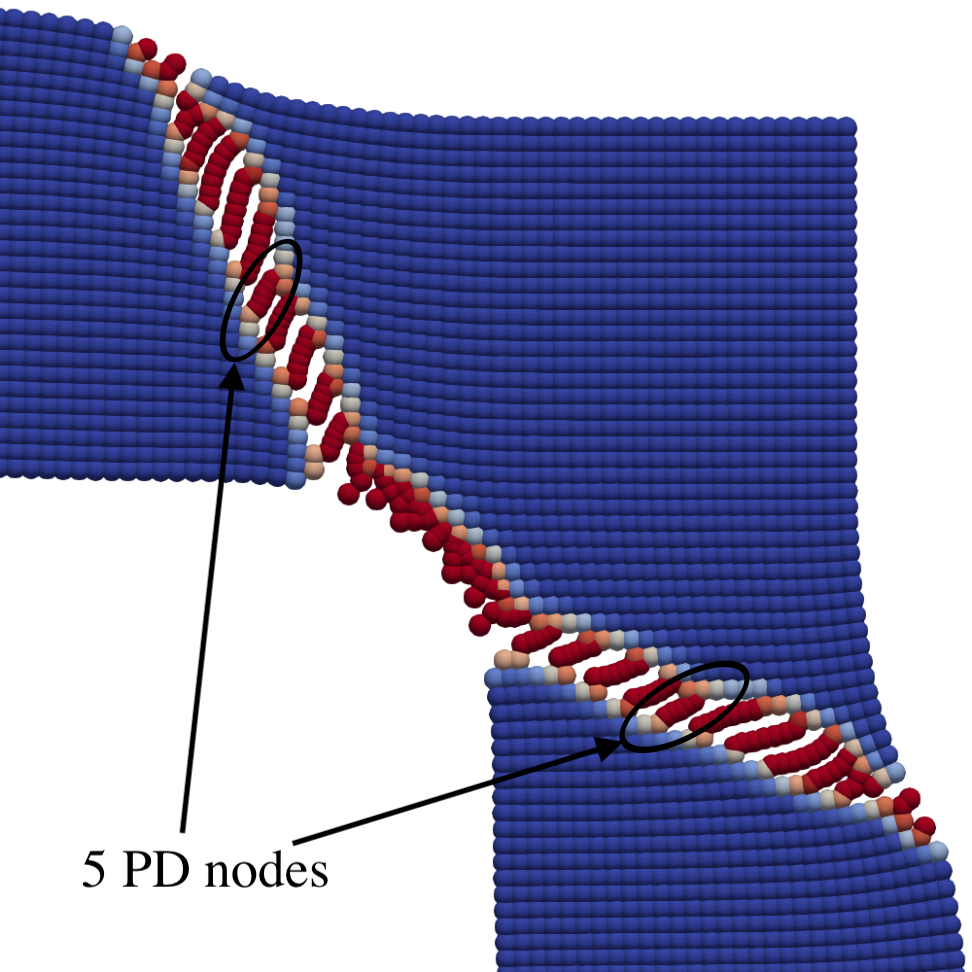}}
  
  \setcounter{subfigure}{0}
  \subfloat[][]{\includegraphics[width=0.32\textwidth,height=0.31\textwidth,trim={0cm 0cm 0cm 0cm},clip]{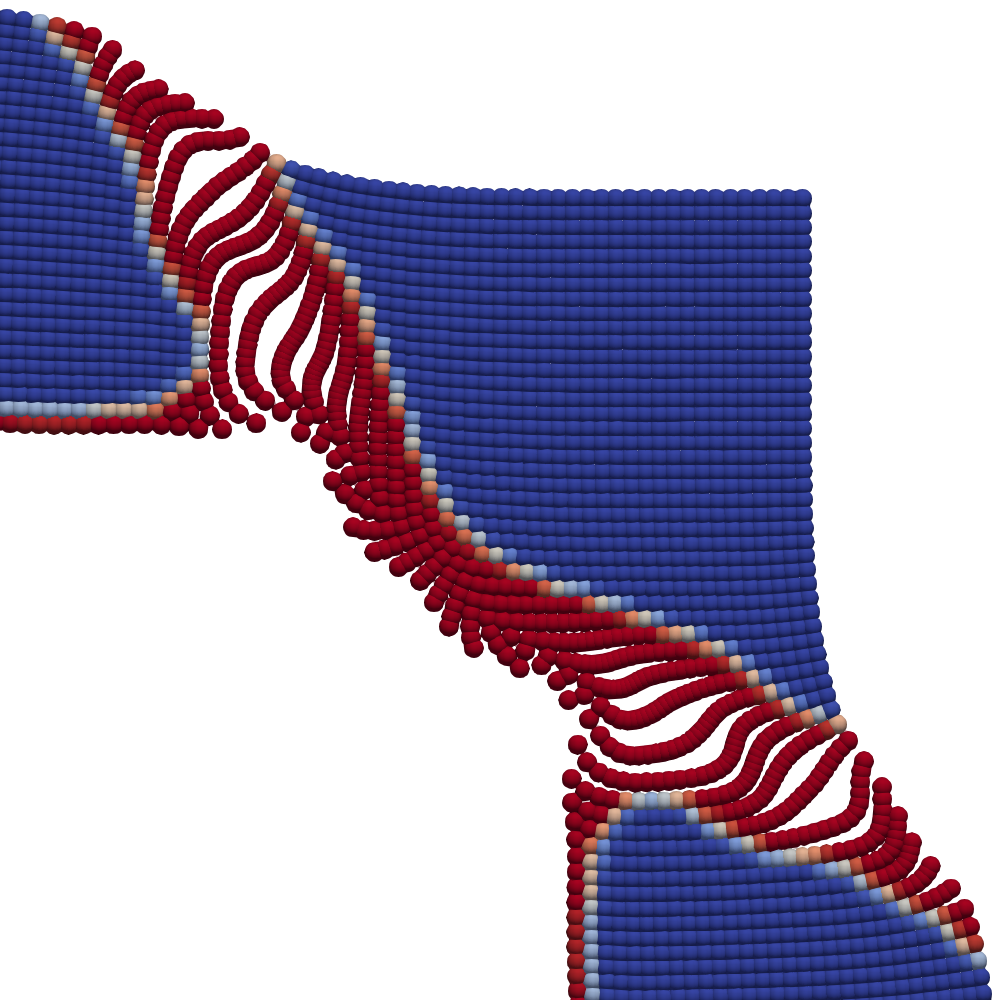}}
  \hspace{6pt}
  \subfloat[][]{\includegraphics[width=0.32\textwidth,height=0.31\textwidth,trim={0cm 0cm 0cm 0cm},clip]{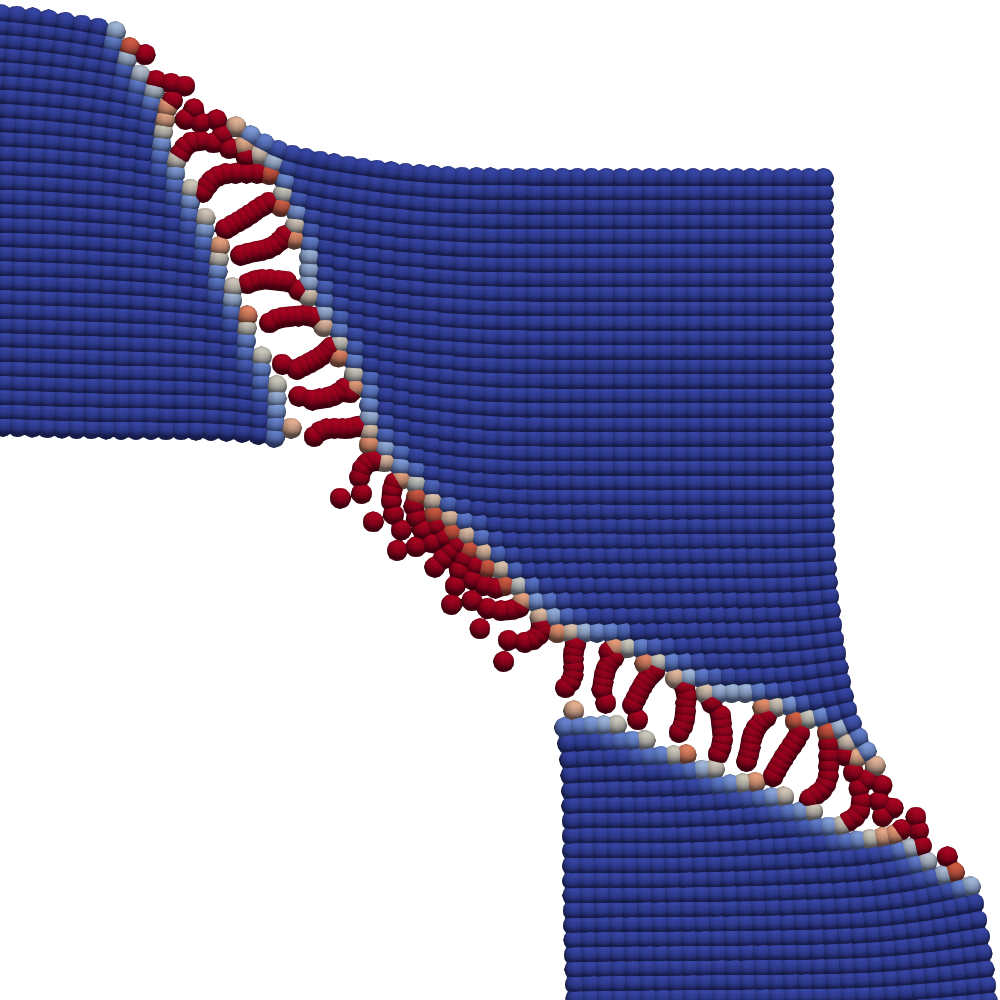}}
  \hspace{6pt}
  \subfloat[][]{\includegraphics[width=0.32\textwidth,height=0.31\textwidth,trim={0cm 0cm 0cm 0cm},clip]{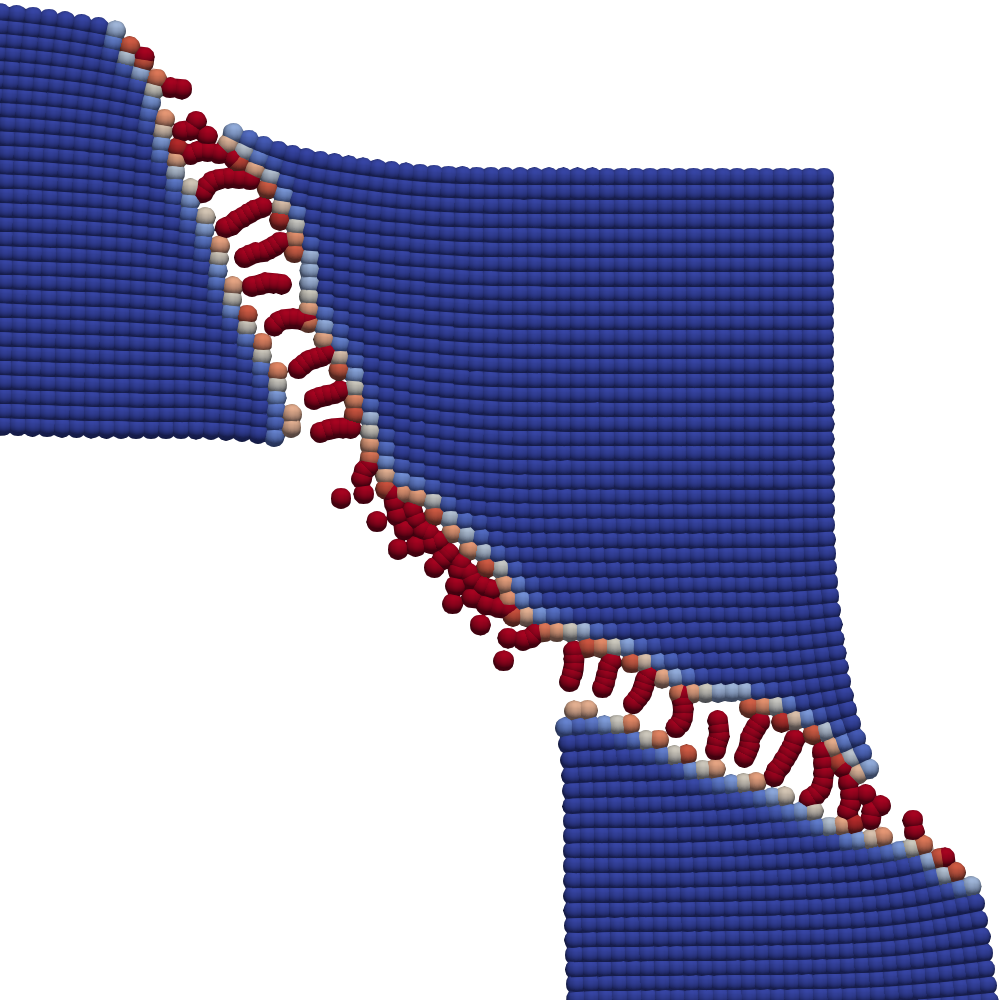}}
  \caption{Ductile fracture problem. Comparison of solid damage for different coupling approaches. Zoomed-in view of the solutions for the finest mesh. Top and bottom rows correspond to \SI{150}{\micro s} and \SI{200}{\micro s}, respectively. The number of fully damaged PD nodes at the fracture interface is indicated. (a): Strong coupling. (b): Weak coupling without damage in the penalty stiffness. (c) Weak coupling with damage in the penalty stiffness.}
  \label{fig:ductile_damage_zoomup}
\end{figure*}

\begin{figure*}[!hbpt]
  \centering
  \subfloat[][]{\includegraphics[width=0.345\textwidth,height=0.31\textwidth,trim={0cm 0cm 0cm 0cm},clip]{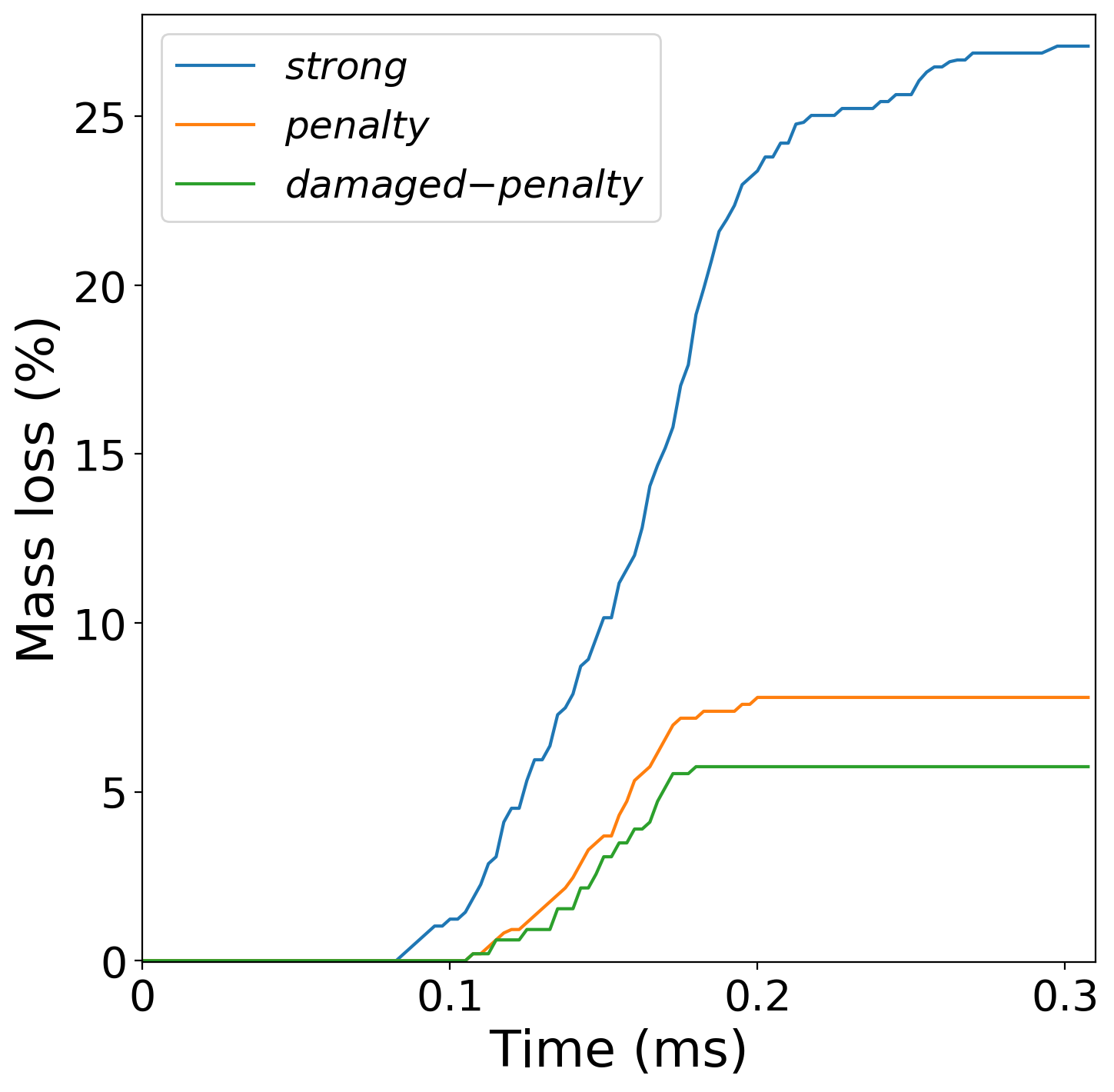}}
  \hspace{5pt}
  \subfloat[][]{\includegraphics[width=0.305\textwidth,height=0.31\textwidth,trim={2.15cm 0cm 0cm 0cm},clip]{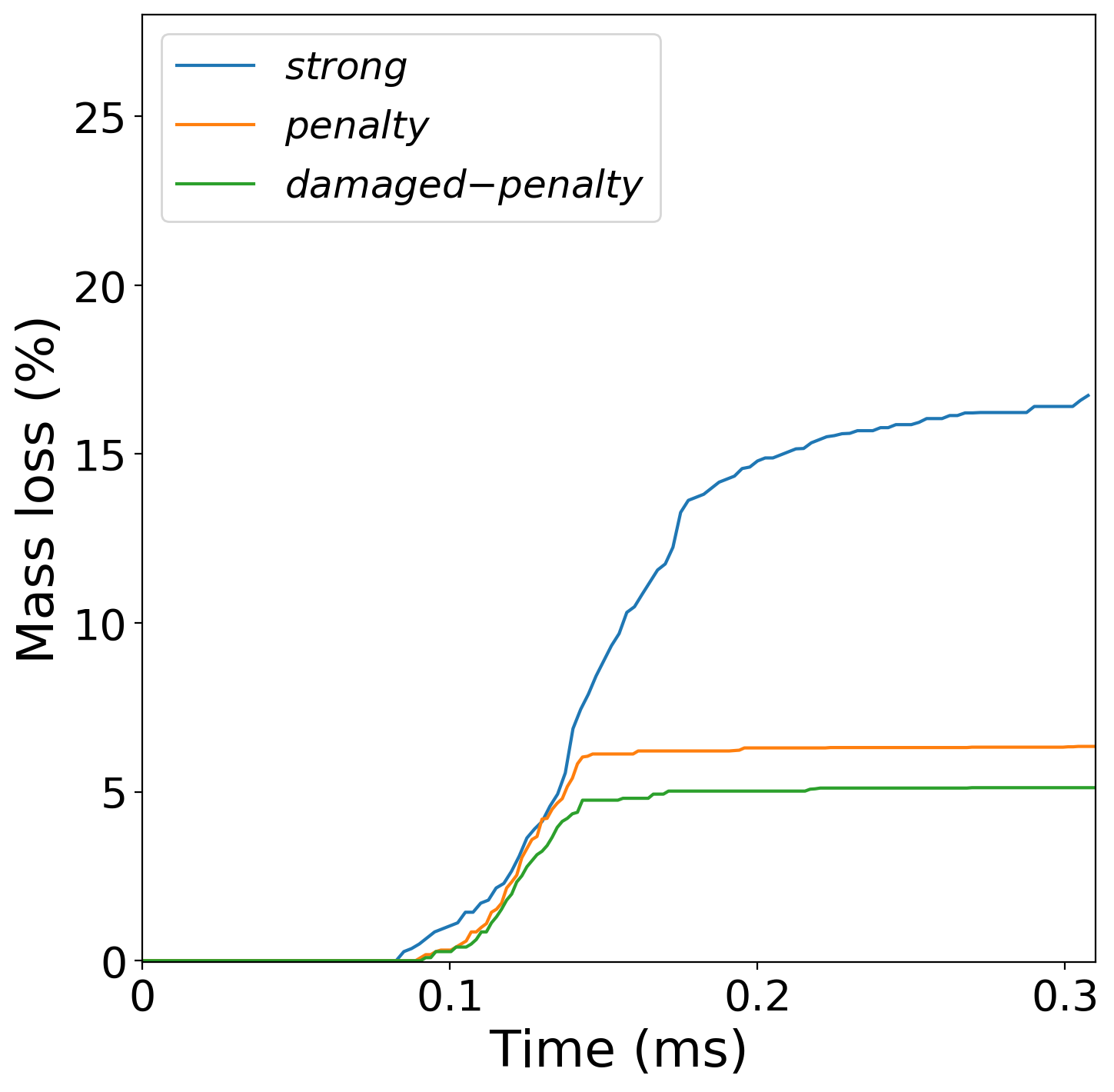}}
  \hspace{5pt}
  \subfloat[][]{\includegraphics[width=0.305\textwidth,height=0.31\textwidth,trim={2.15cm 0cm 0cm 0cm},clip]{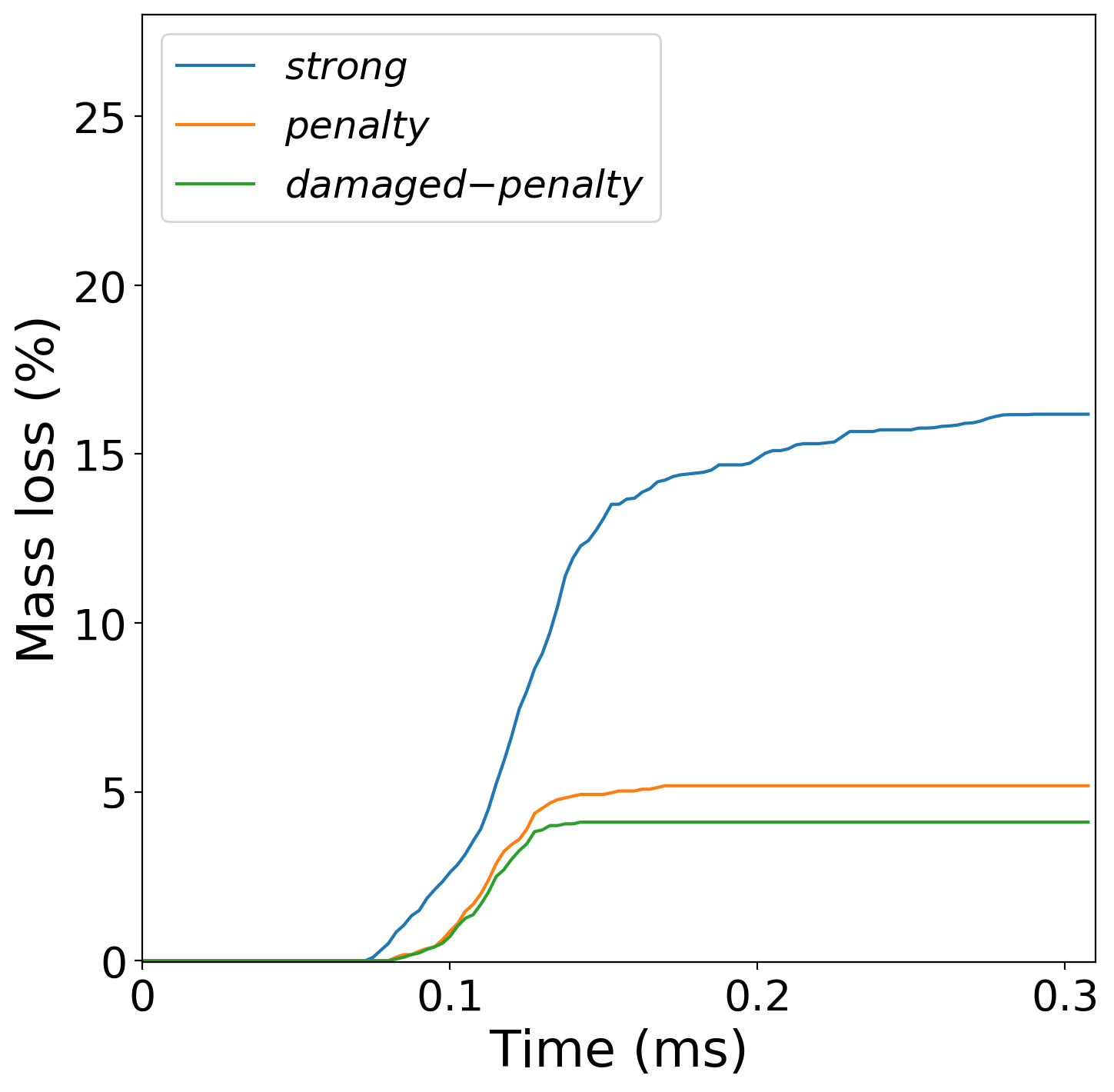}}
  \caption{Ductile fracture problem. Comparison of the normalized solid mass loss for different coupling approaches and discretizations. (a): Coarse mesh. (b): Medium mesh. (c): Fine mesh. Here, \textit{strong} indicates strong coupling, \textit{penalty} indicates Weak coupling without damage in the penalty stiffness, and \textit{damaged-penalty} indicates weak coupling with damage in the penalty stiffness.}
  \label{fig:ductile_mass_loss}
\end{figure*}

\cref{fig:ductile_contours} shows the results of the finest-mesh simulations, where the air speed is plotted on the background mesh while the damage field is plotted on the foreground PD nodes. In all cases the fractures initiate at the interior corners, the locations of stress concentration, and show a very similar final pattern. However, the strongly coupled case exhibits much thicker damage bands and the structural response is sufficiently different to alter the fluid behavior as predicted by the strongly and weakly coupled cases. For the weak coupling the damage bands are much narrower and, as a result, the fractures are much sharper. \cref{fig:ductile_damage_zoomup} examines the thickness of the damage band in more detail and compares the results of the strong coupling to the weak coupling, with and without the damage variable affecting the penalty stiffness. The results with the damage variable affecting the penalty stiffness correspond to the sharpest and cleanest fracture. The time history of the normalized solid mass loss is plotted in \cref{fig:ductile_mass_loss} and gives a quantitative confirmation of the observations in \cref{fig:ductile_contours} and \cref{fig:ductile_damage_zoomup}.

\subsection{Brittle Material Subjected to Internal Explosion}
\label{sec:brittle}

\begin{figure*}[!hbpt]
  \centering
  \subfloat[][]{\includegraphics[width=0.7\textwidth,trim={0cm 0cm 0cm 0cm},clip]{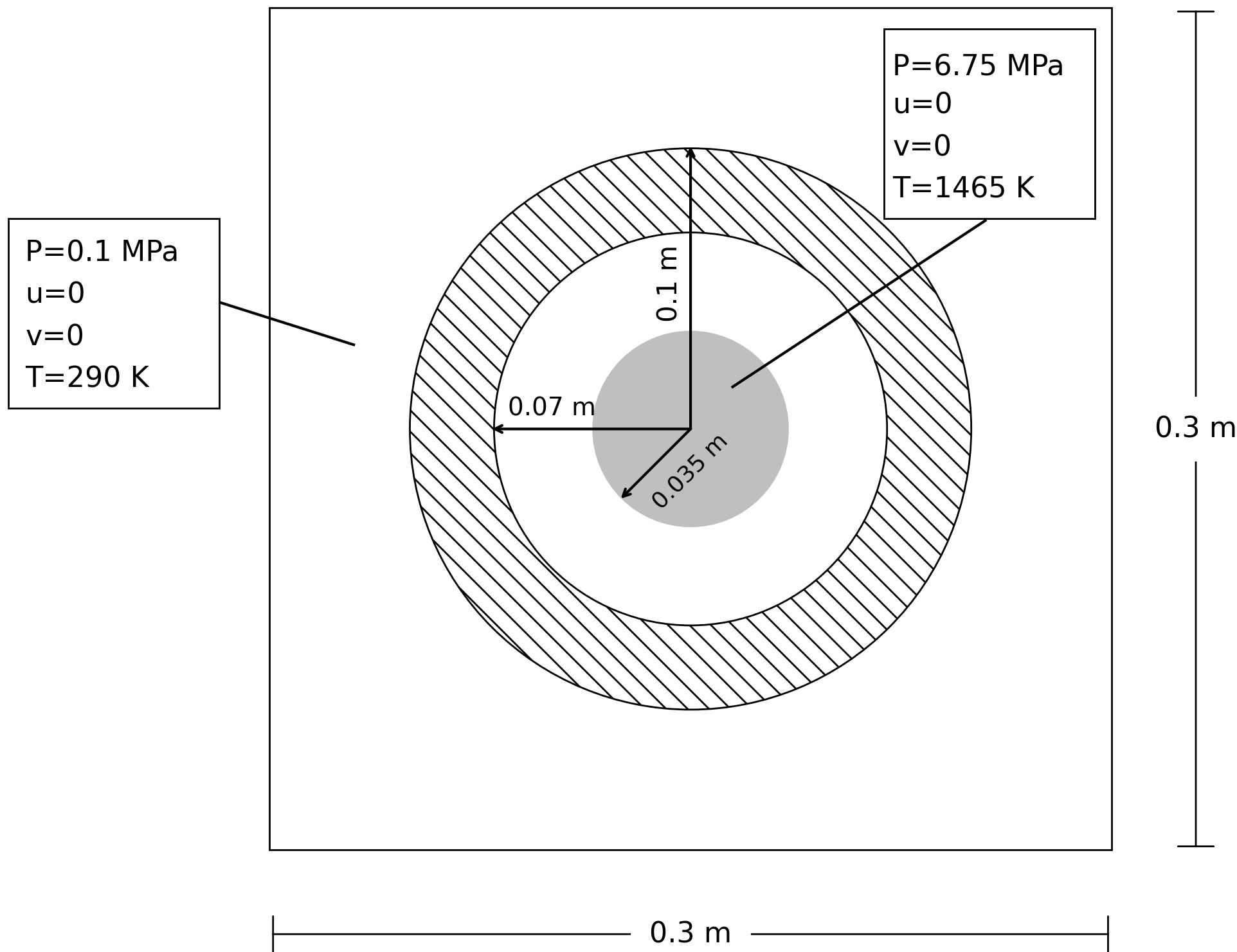}}
  \caption{Brittle fracture problem. A hollow cylinder made of elastic, brittle material subjected to internal explosion. Problem setup and geometry.}
  \label{fig:brittle_setup}
\end{figure*}

Here we consider a hollow cylinder made of elastic, brittle material subjected to internal detonation. As shown in \cref{fig:brittle_setup}, the detonation is initiated at the center of the hollow cylinder with the inner radius of $7 \, {\rm cm}$, outer radius of $10 \, {\rm cm}$, and thickness of $5 \, {\rm mm}$. The background domain is a $30 \, {\rm cm} \, \times \, 30 \, {\rm cm}$ square enclosing the cylinder domain. 
The solid elastic material has Young's modulus ${\rm E}=200 \, {\rm GPa}$, Poisson's ratio $\nu=0.29$, and initial density $\rho^s=7870 \, {\rm kg}/{\rm m}^3$. To simulate brittle fracture using the bond-associative damage model~\cite{behzadinasab2020semi,behzadinasab2020peridynamic}, the maximum principal stress failure criterion is used with the critical stress $\sigma_{\rm cr} = 3 \, {\rm GPa}$. A PD bond is broken once its associated maximum principal stress exceeds $\sigma_{\rm cr}$. Initially, the air is at rest with $p = 0.1 \, {\rm MPa}$ and $T = 290 \, {\rm K}$. The detonation condition is enforced by setting the initial pressure $p = 6.75 \, {\rm MPa}$ and temperature $T = 1465 \, {\rm K}$ in a circular area with a radius of $3.5 \, {\rm cm}$, centered inside the hollow cylinder. 

The background and foreground domains are discretized using a uniform rectangular mesh and a semi-uniform nodal spacing (uniform along the $\theta$-direction), respectively. We again consider three meshes, with the solid node spacing of $h = 2 \, {\rm mm}$, $1.5 \, {\rm mm}$, and $1 \, {\rm mm}$. The fluid mesh size is set to three times that of the solid node spacing in each case. The time step size for the coarse, medium, and fine mesh in the strongly coupled case is set to $\SI{0.4}{\micro s}$, $\SI{0.3}{\micro s}$, and $\SI{0.2}{\micro s}$, respectively. The time step size used for the weakly coupled cases is four times smaller on the corresponding meshes. The penalty constant is set to $\beta = 1$.

\begin{figure*}[!hbpt]
  \centering
  \subfloat{\includegraphics[width=0.6\textwidth,trim={0cm 0cm 0cm 4cm},clip]{air_vel_scale.png}}
  
  \setcounter{subfigure}{0}
  \subfloat[][\SI{50}{\micro s}]{\includegraphics[width=0.235\textwidth,trim={0cm 0cm 0cm 0cm},clip]{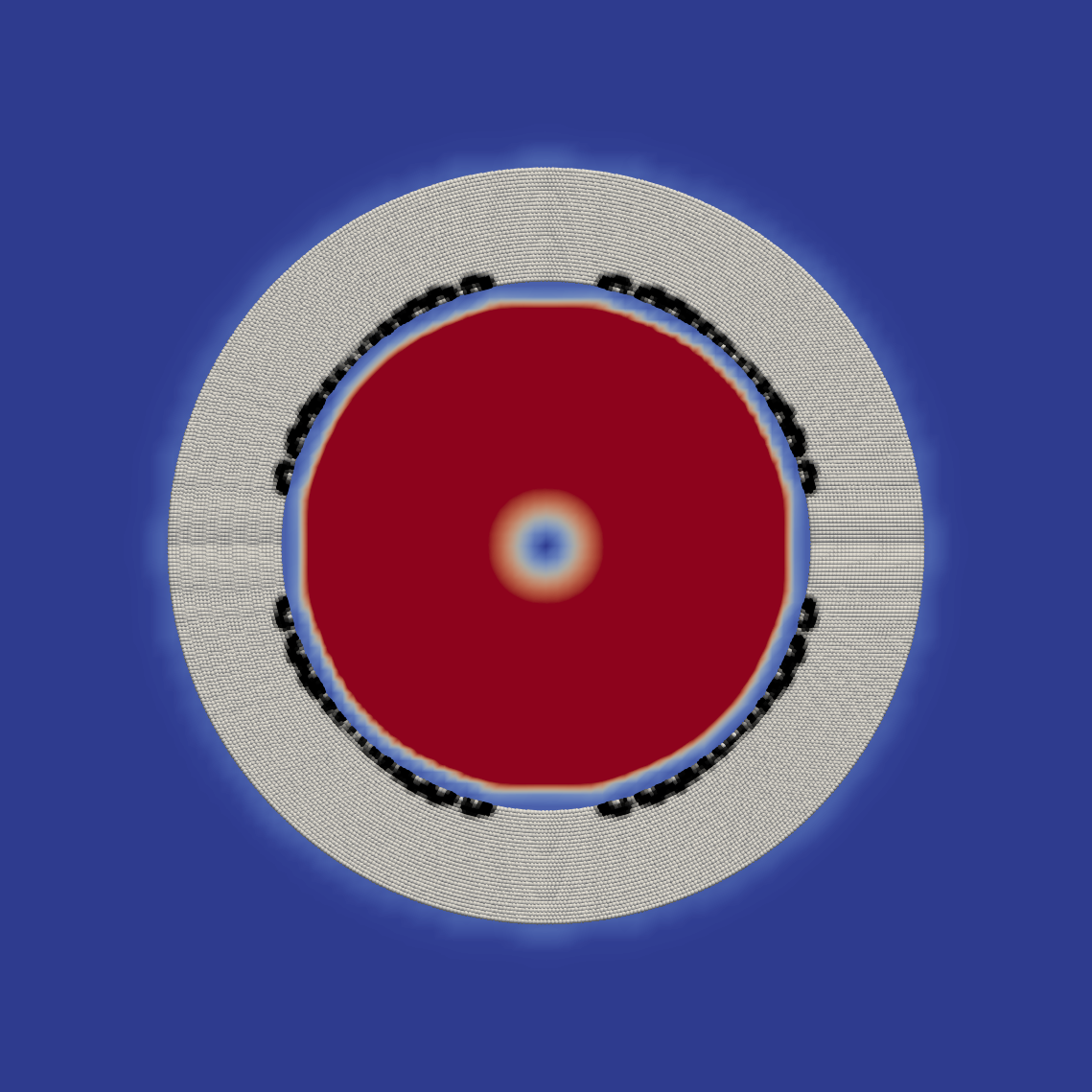}}
  \hspace{5pt}
  \subfloat[][\SI{67}{\micro s}]{\includegraphics[width=0.235\textwidth,trim={0cm 0cm 0cm 0cm},clip]{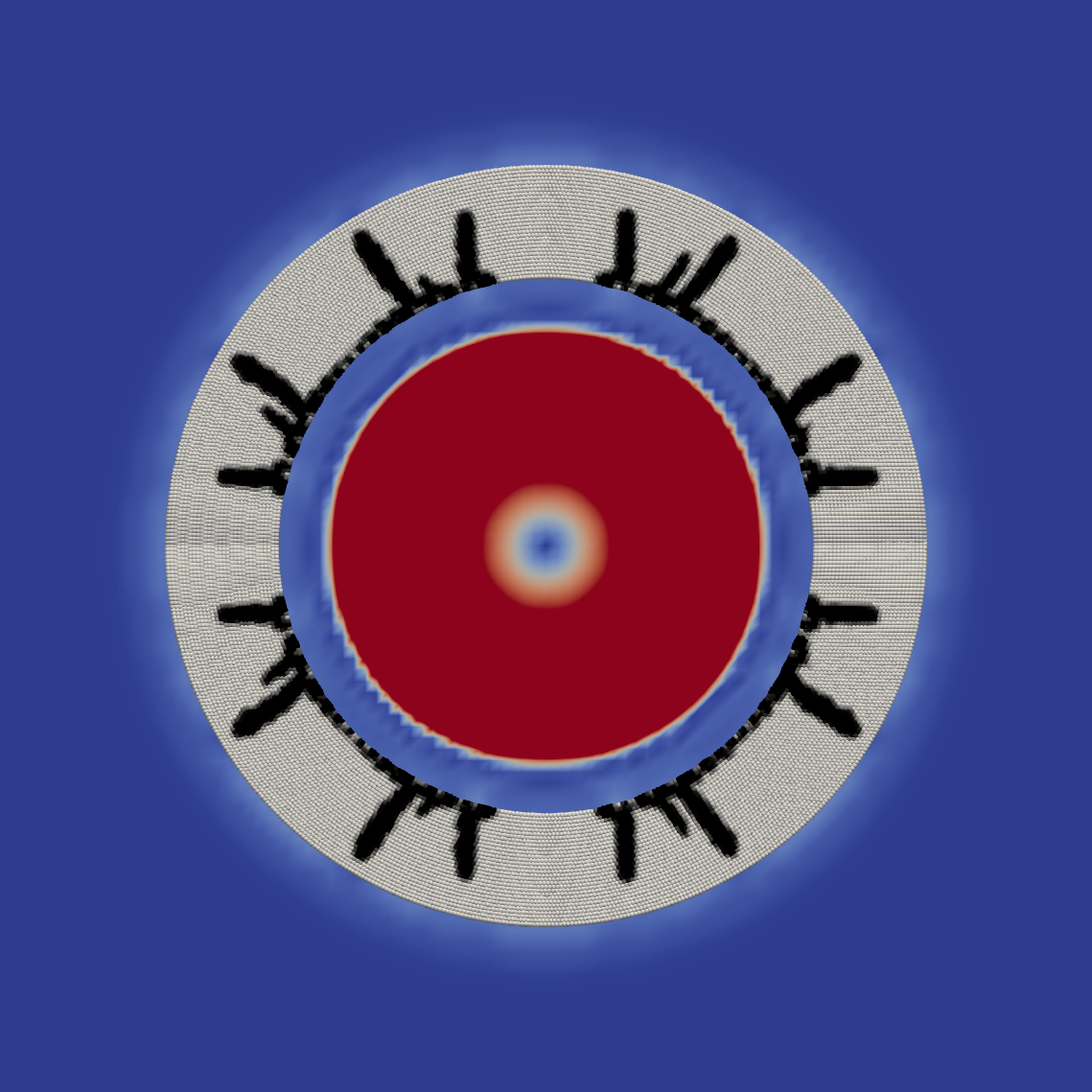}}
  \hspace{5pt}
  \subfloat[][\SI{150}{\micro s}]{\includegraphics[width=0.235\textwidth,trim={0cm 0cm 0cm 0cm},clip]{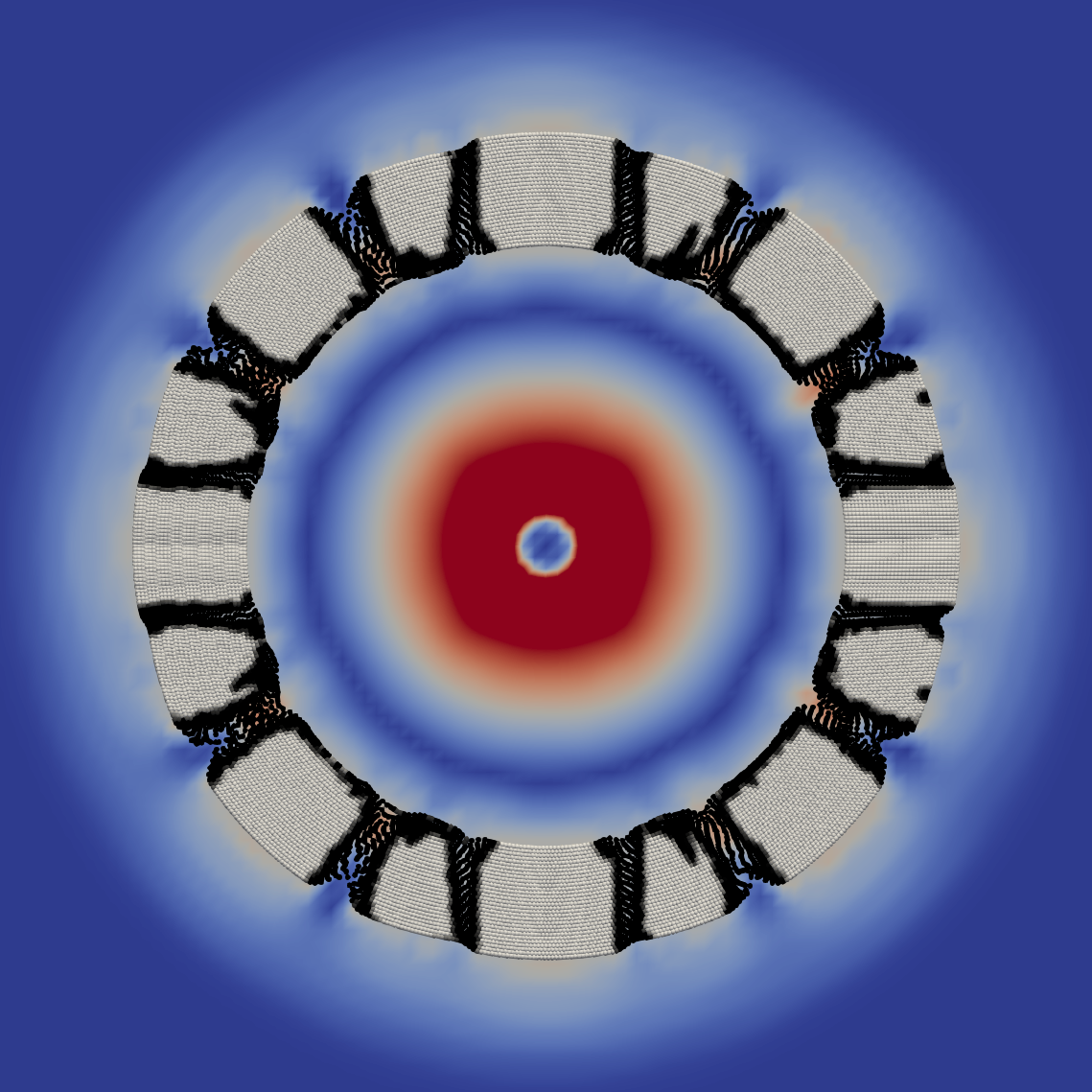}}
  \hspace{5pt}
  \subfloat[][\SI{300}{\micro s}]{\includegraphics[width=0.235\textwidth,trim={0cm 0cm 0cm 0cm},clip]{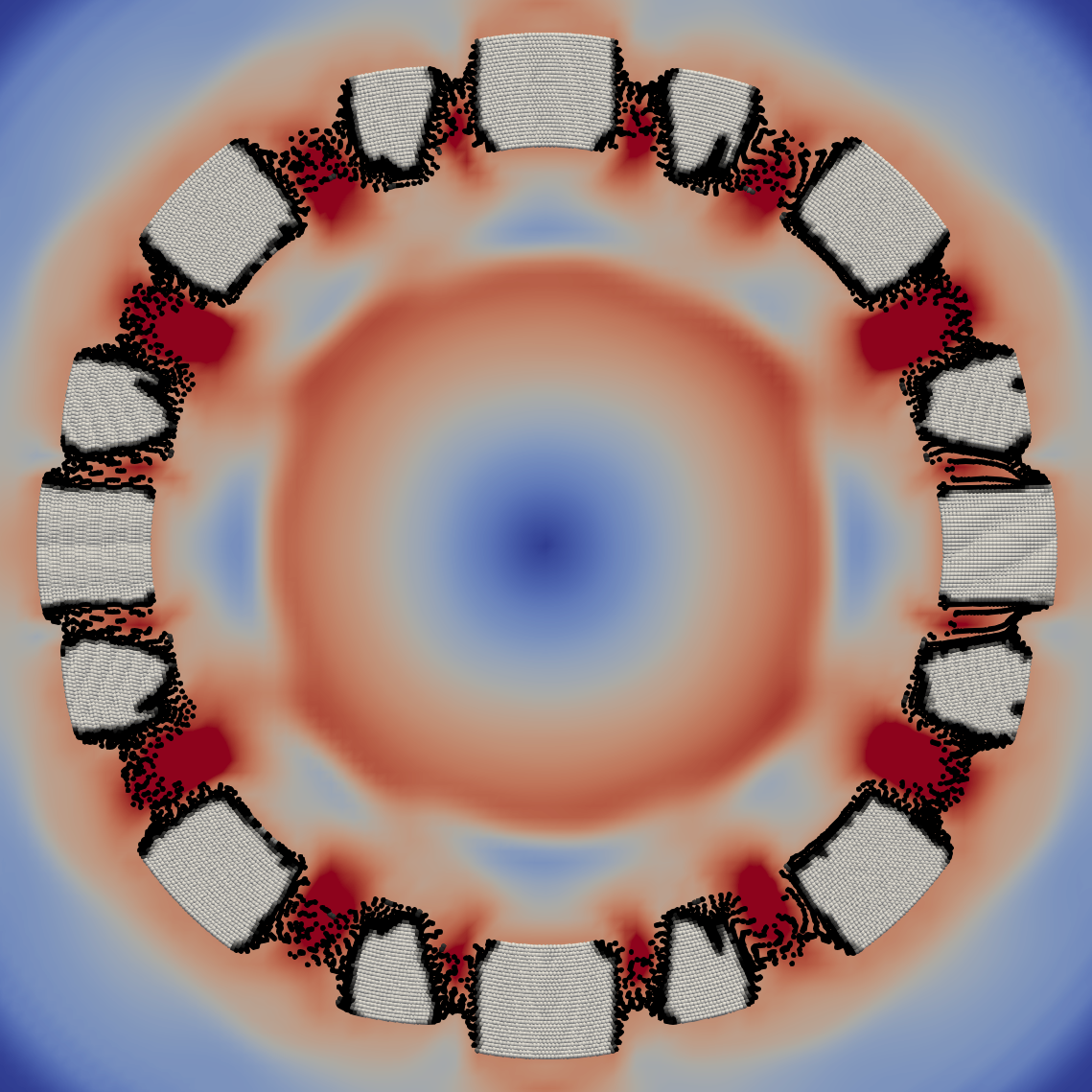}}
  
  \subfloat[][\SI{67}{\micro s}]{\includegraphics[width=0.235\textwidth,trim={0cm 0cm 0cm 0cm},clip]{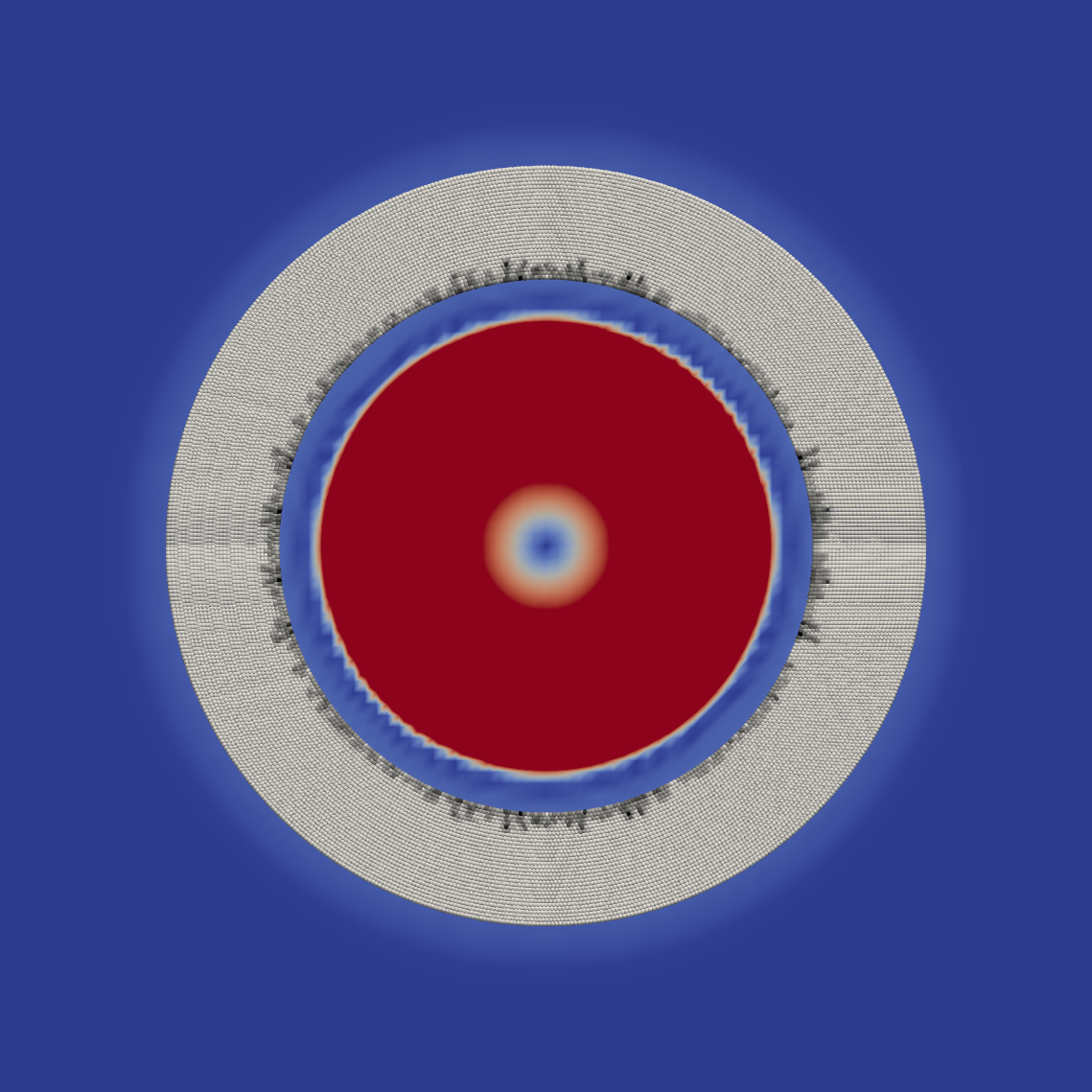}}
  \hspace{5pt}
  \subfloat[][\SI{75}{\micro s}]{\includegraphics[width=0.235\textwidth,trim={0cm 0cm 0cm 0cm},clip]{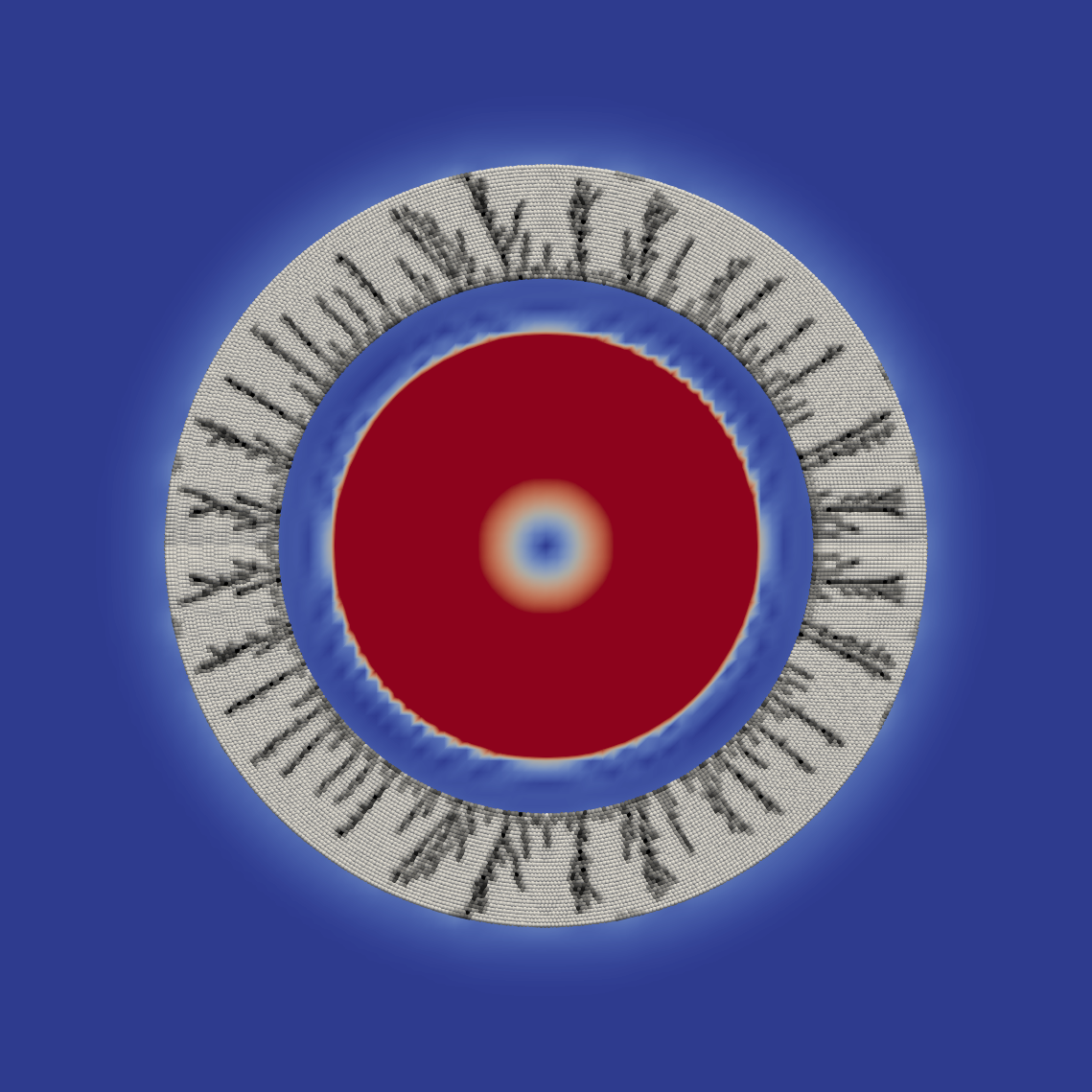}}
  \hspace{5pt}
  \subfloat[][\SI{150}{\micro s}]{\includegraphics[width=0.235\textwidth,trim={0cm 0cm 0cm 0cm},clip]{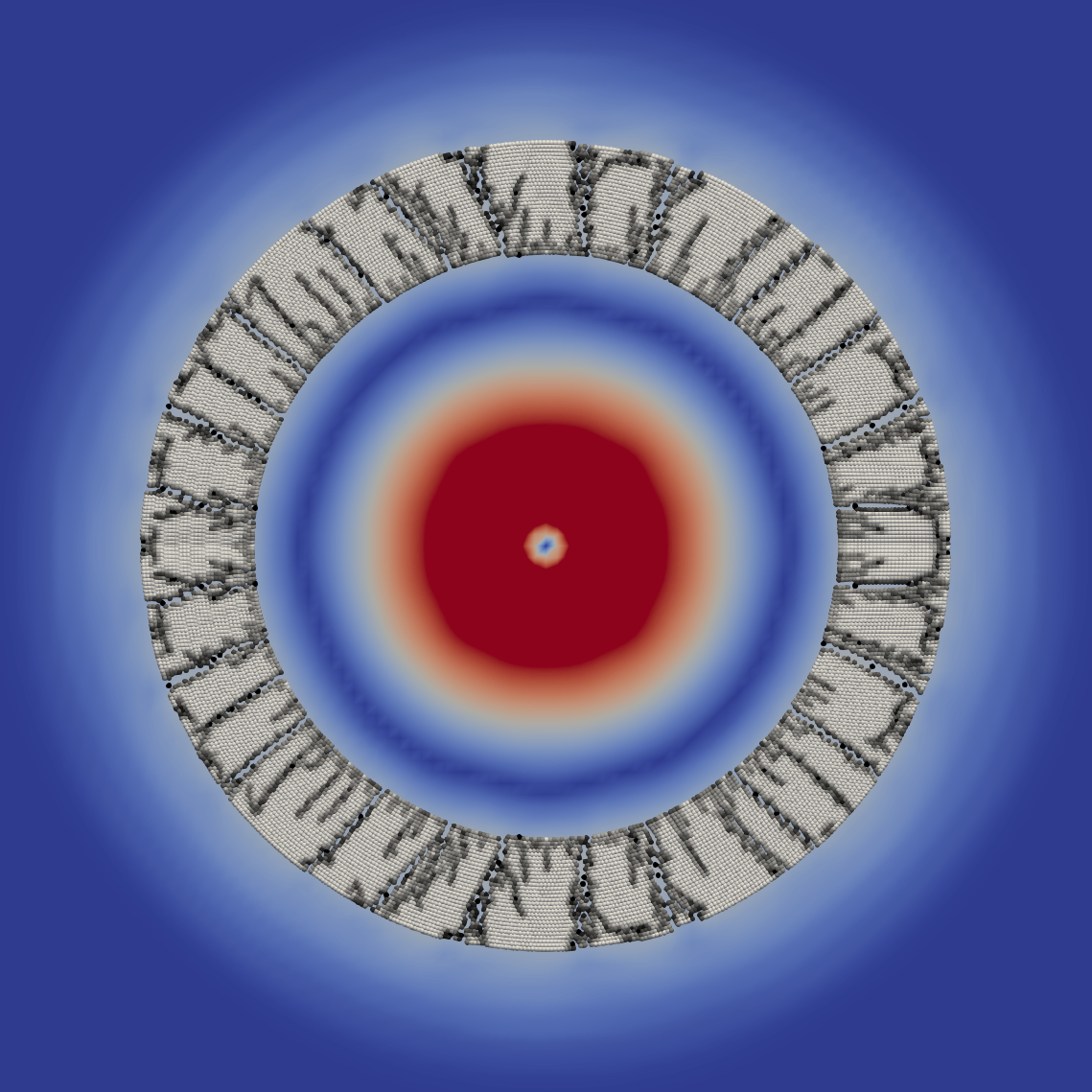}}
  \hspace{5pt}
  \subfloat[][\SI{300}{\micro s}]{\includegraphics[width=0.235\textwidth,trim={0cm 0cm 0cm 0cm},clip]{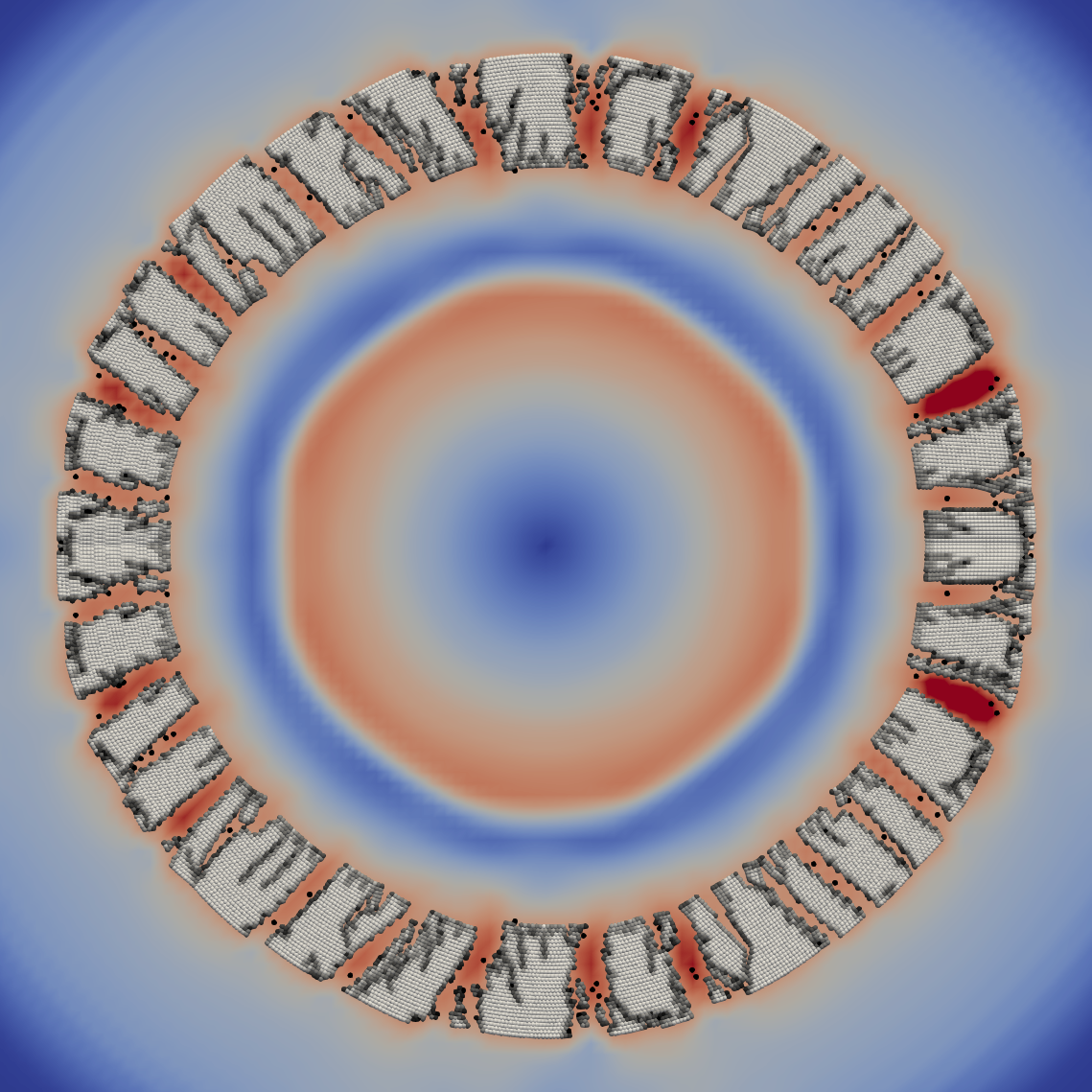}}

  \subfloat[][\SI{67}{\micro s}]{\includegraphics[width=0.235\textwidth,trim={0cm 0cm 0cm 0cm},clip]{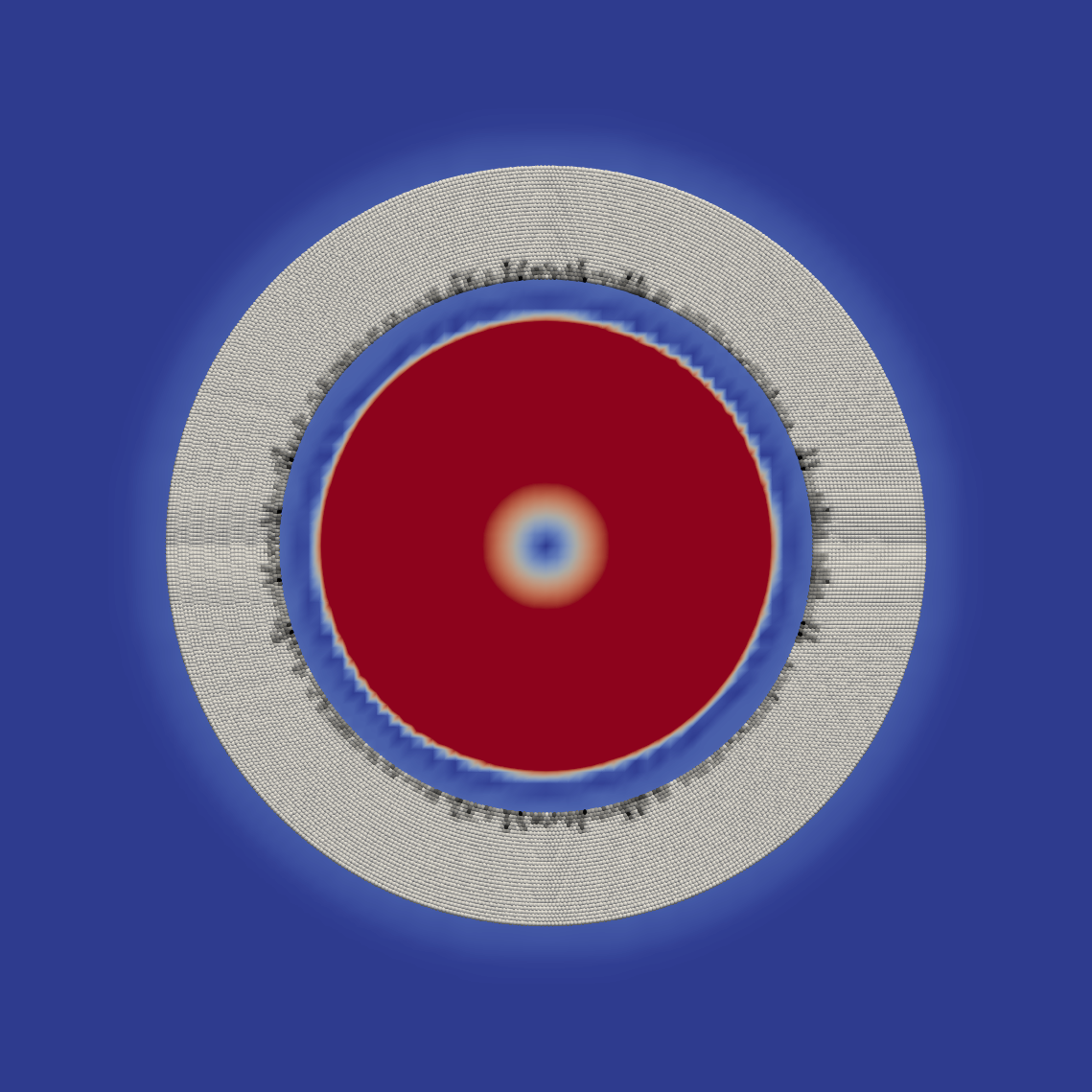}}
  \hspace{5pt}
  \subfloat[][\SI{75}{\micro s}]{\includegraphics[width=0.235\textwidth,trim={0cm 0cm 0cm 0cm},clip]{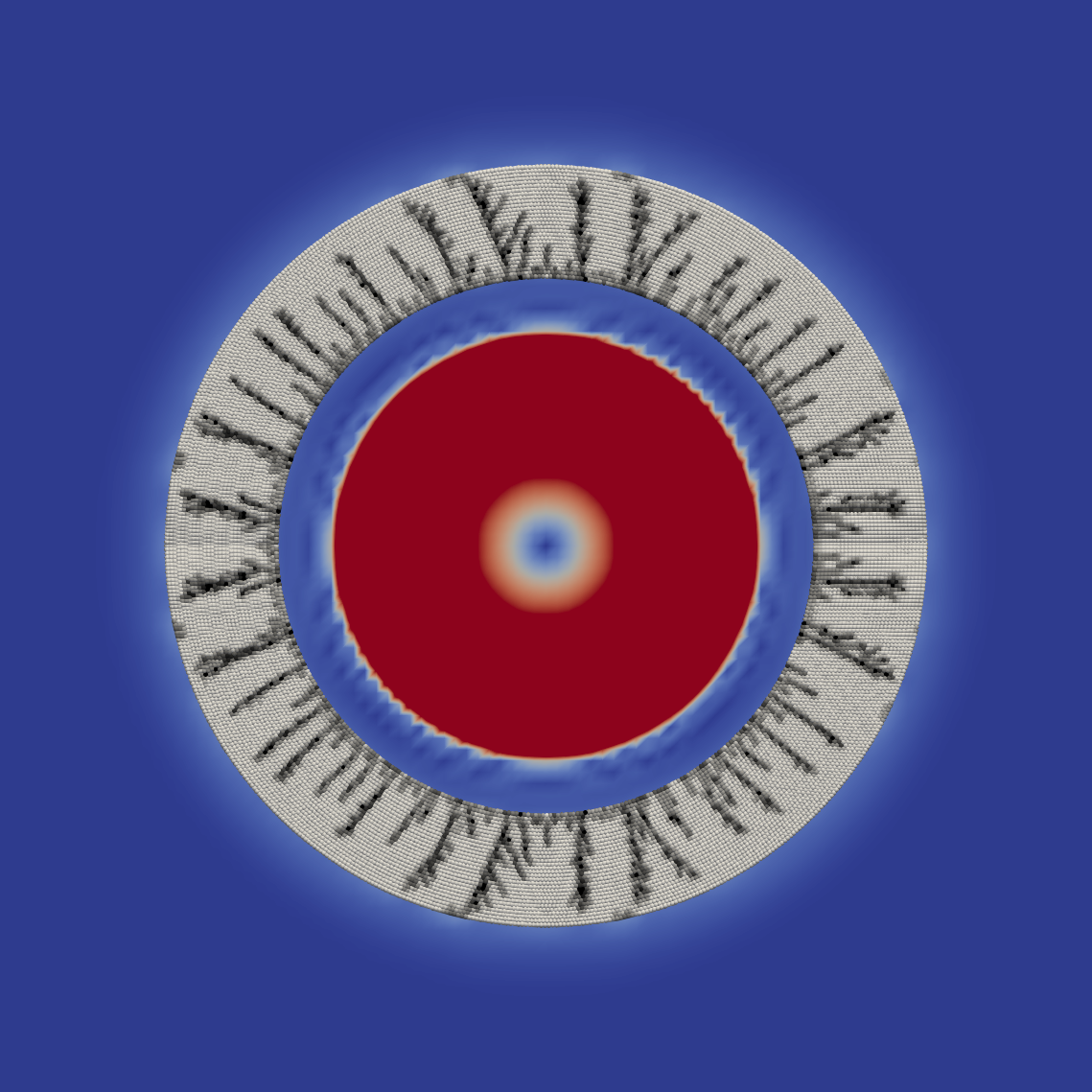}}
  \hspace{5pt}
  \subfloat[][\SI{150}{\micro s}]{\includegraphics[width=0.235\textwidth,trim={0cm 0cm 0cm 0cm},clip]{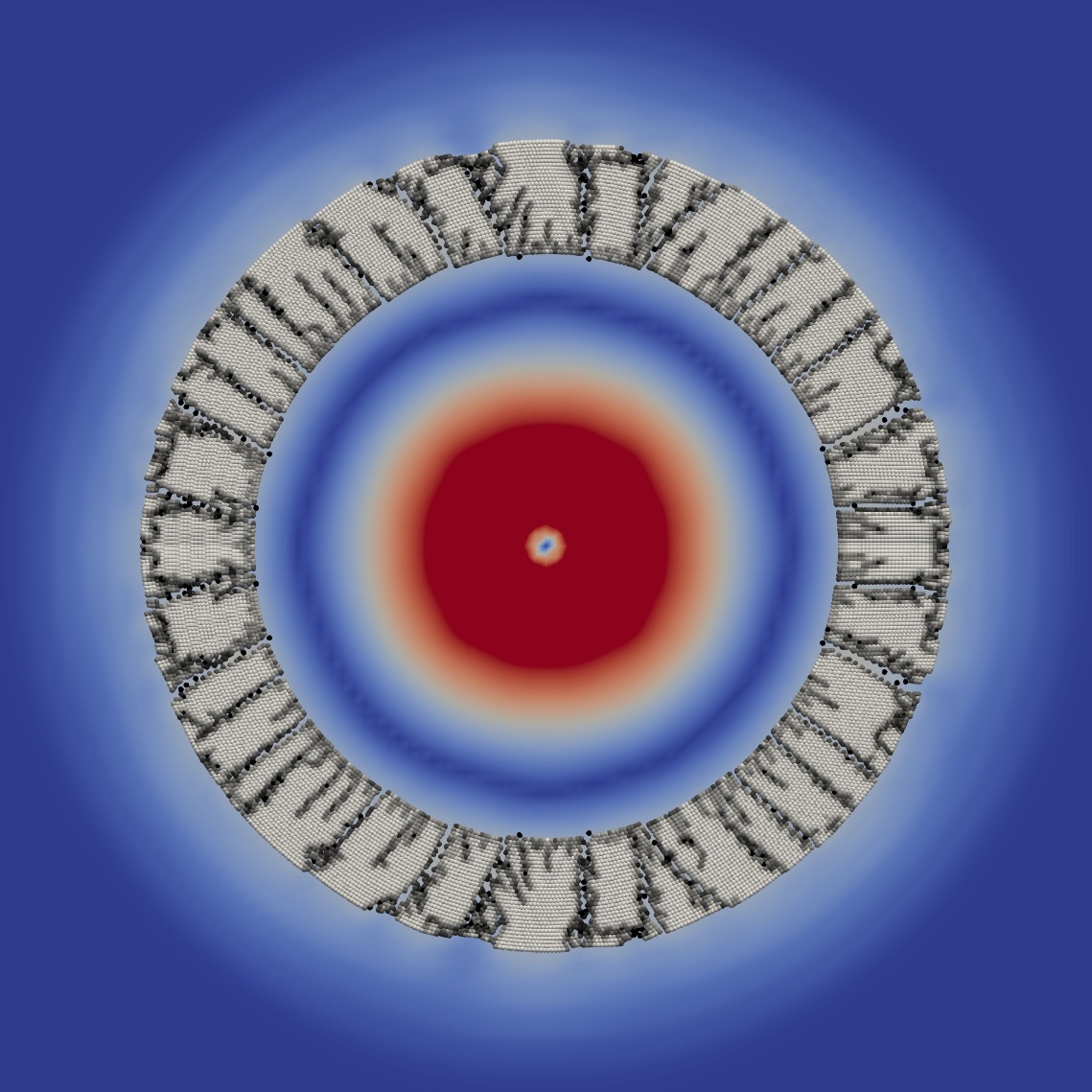}}
  \hspace{5pt}
  \subfloat[][\SI{300}{\micro s}]{\includegraphics[width=0.235\textwidth,trim={0cm 0cm 0cm 0cm},clip]{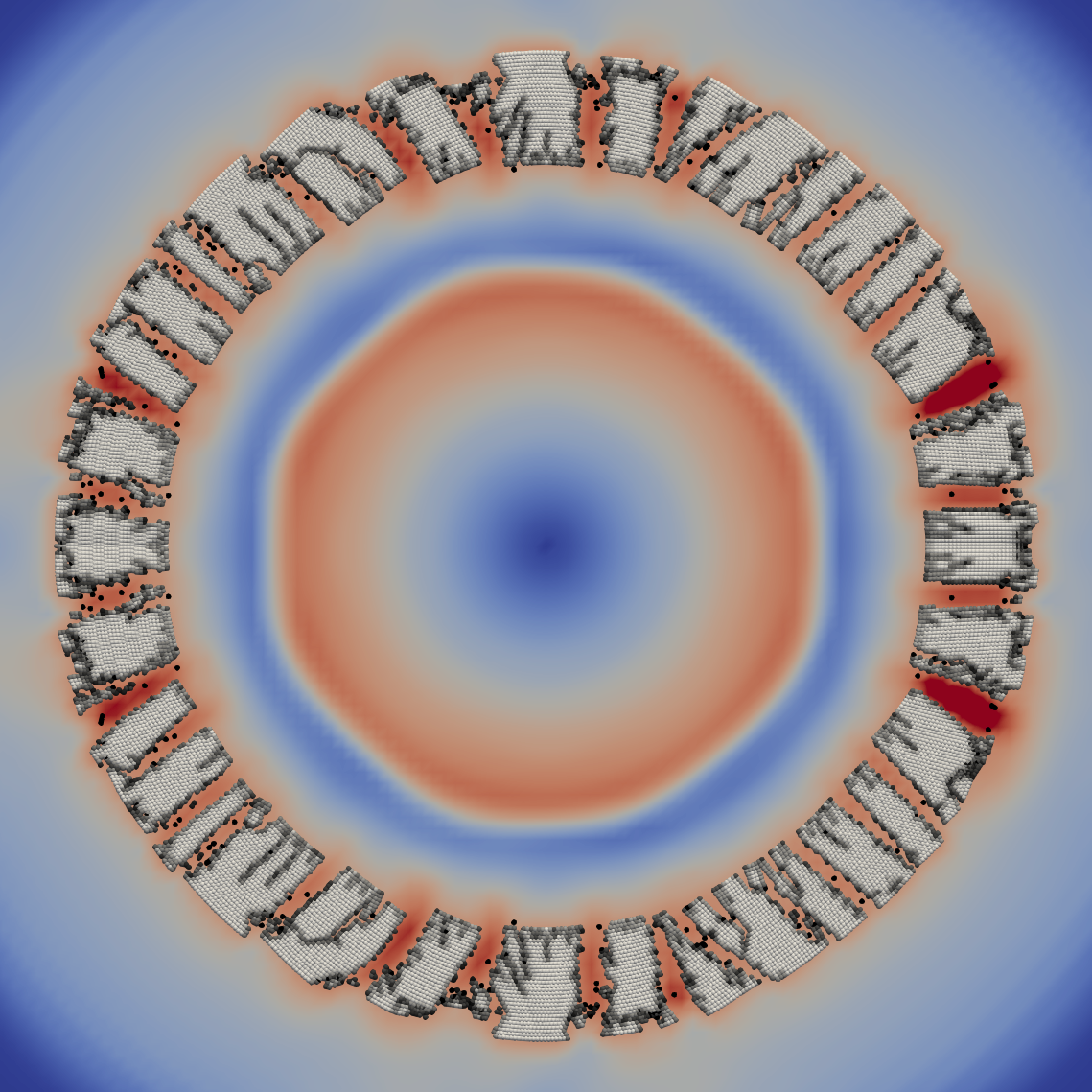}}

  \caption{Brittle fracture problem. Snapshots of air speed (in m/s) and solid damage in the current configuration computed on the finest mesh. Top, middle, and bottom rows correspond to strong coupling, weak coupling without damage in the penalty stiffness, and weak coupling with damage in the penalty stiffness, respectively.}
  \label{fig:brittle_contours}
\end{figure*}

\begin{figure*}[!hbpt]
  \centering
  \subfloat[][]{\includegraphics[width=0.345\textwidth,height=0.31\textwidth,trim={0cm 0cm 0cm 0cm},clip]{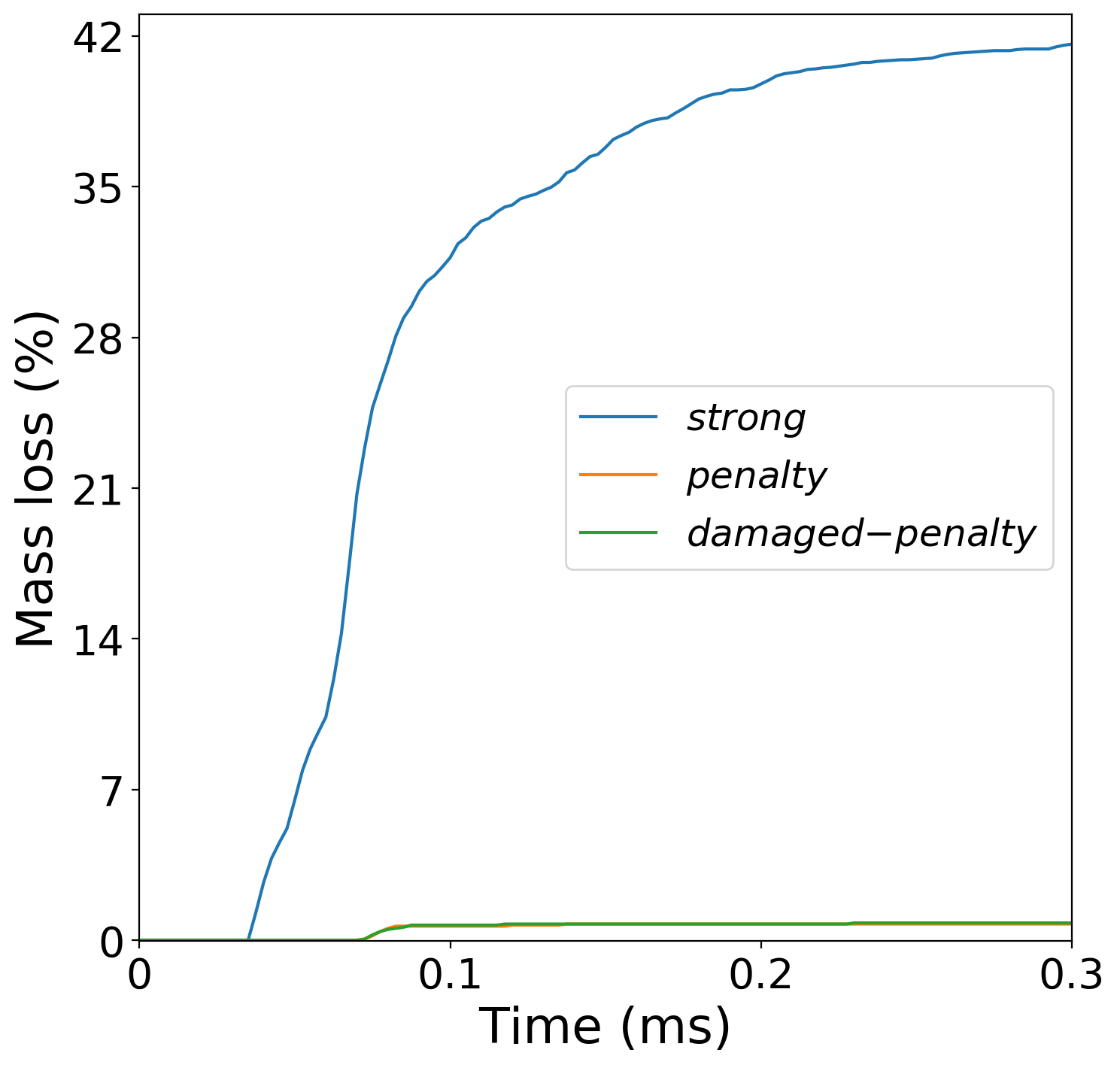}}
  \hspace{5pt}
  \subfloat[][]{\includegraphics[width=0.305\textwidth,height=0.31\textwidth,trim={2.15cm 0cm 0cm 0cm},clip]{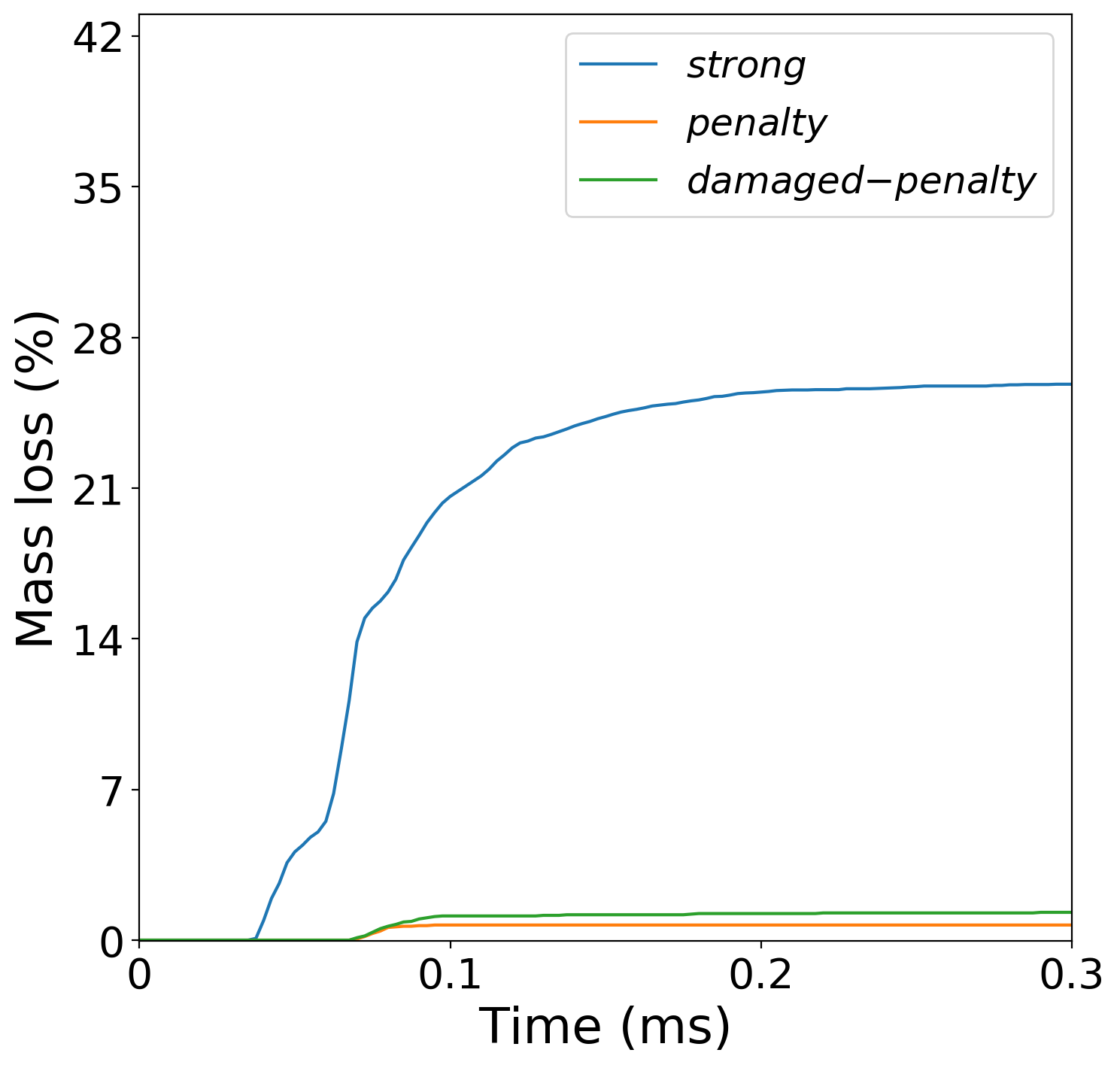}}
  \hspace{5pt}
  \subfloat[][]{\includegraphics[width=0.305\textwidth,height=0.31\textwidth,trim={2.15cm 0cm 0cm 0cm},clip]{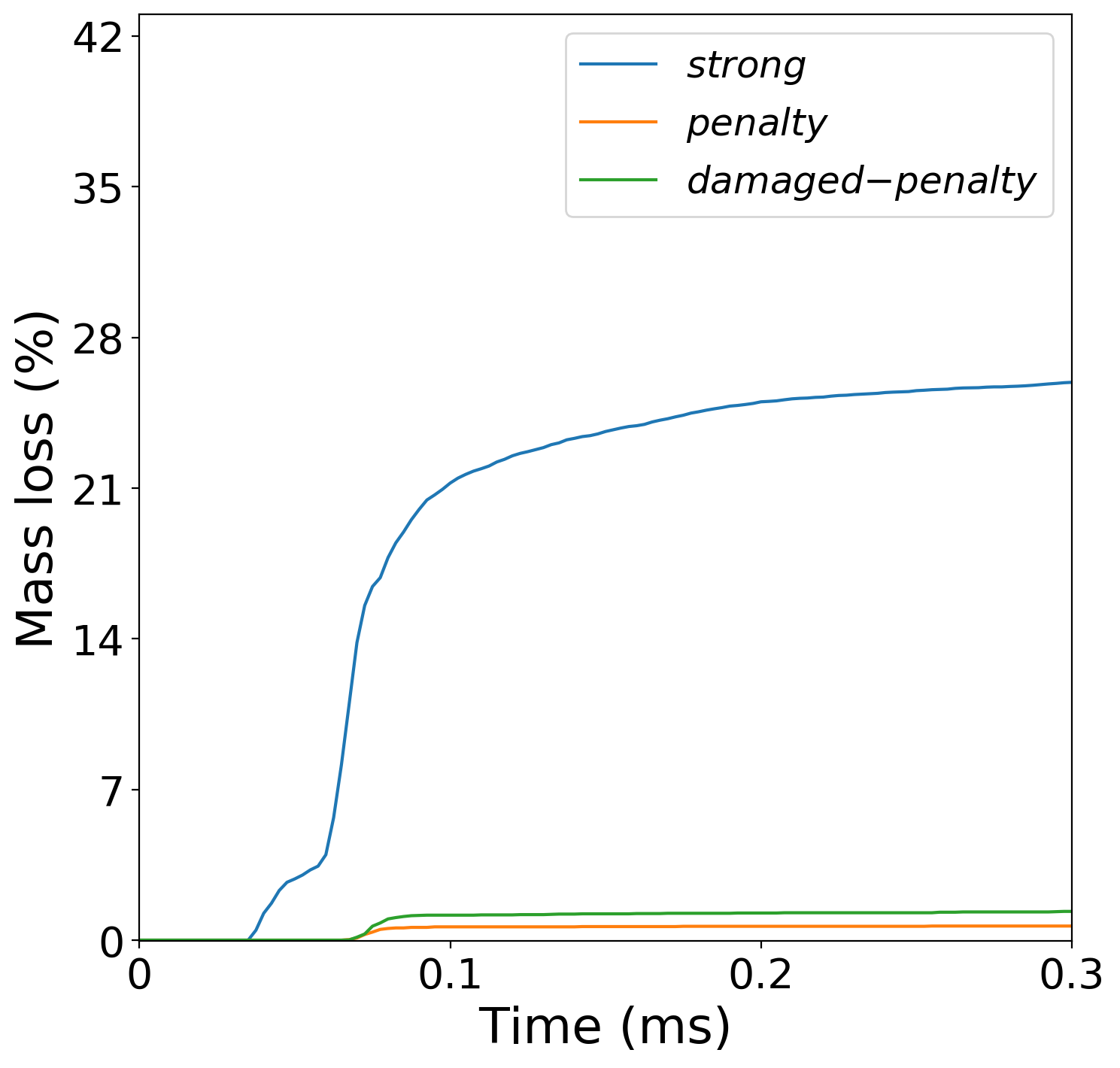}}
  \caption{Brittle fracture problem. Normalized solid mass loss for different coupling approaches and discretizations. (a): Coarse mesh. (b): Medium mesh. (c): Fine mesh. Here, \textit{strong} indicates strong coupling, \textit{penalty} indicates Weak coupling without damage in the penalty stiffness, and \textit{damaged-penalty} indicates weak coupling with damage in the penalty stiffness.}
  \label{fig:brittle_mass_loss}
\end{figure*}

\cref{fig:brittle_contours} shows the air speed on the background grid and the damage field on the PD nodes in the current configuration. Unlike in the ductile case, the blast wave shatters the brittle material into many small fragments as predicted by the weakly coupled simulations. However, in the in the strongly coupled case, the background mesh is not able to support such fine fragments and produces a much smaller number of larger-size chunks. This feature of the weakly coupled methods to enable the solid to fragment into small chunks that are not constrained in size to the resolution of the background grid is remarkable and presents a real breakthrough for the immersed FSI methods. We also note that the fragmentation results are very similar for the cases with and without damage dependence in the definition of the penalty parameter. The time history of the normalized solid mass loss reported in \cref{fig:brittle_mass_loss} shows a dramatic difference between the strongly and weakly coupled cases.

\section{Conclusions}
\label{sec:conclusions}

We developed a practical computational framework that is capable of capturing the mechanics of air blast coupled to solids and structures that undergo large, inelastic deformations, damage and fragmentation. The foundation for the proposed framework is an immersed FSI approach, which does not require explicit tracking of the fluid-structure interfaces and which has no limitations on the solid domain motion and topology. Weak forms of the fluid and structural mechanics equations are discretized on the background and foreground domains, respectively, and are coupled by means of a volumetric penalty operator, which is the main novelty of the proposed approach. We employ IGA based on NURBS in the background domain and a correspondence-based PD solid in the foreground domain using the RK functions to define the nonlocal derivatives~\cite{behzadinasab2021unifiedI,behzadinasab2021unifiedII}. We feel that the combination of these numerical methodologies is particularly attractive for the problem class of interest due to the higher-order accuracy and smoothness of IGA and RK, the benefits of using immersed methodology in handling the fluid--structure interfaces and coupling, and the unique capabilities of PD for modeling fracture and fragmentation. 

Using three numerical examples, the present work illustrates very clearly that strong coupling, a hallmark of Immersed Boundary Methods and Immersed Finite Element Methods, while well suited for FSI with large solid deformations, is not an optimal approach for the modeling of fracture and fragmentation. On the other hand, weak coupling remains accurate for large-deformation FSI and enables the modeling of fracture and fragmentation with a lot less sticking, and with fragment sizes that are not constrained to the resolution of the background mesh. As such, the proposed methodology presents a real breakthrough in the application of immersed methods to FSI with fracture and fragmentation. 

A likely explanation for the observed results is the fact that the volumetric penalty term attempts to minimize the error between the fluid and structural kinematics in the $L_2$-norm over the volume. $L_2$ is a weak norm that does not significantly penalize discontinuities or sharp gradients in the difference between the foreground and background solutions. As a result, the foreground solution is able to develop the discontinuities (i.e., fractures) that remain essentially undetected by the background mesh through the volumetric penalty operator. Conversely, in the strong coupling approach, any discontinuity generated on the foreground mesh is ``overwritten'' by the strong kinematics constraint to the smooth background grid that does not support solution discontinuities.

In the future, it may be worthwhile to investigate more elaborate volumetric Nitsche coupling approaches to reduce the dependence of the overall method performance, including the size of the stable time step in explicit simulations, on the choice of the penalty constant. However, one must be careful to not introduce stronger than necessary coupling and lose the benefits of the present approach.

\section*{Acknowledgments}
\label{sec:acknowledge}
Y.~Bazilevs and M.~Behzadinasab were supported through the ONR Grant No. N00014-21-1-2670. M.~Hillman was supported through the National Science Foundation award number 1826221.

\bibliographystyle{unsrt}
\bibliography{main}

\end{document}